\documentclass[11pt]{article}

\usepackage[margin=1in]{geometry}

\usepackage{amsmath, amsfonts, amsthm, amssymb, commath}

\usepackage{enumitem}

\usepackage{prettyref}

\usepackage[colorlinks,citecolor=blue, linkcolor=blue, urlcolor=blue]{hyperref}

\usepackage{mathrsfs}
\usepackage{natbib}
\usepackage{fancyhdr}
\usepackage{fullpage}
\usepackage{xcolor, framed}
\usepackage{mathtools}
\usepackage{mleftright}

\usepackage{etoolbox}
\usepackage{xcolor}
\usepackage[makeroom]{cancel}
\usepackage{booktabs}
\usepackage{algorithm}
\usepackage{algpseudocode}

\usepackage{esvect}

\newcommand{\Xt}{\widetilde{X}}
\newcommand{\Xto}{\widetilde{X}_o}
\newcommand{\Xho}{\widehat{X}_o}

\newcommand{\fl}[1]{\lfloor#1\rfloor}

\newcommand{\pr}[1]{\left(#1\right)}
\newcommand{\br}[1]{\left[#1\right]}
\newcommand{\cb}[1]{\left\{#1\right \}}
\newcommand{\nm}[1]{\left\|#1\right \|}
\newcommand{\nmhs}[1]{\left\|#1\right \|_{\rm HS}}
\renewcommand{\abs}[1]{\left|#1\right|}

\newcommand{\trbr}[1]{\operatorname{Tr}\left[#1\right]}

\newcommand{\xmin}{x_{\min}}
\newcommand{\xmax}{x_{\max}}
\newcommand{\smin}{\sigma_{\min}}
\newcommand{\smax}{\sigma_{\max}}
\newcommand{\lmin}{\lambda_{\min}}
\newcommand{\lmax}{\lambda_{\max}}

\newcommand{\la}{\lambda}

\newcommand{\e}{\varepsilon}
\newcommand{\de}{\delta}

\newcommand{\R}{\mathbb{R}}

\newcommand{\Q}{\mathbb{Q}}
\newcommand{\N}{\mathbb{N}} 

\newcommand{\Z}{\mathbb{Z}}

\newcommand{\Matrices}[1]{ \mathbb{M}\left( {#1} \right) } 
\newcommand{\SAMatrices}[1]{ \mathbb{M}_{SA}\left( {#1} \right) } 
\DeclareDocumentCommand{\Matrices}{ O{r} O{n} }{ \mathbb{M}_{#2} \left( {#1} \right) }
\DeclareDocumentCommand{\SAMatrices}{ O{r} O{n} }{ \mathbb{M}_{#2,SA} \left( {#1} \right) }

\newcommand{\diag}{\operatorname{diag}}

\newcommand{\E}{\mathbb{E}}
\renewcommand{\P}{\mathbb{P}}
\newcommand{\prob}{\mathbb{P}}

\newcommand{\Prob}{\prob}

\DeclareDocumentCommand \one { o }
{%
\IfNoValueTF {#1}
{\mathbf{1}  }
{\mathbf{1}\left\{ {#1} \right\} }%
}

\usepackage{tabularx, graphicx, wrapfig, tikz, subcaption}
\graphicspath{ {./images/} }

\usepackage{cancel}
\usepackage{cleveref}

\theoremstyle{definition}

\numberwithin{equation}{section}

\newtheorem{definition}{Definition}[section]

\newtheorem{example}[definition]{Example}

\newtheorem{lemma}[definition]{Lemma}

\newtheorem{theorem}[definition]{Theorem}

\newtheorem{proposition}[definition]{Proposition}

\newtheorem{remark}[definition]{Remark}

\usepackage{cleveref}
\crefname{section}{Section}{sections}
\crefname{subsection}{Subsection}{subsections}
\Crefname{section}{Section}{Sections}
\Crefname{subsection}{Subsection}{Subsections}

\Crefname{figure}{Figure}{Figures}

\crefformat{equation}{\textup{#2(#1)#3}}
\crefrangeformat{equation}{\textup{#3(#1)#4--#5(#2)#6}}
\crefmultiformat{equation}{\textup{#2(#1)#3}}{ and \textup{#2(#1)#3}}
{, \textup{#2(#1)#3}}{, and \textup{#2(#1)#3}}
\crefrangemultiformat{equation}{\textup{#3(#1)#4--#5(#2)#6}}%
{ and \textup{#3(#1)#4--#5(#2)#6}}{, \textup{#3(#1)#4--#5(#2)#6}}{, and \textup{#3(#1)#4--#5(#2)#6}}

\Crefformat{equation}{#2Equation~\textup{(#1)}#3}
\Crefrangeformat{equation}{Equations~\textup{#3(#1)#4--#5(#2)#6}}
\Crefmultiformat{equation}{Equations~\textup{#2(#1)#3}}{ and \textup{#2(#1)#3}}
{, \textup{#2(#1)#3}}{, and \textup{#2(#1)#3}}
\Crefrangemultiformat{equation}{Equations~\textup{#3(#1)#4--#5(#2)#6}}%
{ and \textup{#3(#1)#4--#5(#2)#6}}{, \textup{#3(#1)#4--#5(#2)#6}}{, and \textup{#3(#1)#4--#5(#2)#6}}

\crefdefaultlabelformat{#2\textup{#1}#3}

\newrefformat{eq}{(\ref{#1})}
\newrefformat{thm}{Theorem~\ref{#1}}
\newrefformat{th}{Theorem~\ref{#1}}
\newrefformat{chap}{Chapter~\ref{#1}}
\newrefformat{sec}{Section~\ref{#1}}
\newrefformat{seca}{Section~\ref{#1}}
\newrefformat{algo}{Algorithm~\ref{#1}}
\newrefformat{fig}{Fig.~\ref{#1}}
\newrefformat{tab}{Table~\ref{#1}}
\newrefformat{rmk}{Remark~\ref{#1}}
\newrefformat{clm}{Claim~\ref{#1}}
\newrefformat{def}{Definition~\ref{#1}}
\newrefformat{cor}{Corollary~\ref{#1}}
\newrefformat{lemma}{Lemma~\ref{#1}}
\newrefformat{prop}{Proposition~\ref{#1}}
\newrefformat{pr}{Proposition~\ref{#1}}
\newrefformat{property}{Property~\ref{#1}}
\newrefformat{app}{Appendix~\ref{#1}}
\newrefformat{apx}{Appendix~\ref{#1}}
\newrefformat{ex}{Example~\ref{#1}}
\newrefformat{exer}{Exercise~\ref{#1}}
\newrefformat{soln}{Solution~\ref{#1}}
\newrefformat{lmm}{Lemma~\ref{#1}}

\def\ba{{\boldsymbol a}}   \def\bA{\boldsymbol A}  
   \def\bB{\boldsymbol B}  
   \def\bC{\boldsymbol C}  
\def\bd{{\boldsymbol d}}     
     \def\EE{\mathbb{E}}

\def\bm{\boldsymbol m}   \def\bM{\boldsymbol M}

     \def\PP{\mathbb{P}}

   \def\bT{\boldsymbol T}  
\def\bu{{\boldsymbol u}}     
     
\def\bw{\boldsymbol w}     
\def\bx{\boldsymbol x}     
\def\by{\boldsymbol y}

\def\11{\mathbbm{1}}

\def\calC{{\cal  C}} 
 
\def\calE{{\cal  E}} 
 
\def\calG{{\cal  G}}

\def\calS{{\cal  S}}

\newcommand{\bfsym}[1]{\ensuremath{\boldsymbol{#1}}}

             \def\bSigma{\bfsym \Sigma}

\def\P{{\mathbb{P}}}

\def\E{{\mathbb{E}}}
\def\R{{\mathbb{R}}}
\def\N{{\mathbb{N}}}

\def\Tr{{{\rm Tr}}}

\newcommand{\argmin}{{\rm argmin}}
\def\%#1\%{\begin{align}#1\end{align}}
\def\$#1\${\begin{align*}#1\end{align*}}

\graphicspath {{fig/}}

\usepackage{ifthen}
\newboolean{aos}
\setboolean{aos}{FALSE}

\newcommand{\rvprobunion}{2L\exp\left(-{(\sqrt{n}-\sqrt{m})^2 \over 8}\right)}

\usepackage{bm}

\newcommand{\KL}{\operatorname{KL}}

\newcommand{\Sig}{\mathbf{\Sigma}}

\newcommand{\Sohat}{\widehat{\mathbf \Sigma}_o}
\newcommand{\Sotil}{\widetilde{\mathbf \Sigma}_o}
\newcommand{\Stil}{\widetilde{\mathbf \Sigma}_o}
\newcommand{\So}{\mathbf \Sigma_o}

\newcommand{\bde}{\boldsymbol{\delta}}

\newcommand{\av}{{\bf a}}

\newcommand{\xv}{{\boldsymbol x}}

\newcommand{\uv}{{\boldsymbol u}}
\newcommand{\vvv}{{\boldsymbol v}}

\newcommand{\uvh}{\widehat{\boldsymbol u}}

\newcommand{\bxi}{{\boldsymbol{\xi}}}

\newcommand{\lw}{{\vv{\boldsymbol w}}}
\newcommand{\ly}{{\vv{\boldsymbol y}}}
\newcommand{\lv}{{\vv{\boldsymbol v}}}
\newcommand{\netdelta}{{N_{\delta_{\rm net}}(\calC_k)}}

\newcommand{\Cc}{\mathcal{C}}
\newcommand{\Ec}{\mathcal{E}}

\newcommand{\yv}{{\boldsymbol y}}

\newcommand{\zv}{{\boldsymbol z}}

\newcommand{\wv}{{\boldsymbol w}}

\newcommand{\xvh}{\widehat{\xv}}
\newcommand{\xvt}{\widetilde{\xv}}

\ifx\eqref\undefined
\newcommand{\eqref}[1]{~(\ref{#1})}
\fi
\ifx\mod\undefined
\def\mod{\mathop{\rm mod}}
\fi

\usepackage{bm}

\def\exp{\mathop{\rm exp}}

\def\PP{\mathbb{P}}

\def\diag{\mathop{\rm diag}}

\newcommand{\reals}{{\mathbb R}}
\newcommand{\emax}{{E_{\max}}}
\newcommand{\emin}{{E_{\min}}}

\newcommand{\pth}[1]{\left( #1 \right)}
\newcommand{\qth}[1]{\left[ #1 \right]}
\newcommand{\sth}[1]{\left\{ #1 \right\}}

\newcommand{\cnms}{{C_{n,m,\sigma_z}}}

\newcommand{\esing}{{\calE_{\rm sing}}}
\newcommand{\emaxsing}{{\calE_{\rm maxsing}}}
\newcommand{\esingp}{{\calE'_{\rm sing}}}

\usepackage{lineno}

\makeatletter
\def\thanks#1{\protected@xdef\@thanks{\@thanks \protect\footnotetext{#1}}}
\makeatother

\title{Minimax Analysis of Estimation Problems in Coherent Imaging}

\author{Hao Xing$^{\dagger}$, Soham Jana$^\dagger$, and Arian Maleki\thanks{$^\dagger$Equal contributions. H.~X. is with the Department of Mathematics, City University of New York - the Graduate Center, NY, USA (email:\url{hao.xing43@gc.cuny.edu}). S.~J. is with the Department of Applied and Computational Mathematics and Statistics, University of Notre Dame, IN, USA, (email:\url{soham.jana@nd.edu}). A.~M. is with the Department of Statistics, Columbia University, NY, USA.(email: \url{mm4338@columbia.edu}).}}
\date{}

\begin{document}
	\maketitle
	\begin{abstract}
		
		Unlike conventional imaging modalities, such as magnetic resonance imaging, which are often well described by a linear regression framework, coherent imaging systems follow a significantly more complex model. In these systems, the task is to estimate the unknown image ${\boldsymbol x}_o \in \mathbb{R}^n$ from observations ${\boldsymbol y}_1, \ldots, {\boldsymbol y}_L \in \mathbb{R}^m$ of the form
		$
		{\boldsymbol y}_l = A_l X_o {\boldsymbol w}_l + {\boldsymbol z}_l, l = 1, \ldots, L,
		$
		where $X_o = \mathrm{diag}({\boldsymbol x}_o)$ is an $n \times n$ diagonal matrix, 
		${\boldsymbol w}_1, \ldots, {\boldsymbol w}_L \stackrel{\text{i.i.d.}}{\sim} \mathcal{N}(0,I_n)$ represent speckle noise, 
		and ${\boldsymbol z}_1, \ldots, {\boldsymbol z}_L \stackrel{\text{i.i.d.}}{\sim} \mathcal{N}(0,\sigma_z^2 I_m)$ denote additive noise. 
		The matrices $A_1, \ldots, A_L \in \mathbb{R}^{m \times n}$ are known forward operators determined by the imaging system. Here, $m$ can be less than $n$, which makes the linear model ill conditioned, and we collect $L$ measurements of the whole setup and try to quantify the loss.  
		
		The fundamental limits of conventional imaging systems have been extensively studied through sparse linear regression models. However, the limits of coherent imaging systems remain largely unexplored. The goal of this paper is to close this gap by showing, for the first time in the literature, that the minimax risk of estimating ${\boldsymbol x}_o$ scales as
		$\frac{\max\{\sigma_z^4,\, m^2,\, n^2\}\, k \log n}{m^2 n L},$
		where the parameter $k$ quantifies the effective complexity of the image class.
	\end{abstract}

		
		
%
%
	
	\section{Introduction}
	
	\subsection{Motivation and main objective}
	Coherent imaging technology, which uses coherent light sources such as lasers to illuminate the object of interest, underpins many modern imaging systems. 
	Examples include Optical Coherence Tomography (OCT)~\cite{Schmitt}, ultrasound imaging~\cite{Achim}, Synthetic Aperture Radar (SAR)~\cite{Lopez-MartinezFabregas, Dasari}, digital holography~\cite{Bianco}, and near-infrared spectroscopy (NIRS)~\cite{ortega2019contribution}. Compared to other imaging modalities, coherent imaging systems are affected by a complex form of distortion known as speckle noise \cite{racine1999speckle}.
	
	The general mathematical problem that arises in coherent imaging systems is to estimate a signal or image $\xv_o \in \mathbb{R}_+^{n}$ from measurements $\yv$ of the form
	\begin{equation}\label{eq:speckle_mod1}
		\yv = A X_o \wv + \zv,
	\end{equation}
	where $X_o \in \R^{n\times n}$ is a diagonal matrix whose diagonal elements are the same as $\xv_o$, $\wv \in \R^n$, and $\zv\in \R^m$ represent the speckle and additive noises, respectively. In this model, $A\in \R^{m\times n}$ is a known matrix, called the forward operator of the imaging system. Note that, compared to linear regression problems, which are popular for other types of imaging such as MRI and CT, the relationship between $\yv$ and $\xv_o$ is further distorted by the speckle noise $\wv$. In many applications, the speckle noise is "fully developed", which means that the elements of $\wv$ are i.i.d. $N(0,1)$. A standard model for the additive noise $\zv$ is also that it has i.i.d. $N(0, \sigma_z^2)$.

	Before discussing our modeling assumptions, we first highlight an essential technique in coherent imaging systems, namely \emph{multilook} or \emph{multishot} measurements. 
	It is widely recognized in the coherent imaging community that estimating $\xv_o$ from a single measurement of the form in~\eqref{eq:speckle_mod1} is challenging, and in most applications the reconstruction quality is insufficient. 
	Therefore, in many settings, such as SAR and digital holography~\cite{de1998speckle, Argenti, bate2022experimental}, multiple measurements of the same scene are acquired. 
	More specifically, one collects measurements $\yv_1, \ldots, \yv_L$ of the form
	\begin{equation}\label{eq:multilook speckle noise model}
		\yv_l = A_l X_o \wv_l + \zv_l, \quad l=1,2,..,L
	\end{equation}
	where $L$ is referred to as the number of \emph{looks}, and $A_1,\dots,A_L$ represent the forward operators of different shots. There are a few points that we should clarify about this multilook system:
	\begin{itemize}
		\item In practice, effort is made to ensure that $\wv_1, \wv_2, \ldots, \wv_L$ are independent of each other. 
		Similarly, the additive noise vectors $\zv_1, \dots, \zv_L$ are typically assumed to be independent across measurements and also independent of $\wv_1, \ldots, \wv_L$.

		\item The forward models across looks may differ or be the same, depending on the technology used. For instance, if the wavelength of the illuminating light changes, then different $A_i$'s will be observed, but if phase masks are used on the path of the illuminating light, then we will have $A_1= A_2= \cdots=A_L$. 
	\end{itemize}
	
	In our mathematical model, we assume that the variance of the speckle noise is equal to $1$. This is without any loss of generality. Consider the case where the standard deviations of the multiplicative noise and the additive noise are ${\sigma}_w$ and ${\sigma}_z$, respectively. Then, by dividing the sensor measurements $\yv_l$ by $\sigma_w$, we obtain an equivalent system of measurements:
	\[
	\tilde{\yv}_l = A X_o \tilde{\wv}_l + \tilde{\zv}_l.
	\]
	Here, we have defined $\tilde{\yv}_l := \frac{\yv_l}{\sigma_w},
	\tilde{\wv}_l := \frac{\wv_l}{\sigma_w}, \tilde{\zv}_l := \frac{\zv_l}{\sigma_w}$. In view of this transformation, we have 
	\[
	\tilde{\wv}_l \sim \mathcal{N}(0, I), 
	\quad \tilde{\zv}_l \sim \mathcal{N}\!\left(0, \tfrac{\sigma_z^2}{\sigma_w^2} I\right),
	\]
	which is consistent with \eqref{eq:multilook speckle noise model}. Hence, without a loss of generality, we set the variance of the speckle noise to $1$ and discuss the changes in $\sigma_z$. 
	
	{\begin{remark}
		In the formulation $y = AXw + z$, the variance of the speckle noise $\sigma_w^2$ is assumed to be known or fixed (typically set to $\sigma_w^2 = 1$). Because $w$ acts as a multiplicative/scaling noise process on $X$, the term $Xw$ represents a scaled version of the signal. As a result, the signal magnitude $X$ and the speckle variance $\sigma_w^2$ are inherently coupled: scaling $X$ by a factor $c$ while dividing $\sigma_w$ by $c$ yields an equivalent statistical distribution for $Xw$. In several relevant works, e.g., in the area of Bayesian reconstruction frameworks for digital holography and coherent imaging, this scale ambiguity is a standard property. For instance, in \cite{pellizzari2017phase}, using a Bayesian framework, the authors establish a similar forward observation model (e.g., Equation 5) for modeling the image reflectance. In their formulation, the primary modality estimated is a scaled version of the underlying reflection coefficient/signal, and they implicitly assume a normalized speckle variance ($\sigma_w^2 = 1$) to fix the global scale parameter of the signal. In other works where speckle variances are required to be estimated, one considers each image fragments that can be considered homogeneous and then
		estimates of local variance obtained for scanning windows (blocks) that belong to such homogeneous fragments to produce estimates close to the true variance $\sigma_w^2$ \cite{abramova2015blind}. For example, ratio of the standard deviation to the mean of pixel values in the homogeneous featureless areas can be used to estimate the variance parameter \cite{lee1992unsupervised}.
	\end{remark}
}
	
	Despite the widespread use of coherent imaging technology across many applications, the theoretical aspects of the associated estimation problems remain largely unexplored. The main goal of this paper is to help fill this gap by addressing the following questions:
	
	\begin{enumerate}
		\item How do $m, n, L$ and $\sigma_z^2$ affect the accuracy of the estimates of $\xv_o$? 
		
		\item Is there any gain in using different $A_1, A_2, \ldots, A_L$ compared to $A_1= A_2= \ldots = A_L$? 
		
		\item How do additive and multiplicative noise compare in their impact on the estimation error?
		
	\end{enumerate}
	
	To address these questions and to provide insight into the estimation challenges that arise in coherent imaging systems, we seek to characterize the minimax risk associated with this problem. In particular, we aim to quantify 
	\begin{align}\label{eq:minimax:first}
		R_2(\Cc,m,n,\sigma_z^2) 
		:= \inf_{\xvh} \sup_{\xv_o \in \Cc} 
		\E\!\left[ \frac{\|\xvh - \xv_o\|_2^2}{n} \right],
	\end{align}
	where $\xvh$ is any measurable estimate that has access to $\yv_1, \ldots, \yv_L$ and $A_1, A_2, \ldots, A_L$, and $\Cc$ denotes the set of all possible options for $\xv_o$. Notably, \cite{chen2026multilook} provides an upper bound to the above minimax error in the restrictive setup of $m< n$, $\sigma_z =0$, and $\calC$ being particularly chosen as a class of Deep Image Priors \cite{heckel2018deep,ulyanov2018deep}. In contrast, we will completely characterize the minimax rate by obtaining matching lower and upper bounds in a general setup of $m,n,\calC$. Prior work in sparse linear regression and compressed sensing shows that the minimax risk is strongly influenced by the choice of $\Cc$. In the next section, we first describe our choice of $\Cc$, and then discuss our contributions and our responses to the questions we raised above. 
	
	\subsection{Notations}
	Throughout the paper, for the sake of clarity, all matrices of sizes $m\times n$, $m\times m$, and $n\times n$ (without dependence on the number of looks $L$) are represented by uppercase italic letters such as $A_l$, $\Sigma_l$, and $X$. We use boldface uppercase letters, e.g. $\bA, \bB$ and $\Sig$, when the sizes of matrices depend on $L$. These matrices are often constructed by stacking matrices of smaller sizes and may have sizes such as  $mL\times (m+n)L$ and $(m+n)L\times (m+n)L$. For a matrix $A$, $\smax(A)$ and $\smin(A)$ denote the maximum and minimum singular values of $A$. Furthermore $\nm{A}_2=\smax(A)$ and $\nmhs{A}$ denote the spectral norm and Hilbert-Schmidt norms of $A$, respectively. Boldface lowercase letters such as $\xv$ are used for vectors for sizes $m$ or $n$ (again, no dependence on $L$). Arrows above the vectors emphasize that the dimensions of the vectors depend on $L$. Again, such vectors are constructed by stacking $L$ lower-dimensional vectors. For a vector $\xv=(x_1,...,x_n)$, we let $\xv^2:=(x_1^2,...,x_n^2).$ Given two sequences $\{a_n\}_{n\in \N}$ and $\{b_n\}_{n\in \N}$, we say $a_n=O(b_n)$, or equivalently $a_n\ll b_n$, $b_n=\Omega(a_n)$, if there exist constants $C>0$, and $M>0$, such that for all $n>M$, $|a_n|\le C |b_n|$. We say $a_n=o(b_n)$, or equivalently $b_n=\omega(a_n)$, if $\lim_{n\to\infty} {a_n}/{b_n} =0$. We write $a_n = \Theta (b_n)$ if $a_n = O(b_n)$ and $b_n = O(a_n)$. For $a,b\in \R$, we denote $(a,b]_{\Z}:=\{x\in Z: a<x \le b\}.$
	
	\subsection{Organization of the paper}
	
	In Section~\ref{sec:main:contr}, we introduce our model, state the main assumptions, and present our primary contributions. Section~\ref{sec:relatedwork} reviews related work and compares our results with the existing literature. In Section~\ref{sec:preliminaries}, we provide the necessary preliminaries for our analysis. The remaining sections are devoted to the proofs of the theorems stated in Section~\ref{sec:main:contr}.

	\section{Our main contributions}\label{sec:main:contr}
	In this section, we first discuss our choice of the parameter set $\Cc$ in \eqref{eq:minimax:first}, and then present our results to answer the questions we raised in \prettyref{sec:maintheory}.

	\subsection{The choice of $\mathcal{C}$}\label{sec:Cc}
	
	Inspired by developments in the fields of sparse regression and compressed sensing, we note that the structure of $\xv_o$, plays a crucial role in determining the accuracy of the estimates. As will be clarified later in this section, sparsity is not useful for coherent imaging systems.  Hence, in this paper we work with a more general notion of "structure". This notion allows us to cover not only the class of $k$-sparse vectors, but also the more modern classes developed in the field of neural networks, such as the class of untrained networks. For a compact set $\mathcal{C} \subset \mathbb{R}^n$, let $N_{\e}(\Cc)$ denote its covering number under the $\ell_2$ metric, namely the least number of $\ell_2$-balls covering $\Cc$.
	
	\begin{definition}\label{def:polcomplex} 
		We say that $\mathcal{C} \subset \mathbb{R}^n$ satisfies polynomial complexity of order $k$ if there exist constants $a>0,b\ge 0$ independent of $k$ and $n$ such that  
		\begin{equation}\label{eq:degree k polynomial growth}
			N_{\e}(\mathcal C)\le \pr{an^b \over \e}^k.
		\end{equation}
	\end{definition}
	
	{ The polynomial complexity is fundamental in nature, and arises directly from the geometric bounds on the function class and its parametric representation:
		\begin{enumerate}
			\item \textbf{The $n^b$ term:} In an $n$-dimensional space, the radius of any set with bounded coordinates (which here corresponds to the $n$-dimension set of possible signals) scales at most on the order of $\sqrt{n}$. To accommodate this dimensional scaling generally, we incorporate the $n^b$ factor in the bound's numerator where the magnitude of $b$ is of a constant order.
			\item \textbf{The exponent $k$:} For a function class parameterized by $k$ parameters with each having a compact support, , the metric entropy (covering number) standardly scales as ${O}\left(\left(\frac{\text{radius}}{\epsilon}\right)^k\right)$. In our case, $k$ corresponds to dimension of the embedding space of high-dimensional signal $\bx_o$.
		\end{enumerate}
		Combining these geometric and parametric constraints yields a covering number bounded by $\left(\frac{a n^b}{\epsilon}\right)^k$.}
	
	Before proceeding, we review several sets with polynomial complexity of order $k$ to establish the usefulness of this definition. The proof of the following results are provided in the appendix. Our first example can serve as a proxy for images generated by neural network architectures such as deep image priors \cite{ulyanov2018deep}, implicit neural representations \cite{sitzmann2020implicit}, and autoencoders \cite{bank2023autoencoders}. These models have been extensively used as reliable and accurate models for images.

	\begin{example}\label{exam:LipschitzFunc previous}
		Consider $k \ll n$ and let $g: \mathbb{R}^k \rightarrow \mathbb{R}^n$ denote a Lipschitz function with a Lipschitz constant $M$. Define 
		\[
		\mathcal{C} = \{\xv \in \mathbb{R}^n \ | \ \xv= g(\theta) \ \ {\rm for \ some } \ \theta \in [0,1]^k\}.
		\]
		Then, $N_{\e} (\Cc) \leq \left( \frac{2M\sqrt{k}}{\epsilon}+1\right)^k$. Note that when $\epsilon < 2M\sqrt{k}$, we have $\frac{2M\sqrt{k}}{\epsilon}+1 < \frac{3M\sqrt{k}}{\epsilon}$, and hence we can also have $N_{\e} (\Cc) \leq \left( \frac{3M\sqrt{k}}{\epsilon}\right)^k$.
	\end{example}
	
	\prettyref{exam:LipschitzFunc previous} covers a wide range of examples. For instance, in the literature of neural networks, it has been conjectured that the output of certain neural networks, such as implicit neural representation networks \cite{sitzmann2020implicit} and deep image priors \cite{heckel2018deep} can generate all natural images as the parameters of the networks change. Note that the number of parameters of these networks can be interpreted as $k$ in \prettyref{exam:LipschitzFunc previous}. Often times the number of parameters is much smaller than the ambient dimension of the signal that is generated by the network. {
		
		\begin{remark}
			Prior works on architectures such as the Deep Decoder, such as \cite{heckel2018deep} and \cite{heckel2019denoising}, suggest that, in many realistic noisy settings, reducing parameter count through architectural constraints or fewer channels can still maintain strong performance, although the precise effect depends on the task and reconstruction setting.
			More specifically, as noted in \cite{heckel2018deep} (p.~5), the authors model the imaging output $\bx$ as $\bx=G(\bC)$, where $\bC$ denote the network parameters, corresponding to $g(\theta)$ in our example. The effective reconstruction of an RGB image with $786{,}432$ pixels ($n$ in our example) can be achieved using the Deep Decoder with a constant number of channels across convolutional layers. When the number of channels is set to 64 and 128, the corresponding parameter counts ($k$ in our example) are $25{,}536$ and $100{,}224$, respectively, both significantly smaller than the number of image pixels. This observation suggests that even when $k \ll n$ (where $k$ is the number of network parameters and $n$ the number of pixels), the network can still represent natural images effectively.
		\end{remark}
	}
	Our second example offers upper bounds for the covering numbers of $k$-sparse vectors.
	
	\begin{example}\label{exam:sparse}
		For the set $C = B_2(1) \cap \mathcal{S}_k$, where $\mathcal{S}_k = \{\xv \in \mathbb{R}^n \ | \ \|\xv\|_0 \leq k \},$ we have
		\[
		\Big(\frac{1}{\epsilon} \Big)^k \leq N_{\epsilon} (\Cc) \leq {n \choose k} \pr{\frac{2}{\e}+1}^k \leq \pr{\frac{2n}{\e}+n}^k \le \pr{3n \over \e}^k.
		\]
	\end{example}
	
	The following example generalizes the above result to a wide range of models.
	
	\begin{example}\label{exam:singular value and Minkowski}
		Let $D \in \mathbb{R}^{n \times n}$ denote a matrix with the maximum singular value $\sigma_{\max}(D)$. Suppose that $\Cc \subset \{D \theta  \ | \ \theta \in \mathcal{S}_k \cap B_2(0,1)\}$. Then, $N_{\e} (\Cc) \leq \Big(\frac{3\sigma_{\max} (D) n}{\epsilon} \Big)^k$. 
	\end{example}
	
	Simple applications of the above result are piecewise constant and piecewise polynomial functions. The following lemma, proves this claim for the class of piecewise constant vectors.
	
	\begin{example}\label{exam:piececonstant_1}
		Define $D_{(1)} \in \mathbb{R}^{n \times n}$ as a matrix whose diagonal elements are equal to $1$ and the immediate super diagonal elements are equal to $-1$ (all other entries are zero). Suppose that for every $\xv \in \Cc$, $D_{(1)} \xv \in \mathcal{S}_k \cap B_2(0,1)$. Then, we have $N_{\e}(\Cc) \leq \left(\frac{3n^2}{\e}  \right)^k$.
	\end{example}
	
	Note that if $\xv$ is a constant vector, i.e., all its entries have the same value, then all elements of $D_{(1)} \xv$ except for the last one are equal to zero. Therefore, if we assume that $D_{(1)} \xv \in \mathcal{S}_k$, it follows that $\xv$ is a piecewise constant vector with at most $k$ jumps (changes) in its values.

	There are several ways to extend the above example to piecewise polynomial functions of degree at most $P$. We do the following simple extension.

	\begin{example}\label{exam:piecepolynomial_1}
		Define $D_{(1)} \in \mathbb{R}^{n \times n}$ as a matrix whose diagonal elements are equal to $1$ and the immediate super diagonal elements are equal to $-1$. Define $D_{(P+1)} = (D_{(1)})^{P+1}$.  Suppose that for every $\xv \in \Cc$, $D_{(P+1)} \xv \in \mathcal{S}_k \cap B_2(0,1)$. Then, we have $
		N_{\e}(\Cc) \leq \left(\frac{3n^{P+2}}{\e}  \right)^k$.
	\end{example}
	Note that if $f: [0,1] \rightarrow \mathbb{R}$ is a polynomial of degree $P$, and $\bm{x}_i = f(i/n)$, it follows from a well-known fact in the theory of forward difference operators (see \cite[Section 5.3]{GrahamKnuthPatashnik1994}) that all elements of $D_{(P+1)}(\bx)$ are equal to zero, except possibly for the last $P+1$ elements. Hence, the set $\mathcal{C}$ in \prettyref{exam:piecepolynomial_1} can be viewed as discretized piecewise polynomial vectors. 
	
	Inspired by all the examples above, in our theoretical results we will be assuming that $\bx_o$ in \eqref{eq:speckle_mod1} is from a set $\mathcal{C}$ that has a polynomial complexity of order $k \ll n$. 
	
	\begin{remark}
		{ In standard image processing, the signal of interest, e.g., two-dimensional digital images are frequently vectorized into a one-dimensional signal representation. Consequently, the simplest modeling paradigm involves imposing structural constraints directly on these vectorized 1D representations. While this simplified 1D framework serves to demonstrate the core mechanics of our approach without unnecessary visual or computational complexity, we recognize the importance of explicitly extending these constraints to native 2D image structures. We leave the generalization to 2D signals for future work.}
	\end{remark}
	
	\subsection{Main theoretical result for independent $A_i$s}\label{sec:maintheory}

	As we discussed before, we consider the problem of estimating $\xv_o$ from the observations 
	\begin{equation}\label{eq:model}
		\yv_l = A_l X_o \wv_l + \zv_l, \text{ 
			for } l=1,\ldots,L,
	\end{equation}
	under the assumption $\wv_1, \wv_2, \ldots, \wv_L \overset{i.i.d.}{\sim} N(0, I)$, and $\zv_l \overset{i.i.d.}{\sim} N(0, \sigma_z^2 I)$. Our main goal is to characterize the minimax risk of the estimation problem in \eqref{eq:model} defined as:  
	\begin{align} \label{eq:minimaxrisk:varying}
		R_2(\Cc,m,n,\sigma_z):=\inf_{\xvh}\sup_{\xv_o\in \mathcal{C}}\E\br{\frac{\nm{\xvh-\xv_o}_2^2}{n}}.
	\end{align}
	Our first theorem obtains an upper bound for this quantity:
	Let $\mathcal{F}_{a,b,k,n}$ denote all subsets of $[x_{\min}, x_{\max}]^n$ whose $\e$-covering number is upper bounded by $\pr{an^b \over \e}^k$.
	
	\begin{theorem}\label{thm:main-genr}
		Suppose that $A_1,...,A_L$ are independent $m\times n$ matrices and have i.i.d.~$N(0,1)$ entries. Suppose that $\mathcal{C}_k \in \mathcal{F}_{a,b,k,n}$.  If $mL\le n^4 k\log n$, then
		\begin{equation*}
			R_2(\Cc_k,m,n,\sigma_z)
			= O_{x_{\max},x_{\min},a,b}\pr{\min\pr{\frac{\max(\sigma_z^4,m^2,n^2)k\log n}{m^2 n L},1}}.
		\end{equation*}
		
	\end{theorem}
	
	The above theorem obtains an upper bound for the minimax risk that holds for any $\mathcal{C}_k \subset [x_{\min}, x_{\max}]^n$ that satisfies polynomial complexity of order $k$. Before discussing the assumptions made in the above theorem, let us discuss the sharpness of this upper bound.  
	
	\begin{theorem}\label{thm:lower_bnd}
		Suppose $a\ge \xmax-\xmin$ and $b\ge 1$. If 	$\log m =\Theta(\log n)$, $\log L=O(\log n)$, and there exists $\e\in (0,1/2)$ such that $k \le n^{1-2\e}$, and that $\max(\sigma_z^4,m^2,n^2)k\log n \le m^2n^{1-\e}L$, then we have
		\begin{equation*}
			\sup_{\Cc  \in \mathcal{F}_{a,b,k,n}} R_2(\Cc,m,n,\sigma_z)
			= \Omega_{\e,x_{\max},x_{\min}}\pr{\frac{\max(\sigma_z^4,m^2,n^2)k\log n}{m^2 n L}}.
		\end{equation*}
	\end{theorem}
	
	 { The lower bound here holds on the case where $a\geq x_{\max}-x_{\min},b\geq 1$. Our lower bounds hold for any $a,b$ in that region and is a generalized bound, compared to the result regarding the upper bound, where the bounds may depend on the particular values of $a,b$}. Note that \prettyref{thm:lower_bnd} can be established for any $a,b>0$; see the appendix for discussions.
	 
	 Before we discuss the implications of our result, we discuss some of the assumptions we have made in the above theorems. A natural question is why we did not adopt the standard notion of sparsity widely used in sparse regression and imaging systems that fit well within the framework of linear regression. We mention two reasons below:
	
	\begin{enumerate}
		\item In imaging sciences, it is often the case that the vector $\xv$ is not sparse itself. In fact, some linear transformation of the vector, e.g. wavelet or Fourier transform of $\xv$ is sparse \cite{donoho1995wavelet, donoho1998minimax}. Suppose that $\xv = F \uv$, where $\|\uv\|_0 \leq k$. Then, in linear regression,  one can write the measurement $\yv= A\xv +\zv$ as $\yv= \tilde{A} \uv + \zv$, where $\tilde{A} = AF$. Hence, the problem of imaging when linear model is accurate, is equivalent to the problem of sparse linear regression. As is clear, because of the nature of the speckle noise, we \textbf{cannot} transform the estimation of $\xv$ from the observation $\yv_1, \yv_2, \ldots, \yv_L$ to an estimation of a sparse vector (with a different design matrix).  
		
		\item Because of the nature of speckle noise (that is multiplied by $\xv$), the estimation of sparse signals are easier than the estimation of the non-sparse signals. Intuitively speaking, this is due to the fact that sparse vectors, automatically remove most of the speckle noises. In other words, out of the $n$-speckle noise elements that are often present in these systems, $n-k$ of them will be multiplied by zeros during the measurement process and will not have a major impact on the estimation problem. To confirm this intuition rigorously, the next theorem shows that the minimax risk of estimating sparse vectors from  \eqref{eq:model} is much smaller than the bounds presented in \prettyref{thm:main-genr} and \prettyref{thm:lower_bnd}. 
		
	\end{enumerate}
	Let
	\begin{align}\label{sparse signal class}
		\mathcal{S}_k^{\rm bdd}:=\cb{\xv\in \R^n: \nm{\xv}_0\le k, x_i=0 \text{ or } 0<x_{\min}\le x_i \le x_{\max}}.
	\end{align}
	
	\begin{theorem}\label{thm:ksparse}
		There are constants $c_{\xmax,\xmin}, C_{\xmax,\xmin}$ such that for $k\log(en/k)\le m$
		\begin{align*}
			c_{\xmax,\xmin}{k \over nL}\le R_2(\mathcal{S}_k^{\rm bdd}, m,n,k,L, \sigma_z)
			\le C_{\xmax,\xmin}\pr{{k \over nL} + \frac{\sigma_z^2 k\log(n/k)}{mn}}
		\end{align*}
		In particular, if $\sigma_z^2 L\log(n/k)\le m$, the upper and lower bounds have the same order $\frac{k}{nL}$.
	\end{theorem}
	
	{ 
		\begin{remark}
			The above result displays that for the case $\sigma^2_zL \log(n/k) \leq m$ the estimation rate is $\Theta(k/Ln)$ which is very small, implying that the sparse signal setting is much simpler to reconstruct. The above regime is standard in the literature as the number of looks $L$ is often smaller than the number of sensors $m$, due to the cost involved in replicating the experimental setups. In particular in the case where the number of looks $L$ is a constant, the estimation rate is $\Theta(k/n)$, and is achieved by any $k$-sparse estimator. This shows that in the above standard regime, naive estimation strategies are sufficient, however, they do not provide the optimal error rate in case of polynomial complexity image space.
		\end{remark}
	}
	The proof of this result is provided in the appendix. By comparing \prettyref{thm:ksparse} with 
	Theorems \ref{thm:main-genr} and \prettyref{thm:lower_bnd}, we can see that estimating $k$-sparse signals can be much easier than estimating other signal classes satisfying polynomial growth of order $k$ that includes for instance, the signals that are sparse in an orthonormal basis such as wavelet (see \prettyref{exam:singular value and Minkowski}).

	The remaining assumptions of \prettyref{thm:main-genr} and \prettyref{thm:lower_bnd} are only technical and relatively minor; they could likely be removed, though doing so would make the proof less transparent. We discuss these assumptions briefly below:
	\begin{enumerate}
		\item $mL \leq n^4 k \log n$: This assumption appeared in \prettyref{thm:main-genr}. 
		In practice, obtaining more than $L>100$ independent looks is rarely feasible. 
		Since in most imaging applications $n$ is on the order of hundreds of thousands to millions, 
		this condition is typically satisfied.
		
		\item $\log (L) = O(\log n)$: This assumption appears in  \prettyref{thm:lower_bnd}. As we discussed before, in all applications, $L$ is much smaller than $n$. Hence, the assumption that $L$ is not growing too fast in terms of $n$ is a natural assumption in practice.  
		
		\item $\log m = \Theta (\log n)$: As will be discussed in the proof of  \prettyref{thm:lower_bnd} and might be even clear from the formulation of the problem, increasing $m$ beyond $n$ does not help in removing the speckle noise, and it only helps in removing the additive part of the noise. Hence, increasing $m$ to a very large number is not particularly helpful in reducing the risk, since unless the additive noise is too large, the errors induced by the speckle noise are the dominant part of the risk.  Note that increasing $m$ in imaging applications is equivalent to increasing the number of sensors which is costly.  As a result, there is no reason to increase $m$ much beyond $n$ in real-world applications, and again the assumption $\log m = \Theta (\log n)$ is a mild assumption. 
		
		\item $k \le n^{1-2\e}$: This assumption is used in  \prettyref{thm:lower_bnd}. It is always the case that $k \ll n$. Hence, again this is a mild assumption. However, at this stage it is unclear, whether this assumption is necessary.  
		
	\end{enumerate}

	\subsection{Interpretation of results}
	
	The following remarks discuss \prettyref{thm:main-genr} and \ref{thm:lower_bnd}.
	
	\begin{remark}[Difference of upper and lower bounds] \prettyref{thm:main-genr} holds for any set $\Cc$ that satisfies the polynomial complexity of order $k$. In contrast, the lower bound is obtained by taking the supremum of the minimax risk over all sets that satisfy polynomial complexity of order $k$. As \prettyref{thm:ksparse} illustrates, due to the nature of the speckle noise, which is multiplied by the entries of $\xv_o$, the estimation problem is easier for sparse vectors. Nevertheless, the supremum in the lower bound demonstrates that for certain sets $\Cc$ that satisfy polynomial complexity of order $k$, the upper bound established in \prettyref{thm:main-genr} is in fact sharp.
	\end{remark}
	
	Depending on the relative values of $m,n, \sigma_z$, the minimax risk can be simplified further:
	\begin{align*}
		R_2(\Cc_k,m,n,\sigma_z)=&\begin{cases}
			\Theta_{x_{\max},x_{\min}}\pr{\frac{kn\log n}{m^2L}}, &\text{ if }n\ge \max(m, \sigma_z^2);\\
			\Theta_{ x_{\max},x_{\min}}\pr{\frac{k\log n}{n L}}, &\text{ if }m \ge \max(n, \sigma_z^2);\\
			\Theta_{ x_{\max},x_{\min}}\pr{\frac{\sigma_z^4 k\log n}{m^2L}}, &\text{ if }\sigma_z^2 \ge \max(m, n). 
		\end{cases}
	\end{align*}
	In what follows, we offer some intuition to explain these bounds.
	
	\begin{remark}[When does the multiplicative noise dominate the additive noise?]
		Suppose that $\sigma_z^2$ is much smaller than $\max(m,n)$. In this case, the minimax risk is unaffected by the additive noise. To understand this phenomenon, suppose that we are working in the setting $m<n$. If we consider the model, $\yv= AX_o\wv + \zv$, let $\yv_i$ denote, the $i^{\rm th}$ element of $\yv$, then we have $\yv_i= \mathbf{a}_i^T X_o \wv+ \zv_i$. In this measurement, the variance of ${\rm var}({\mathbf{a}_i^T X_o \wv})$ is $\sum_{i=1}^n \xv_i^2 = \Theta (n)$, and ${\rm var} (\zv_i)= \sigma_z^2$. Hence, when $\sigma_z^2 = o(n)$, one would expect the additive noise to be negligible.  However, this heuristic does not explain why the additive noise does not matter when $m>n$ and $n \ll \sigma_z \ll m$. Again to provide some intuition on the impact of $m$, when $m>n$. Note that in this case, since $A^TA$ is an invertible matrix, we can calculate: $\tilde{\yv} = (A^TA)^{-1}\yv = X\wv + \tilde{\zv}$, where $\tilde{\zv} \sim N(0, \sigma_z^2 (A^TA)^{-1})$. While the additive noise $\tilde{\zv}$ is colored, and discussing signal-to-noise ratio on the individual elements does not necessarily provide an accurate information, note that in $\tilde{y}_i$ we have ${\rm var} (x_i w_i) = \Theta (1)$, and we can prove that ${\rm var}(\tilde{{z}_i}) = \Theta (\sigma_z^2/m)$ (See for instance \prettyref{lmm:singvalues} for the eigenvalues of $A$). Hence, in this case, again we can see that when $\sigma_z^2 \ll m$, the additive noise becomes negligible.   
	\end{remark}
	
	\begin{remark}[Comparison with \em{linear} imaging systems ] Again, consider the case $\sigma_z^2 \leq \max(m,n)$. In the classical regimes where the sparse linear regression problem is studied (e.g., \cite{bickel2009simultaneous}), namely $k \ll m \ll n$, the minimax risk of coherent imaging systems reduces to $
		\frac{kn \log n}{m^2 L}$.
		If $L$ is not too large, achieving a small risk requires $m \gg \sqrt{n}$. 
		This contrasts sharply with imaging systems based on linear regression, where obtaining a small minimax risk typically requires only $m \gg k \log n$. This is consistent with the general belief in the coherent imaging community that recovering images from measurements in coherent imaging systems is much more challenging than in imaging systems based on linear models.  
	\end{remark}

	\subsection{Fixed $A$ model}\label{sec:fixedA}
	In \prettyref{sec:maintheory}, we considered the setting in which the forward models $A_i$s across looks are independent. However, as we discussed before, in some multilook systems the forward models do not change across looks and we have
	\[
	A_1 =A_2 = \ldots = A_L.
	\]
	The main question that we would like to address in this section is whether either of these two multilook systems have an advantage over each other. To respond to this question, we aim to study the minimax estimation rate under the setting $A_1 = A_2 =\ldots = A_L$ and compare it with the result of \prettyref{sec:maintheory}. We define the minimax risk for this setting similar to what we defined before for the vase of different forward models. 
	\begin{equation}\label{eq:minimaxrisk:def:fixedA}
		R_2^{\dagger}(\Cc,m,n,\sigma_z) = \inf_{\xvh}\sup_{x_o\in \mathcal{C}}\E\br{{\nm{\xvh-\xv_o}_2^2/n}}.
	\end{equation}

	\begin{theorem}\label{thm:main-genr2}
		 
		Suppose $A_1 = A_2 = \ldots= A_L =A \in \mathbb{R}^{m \times n}$ and that $A_{ij} \overset{i.i.d.}{\sim} N(0,1)$. Furthermore, assume that $ \mathcal{C}_k \in \mathcal{F}_{a,b,k,n}$. Furthermore, assume $mL\le n^4 k\log n$. Then,
		\begin{equation*}
			R_2^{\dagger}(\Cc_k,m,n,\sigma_z)		=O_{x_{\max},x_{\min}}\pr{\min\pr{\frac{\max(\sigma_z^4,m^2,n^2)}{m^2 n} \frac{k\log n}{L}+
			{\frac{\max(\sigma_z^4,m^2,n^2)k \log n}{nm^3}},1}}.
		\end{equation*}
		
		{ In particular, if $m>L$}, then 
		\begin{equation*}
			R_2^{\dagger}(\Cc_k,m,n,\sigma_z)
			=O_{x_{\max},x_{\min}}\pr{\min\pr{\frac{\max(\sigma_z^4,m^2,n^2)k\log n}{m^2 n L},1}}.
		\end{equation*}
		
	\end{theorem}
	
	Similar to \prettyref{sec:maintheory}, the above theorem obtains an upper bound on the minimax risk for any $\mathcal{C} \subset [x_{\min}, x_{\max}]^n$ that  satisfies polynomial complexity of order $k$. Before discussing the implications of the above result, let us discuss the sharpness of this upper bound.

	\begin{theorem}\label{thm:lower_bnd2}
		Assume (i) $a\ge \xmax-\xmin, b\ge 1,\log m =\Theta(\log n),\log L=O(\log n)$, (ii) there exists $\e\in (0,1/2)$ such that $k \le n^{1-2\e}$, and (iii) $\max(\sigma_z^4,m^2,n^2)k\log n \le m^2n^{1-\e}L$. Then we have
		\begin{equation}\label{eq:main theorem lower bound for fixed A}
			\sup_{\Cc  \in \mathcal{F}_{a,b,k,n}} R_2^{\dagger}(\Cc,m,n,\sigma_z)
			=\Omega_{\e,x_{\max},x_{\min}}\pr{\frac{\max(\sigma_z^4,m^2,n^2)k\log n}{m^2 n L}}.
		\end{equation}
		
	\end{theorem}
	
	 { 
     
	Since the assumptions in the above two theorems are similar to the ones in \prettyref{thm:main-genr} and \prettyref{thm:lower_bnd} we will not discuss the assumptions again. However, there is one condition that we did not discuss before. This condition appears in the second part of \prettyref{thm:main-genr2} and indicates that if $m\geq L$, then 
	\begin{align}\label{eq:upper-fixed-match}
		R_2^{\dagger}(\Cc,m,n,\sigma_z)
		=&O_{x_{\max},x_{\min}}\pr{\min\pr{\frac{\max(\sigma_z^4,m^2,n^2)k\log n}{m^2 n L},1}}.
	\end{align}
	Note that the upper bound in  \prettyref{thm:main-genr2} matches the lower bound in \prettyref{thm:lower_bnd2} when $m\geq L$. Hence, this leads to the following questions
	\begin{enumerate}
		\item How strong is the assumption $m\geq L$?
		\item Is there any intuition why the $L$ independent term $\frac{\max(\sigma_z^4,m^2,n^2)k \log n}{nm^3}$ appears in  \prettyref{thm:main-genr2}, while not in \prettyref{eq:upper-fixed-match}?
	\end{enumerate}

	In response to the question above, note that in practice, $m$ is often chosen to be a large number that increases with the image dimension to efficiently capture its complexity. On the other hand, $L$ is often chosen to be a large constant between 2 and 100. Hence, the condition holds in general.
	
	Regarding the second question raised above, the necessity of the $L$-independent term is still unclear. However, an intuitive argument sheds some light on the difference between the fixed-$A$ and varying-$A$ cases. Suppose that the additive noise in \eqref{eq:model} is zero.  
	In the fixed design setting, it is straightforward to show that the statistics $\frac{1}{L} \sum_{i=1}^L \yv_l {\yv_l}^T$ is a sufficient statistic for $X_o$. If we fix $m,n$ and let $L \rightarrow \infty$, the sufficient statistic converges to $AX_o^2A^T$ in probability. In other words, the sufficient statistic converges to a linear transformation of $X_o^2$. However, note that recovering the exact $X_o^2$ from $AX_o^2A^T$ is not possible. In fact, in the most optimistic setting $AX_0^2 A^T$ offers $m(m+1)/2$ linearly independent observations of $X_o^2$. Hence, if $m (m+1)/2<k$, then $X_o^2$ is not exactly recoverable from $AX_o^2A^T$. When $k>m^2$, the term $\frac{\max(\sigma_z^4,m^2,n^2)}{nm^3}\geq \frac{(m^2+n^2)k\log n}{2nm^3}\geq {k\log n\over m^2}$ is quite large. Hence, this term is consistent with the intuition that when $k\geq m^2$ the error has to be large. 
	}
	
	\subsection{Fixed forward model or varying forward model? }
	
	By comparing the results of \prettyref{sec:maintheory} and \prettyref{sec:fixedA}, we observe that, in terms of minimax rates, there is no significant difference between fixed and varying forward models with respect to estimation accuracy. The upper bound that we have derived for the fixed-A model, has an extra $\frac{k \log m \log n}{m^2}$. However, this extra term is quite small for most practical settings and does not seem to be important. 
	
	{
	\begin{remark}
		The situation in coherent imaging systems is closely related to that in linear imaging problems (studied in compressed sensing), and in phaseless imaging problems (studied in phase retrieval). In many practical imaging systems, the underlying forward model involves structured Fourier-type measurements rather than i.i.d. Gaussian measurements. Nevertheless, Gaussian measurement matrices are widely used in theoretical analyses as analytically tractable proxies for real-world forward models. In particular, Gaussian designs are often viewed as representative of a broader class of sufficiently random and incoherent sensing operators.
		
		Substantial progress has been made in compressed sensing and phase retrieval in extending results initially developed for Gaussian measurements to more realistic structured sensing models. Examples include randomly subsampled Fourier measurements \cite{candes2008introduction} and coded diffraction patterns, where appropriate randomization can lead to properties (and proof techniques) analogous to those of Gaussian sensing matrices. 
		We expect that a similar phenomenon may hold for coherent imaging systems. Establishing the extent to which our theoretical results extend from Gaussian measurements to structured Fourier-type forward operators is an important and interesting direction we aim to pursue for our future research.  
	\end{remark}
	
	}
	
	\section{Related works}\label{sec:relatedwork}
	In this paper, we make the first attempt to establish the rate of the minimax risk for coherent imaging systems. 
	
	There is a substantial body of research on the theoretical characterization of the minimax risk for other imaging systems, such as MRI and CT, which are modeled by linear regression problems \cite{Tsybakov86, KT93, donoho1994minimax, bickel2009simultaneous, RWY11, candes2015slope, su2017false, weng2016phase, donoho2009message, metzler2016denoising, guo2024note,ghosh2025signal}, as well as crystallography and astrophotography, which are modeled by the phase retrieval problem \cite{chen2019gradient, chen2017solving, candes2015phase, zhang2017nonconvex, zhang2016provable, CaiLiMa16, hand2018phase, bakhshizadeh2020using, zhang2017nonconvex, ma2019optimization}.
	
	The theoretical properties of coherent imaging systems have only recently been explored in a few papers by subsets of the authors and their collaborators \cite{zhou2022compressed, chen2024bagged, chen2026multilook, malekian2025speckle}. We discuss the contributions of these papers, and compare our contributions with what are offered in those papers:
	
	\begin{enumerate}
		\item Speckle noise in nonparametric settings: \cite{malekian2025speckle} studies the speckle noise under the nonparametric regression settings, considering the model:
		\begin{equation*}
			y_i = f(x_i) \xi_i + \tau_i, \quad i = 1, 2, \ldots, n.
		\end{equation*}
		where $\xi_i$'s are i.i.d. $\mathcal N(0,1)$ and $\tau_i$'s are i.i.d. $\mathcal N(0,\sigma_{\tau}^2)$ random variables, $x_i=i/n, i=1,2,...,n$ are fixed design points, and unknown $f$ is a smooth function assumed to be in a Holder class $\mathscr S$. Then, the authors characterized the minimax risk:
		\begin{equation*}
			R_2(\mathscr S,\sigma_\tau):=\inf_{\hat{f}} \sup_{f \in \mathscr S} \mathbb{E}_f \|f- \hat{f}\|^2_2,
		\end{equation*}
		where $\mathscr S$ denotes a Holder class of function. Note that this problem reduces to the standard problem of nonparametric regression when $\xi_i$ is equal to $1$, on which there exists a large body of work \cite{Tsybakov86,KT93, arias2012oracle, maleki2012suboptimality, maleki2013anisotropic, KP92, donoho1998minimax, donoho1999wedgelets, devore2025optimalrecoverymeetsminimax}.

		Compared to our paper, we should emphasize that \cite{malekian2025speckle} has assumed that the forward operator $A$ is given by $I$. As expected, many complications in our derivations arise because of existence the forward operator in our model. Hence, we need completely different techniques (and different algorithms for obtaining upper bounds) from the ones presented in \cite{malekian2025speckle}. 
		
		\item Fixed forward models: The authors of \cite{zhou2022compressed, chen2024bagged, chen2026multilook} have studied a problem similar to the one presented in \prettyref{sec:fixedA}. However, there are several major differences between their work and ours.
		
		\begin{enumerate}
			
			\item None of these three papers establish lower bounds for the minimax risk. As will be clarified later, one of the main technical contributions of this paper has been to develop lower bounds for different regimes. For example, $m < n$ versus $m > n$, require distinct lower-bounding techniques. Studying exactly sparse signals again requires new techniques. In addition, for studying the lower bound in the singular case $m>n$ (i.e., the number of sensors $m$ is larger than the dimension $n$ of the signal $\bx_o$, which forces the data $\by_l\in \reals^m$ to lie in a significantly higher dimensional space) we introduce novel ideas involving the theory of Rao-Blackwell theorem and sufficient statistics, to achieve matching lower bounds. 

			\item In their models, \cite{zhou2022compressed, chen2024bagged, chen2026multilook} did not account for varying measurement scenarios. They also studied the ideal setting where the additive noise was set to zero, and imposed the assumption $m < n$; in fact, \cite{zhou2022compressed} further required $m = \Theta(n)$. As our proofs will demonstrate, relaxing each of these assumptions and deriving sharp bounds necessitate new technical contributions. For example, our proof strategy for the case $m > n$ is fundamentally different from that for $m < n$.

		\end{enumerate}

	\end{enumerate}

	Prior to this work, coarser high probability upper bounds for $\frac{\nm{\xvh-\xv_o}_2^2}{n}$ were obtained in \cite{zhou2022compressed} for the single-look speckle noise model where the signal class comes from a structured compression codebook, and \cite{chen2024bagged, chen2026multilook} for multi-look unvarying measurement speckle noise model where the signals are considered as images of a (bi-)Lipchitz function. In addition, the main theorems in these works assume the undersample regime $m\ll n$ and additive noise $\sigma_z=0$. This paper overcomes all these limitations. { The primary difference of the current results compared to the previous related works \cite{chen2024bagged,chen2026multilook} is that the above studied the special case where $A_1=\dots,A_L=A$, and only provided an estimation upper bound rather than completely characterizing the error rate. In comparison, the current work improves the error upper bound (Th II.13) and then establishes an estimation lower bound (Th II.14). Note that even then the estimation error for the fixed $A$ case rate is only characterized exactly in the special regime $\max(\sigma_z^4 ,m^2, n^2) \geq nL \log m$. In contrast, the major result in the current manuscript regarding varying $A_1,\dots,A_L$ matrices shows that we can exactly characterize the minimax error rate when the sensor matrices are independently chosen over multilooks. This not only provides a characterization of the error in the important case of varying multilooks, but also provides a complete analysis of the setup. This is where our current contribution significantly deviates from the existing results.}

	
	The classical subjects of study for $\mathscr S$ in relevant nonparametric estimation problems include H\"older classes \cite{KT93, maleki2012suboptimality, maleki2013anisotropic}, Soblev classes \cite{Ne85,NPT85}, and Besov classes \cite{KP92, donoho1999wedgelets, devore2025optimalrecoverymeetsminimax, donoho1998minimax}. Recently,  \cite{malekian2025speckle} studied the minimax risk under the speckle noise model
	\begin{equation*}
		y_i = f(x_i) \xi_i + \tau_i, \quad i = 1, 2, \ldots, n.
	\end{equation*}
	where $\xi_i$'s are i.i.d. $\mathcal N(0,1)$ and $\tau_i$'s are i.i.d. $\mathcal N(0,\sigma_n^2)$ random variables, $x_i=i/n, i=1,2,...,n$ are fixed design points, and $\mathscr S$ is the space of functions with uniform upper and lower bounds in a H\"older class. 
	
	Moreover, we treat both undersample ($m\le n$) and oversample ($m\ge n$) regimes and demonstrate how $m$ and $n$ determine the thresholds with respect to which the noise level $\sigma_z$ from $\zv_1,...,\zv_L$ can affect the minimax rates. This provides a complete picture to the minimax error estimation of this problem.
	
	\section{Preliminaries}\label{sec:preliminaries}
	
	In this section we summarize technical results used in this paper.

	\subsection{Results regarding the minimax risk}

	In the proofs of some our main results we will need some of the basic monotonicity properties of the minimax risk. While such results are intuitive and well-known on simpler problems such as in the estimation of the mean of a Gaussian random vector \cite{donoho1994ideal}, for completeness, we prove them for the estimation problems we discuss in this paper. Our next results suggest that increasing the number of observations $m$ makes the statistical problem easier and increasing the variance of the additive noise only makes the estimation problem harder.
	
	\begin{lemma}\label{lmm:monotonicity in m}
		$R_2(\Cc,m,n,\sigma_z)$ in \eqref{eq:minimaxrisk:varying} and $R_2^\dagger(\Cc,m,n,\sigma_z)$ in \eqref{eq:minimaxrisk:def:fixedA} are non-increasing in $m$.
	\end{lemma}
	
	\begin{lemma}\label{lmm:monotonicity in sigma}
		$R_2(\Cc,m,n,\sigma_z)$ in \eqref{eq:minimaxrisk:varying} and $R_2^\dagger(\Cc,m,n,\sigma_z)$ in \eqref{eq:minimaxrisk:def:fixedA} are non-decreasing in $\sigma_z$.
	\end{lemma}

	We will use the following version of Fano's method to obtain the minimax lower bounds.
	
	\begin{lemma} \cite[Lemma 3]{Yu1997}\label{lmm:fano}
		Let $\mathcal{P}$ be a space of probability measures such that for each $\P\in \mathcal{P}$, there is an associated parameter $\theta(\P)$ of interest. Let $d$ be a pseudo-metric on the space $\theta(\mathcal{P})$. Suppose there exists an integer $r\ge 2$ and parameters $\alpha_r$ and $\beta_r$ satisfying
		\begin{enumerate}
			\item $\cb{\theta(\P_1),...,\theta(\P_r)}$ is an $\alpha_r$-separated subset in $(\theta(\mathcal{P}),d)$, namely for all $1\le i\ne j\le r$,
			\begin{align*}
				d(\theta(\P_i),\theta(\P_j))\ge \alpha_r.
			\end{align*}
			
			\item For all $1\le i\ne j\le r$, we have the upper bound for Kullback–Leibler divergence
			\begin{align*}
				\KL(\P_i \parallel
				\P_j):=\int \log(\P_i/\P_j) d\P_i \le \beta_r
			\end{align*}
		\end{enumerate}
		
		Then for any $\widehat \theta \in \theta(\mathcal P)$, we have the lower bound
		$
			\max_{1\le j \le r}\E_{\P_j}\br{d\pr{\widehat \theta,\theta(\P_j)}}\ge \frac{\alpha_r}{2}\pr{1-\frac{\beta_r+\log 2}{\log r}}.
		$
	\end{lemma}

	\subsection{Results on covering numbers}
	We use the following covering number bounds.

		%



	\begin{proposition}{\cite[Proposition 4.2.12]{vershynin2018high}}\label{prop:estimate of covering numbers}
		For any Euclidean ball $B_R\subset \R^n$ of radius $R>0$, its $\de$-covering number $N_{\de}(B_R)$ (in $\ell_2$ norm) satisfies:
		\begin{equation*}
			\pr{R \over \de}^n \le N_{\de}(B_R) \le \pr{{2R \over \de} + 1}^n.
		\end{equation*}
	\end{proposition}

	\subsection{Matrix concentration}
	We extensively use the following result on random matrices.
	\begin{lemma}\cite[Theorem 2.6]{RudelsonVershinin2010}\label{lmm:singvalues} 
		Let $A$ be an $m \times n$ random matrix with elements drawn i.i.d. from $N(0,1)$. Then $\Prob(\sigma_{\max}(A) \leq \sqrt{n}+\sqrt{m}+ t) \geq 1-  2 {\rm e}^{-\frac{t^2}{2}},
		t>0$.
		Moreover, if $m<n$, then for any $t>0$,
		\begin{align*}
			\Prob(\sqrt{n}-\sqrt{m}- t \leq \sigma_{\min} (A) \leq \sigma_{\max}(A) \leq \sqrt{n}+\sqrt{m}+ t) \geq 1-  2 {\rm e}^{-\frac{t^2}{2}}.
		\end{align*}
	\end{lemma}

	Throughout the proof for the varying forward operators, we will use this result in the following way.  
	For i.i.d. $m\times n$, $N(0,1)$ Gaussian matrices $A_1,...,A_L$, we define the event 
	\begin{align}
		\Ec_{\rm maxsing}(t,L):=\bigcap_{l=1}^L \cb{\sigma_{\max}(A_l) \leq \sqrt{n}+\sqrt{m}+ t}.
	\end{align}
	With a slight overloading of our notation we also define:
	\begin{equation*}\label{eq:emaxsing}
		\Ec_{\rm maxsing}:=\Ec_{\rm maxsing}\pr{{\sqrt{n}+\sqrt{m}\over 2},L}.
	\end{equation*}
	If $n\ge 4m$, we define the event 
	\begin{align*}
		\Ec_{\rm sing}(t,L)=\bigcap_{l=1}^L \cb{\sqrt{n}-\sqrt{m}- t \leq \sigma_{\min} (A_l) \leq \sigma_{\max}(A_l) \leq \sqrt{n}+\sqrt{m}+ t},
	\end{align*}
	and again define a slightly overloaded notation:
	\begin{align}\label{eq:esing}
		\Ec_{\rm sing}:=\Ec_{\rm sing}\pr{{\sqrt{n}-\sqrt{m}\over 2},L},
		\quad
		\PP[\esing]
		\geq 1- \rvprobunion,
	\end{align}
	where the last inequality followed from \prettyref{lmm:singvalues}. Similarly, for $m\ge 4n$, define the event
	\begin{align}\label{eq:E-sing}
		\Ec_{\rm sing}'(t,L)=\bigcap_{l=1}^L \cb{\sqrt{m}-\sqrt{n}- t \leq \sigma_{\min} (A_l^\top) \leq \sigma_{\max}(A_l^\top) \leq \sqrt{m}+\sqrt{n}+ t}.
	\end{align}
	
	For the special case when $t={\sqrt{m}-\sqrt{n}\over 2}$, we again use \prettyref{lmm:singvalues} to get
	\begin{align}\label{eq: definition of esing prime}
		\Ec_{\rm sing}':=\Ec_{\rm sing}'\pr{{\sqrt{m}-\sqrt{n}\over 2},L}, \quad 
		\PP[\esingp]
		\geq 1-\rvprobunion.
	\end{align}
	
	The following lemma, proved in the appendix, is a generalization and more accurate version of Lemma 4 and 5 from \cite{zhou2024corrections}:
	
	\begin{lemma}\label{lmm:decoupling-genr}
		Let $\{A_1\}_{l=1}^L\in \reals^{m\times n}$ be Gaussian matrices. For $\bd \in \R^n$, let $D = \diag (\bd)$. Define $\widetilde \Ec_{\rm maxsing}:=\bigcap_{l=1}^L \bigcap_{i=1}^m \cb{\smax(\tilde{A}_{l,\backslash i})\le \frac{3}{2}(\sqrt{m}+\sqrt{n})}$ where $\tilde{A}_{l,\backslash i}$ is an independent copy of $A_l$ with $i$-th row removed. Then $\P(\widetilde \Ec_{\rm maxsing})\ge 1-2mL e^{-cn}$, and we have
		\begin{enumerate}
			\item The upper tail probability: For any fixed $\bd \in \R^n$
			\begin{equation*}
				\begin{split}
					&\Prob \pr{ \br{\sum_{l=1}^L\|A_l DA_l^\top \|^2_{\rm HS} > Lm\pr{\Tr(D)+t_1}^2+Lm(m-1)\|\bd\|_{2}^2 + t_2}\cap \widetilde \Ec_{\rm maxsing}}\\
					&\leq 2mL\exp\left(-c\min\left({t_1^2\over K^4\|\bd\|_{2}^2 },{t_1 \over K^2\|\bd\|_{\infty} }\right)\right)
					+
					2 m  \exp\left(-c\min\left({t_2^2\over K^4L m^3 \|\bd^2\|_2^2 },{t_2  \over K^2m \|\bd^2\|_\infty }\right)\right)\\
					&+ 2C  \exp \left( -c \min \pr{\frac{4t_2^2}{81C^2 K^4  \|\bd^2\|_{2}^2  m L ( \sqrt{n} + \sqrt{m})^4 }, \frac{2t_2}{9C K^2\| \bd\|_{\infty}^{2} ( \sqrt{n} + \sqrt{m})^2} }\right).
				\end{split}
			\end{equation*}
			
			\item The lower tail probability: For any fixed $\bd \in \R^n$
			\begin{equation*} 
				\begin{split}
					&\Prob \pr{ \br{\sum_{l=1}^L\|A_l DA_l^\top \|^2_{\rm HS} < Lm(m-1) \|\bd\|_2^2 - t}\cap \widetilde \Ec_{\rm maxsing}}\\
					&\leq  2C  \exp \left( -c \min \pr{\frac{4t^2}{81C^2 K^4  \|\bd^2\|_{2}^2  m L ( \sqrt{n} + \sqrt{m})^4 }, \frac{2t}{9C K^2 \| \bd\|_{\infty}^{2} ( \sqrt{n} + \sqrt{m})^2} }\right)  \\
					&+ 2 m  \exp\left(-c\min\left({t^2\over K^4L m^3 \|\bd^2\|_2^2 },{t  \over K^2m \|\bd^2\|_\infty }\right)\right).
				\end{split}
			\end{equation*}
		\end{enumerate}
		where $C$ and $c$ are absolute constants. Here $1\le K \le 2$ denotes the subgaussian norm of $Z\sim N(0,1)$.
	\end{lemma}

	\subsection{Linear algebraic results}
	We use the following result that is a generalization of \cite[Lemma 5]{zhou2022compressed}. 
	
	\begin{lemma}\label{lmm:bound-trace}
		Denote $\widetilde\Sigma=(\sigma_z^2I_m+A\widetilde X^2A^\top)^{-1}, \Sigma = (\sigma_z^2I_m+A X^2A^\top)^{-1}$. Then,
		\begin{align*}
			&\nmhs{A(\widetilde{X}^2-X^2)A^\top}^2{\left(\sigma_z^2+x_{\min}^2 \lmin(AA^\top) \right)^2 \over \left(\sigma_z^2+x_{\max}^2 \lmax(AA^\top)\right)^4}
			\leq
			{\trbr{(\Sigma^{-1}(\widetilde\Sigma-\Sigma) \Sigma^{-1}(\widetilde\Sigma-\Sigma)}}
			\nonumber\\
			&
			\leq
			{\left(\sigma_z^2+x_{\max}^2 \lmax(AA^\top) \right)^2 \over \left(\sigma_z^2+x_{\min}^2 \lmin(AA^\top)\right)^4}\nmhs{A(\widetilde{X}^2-X^2)A^\top}^2. 
		\end{align*}
	\end{lemma}
	
	\section{Proof of the lower bound in \prettyref{thm:lower_bnd} in the case $m< \frac n4$}
	\label{sec:proof-of-lower-bound-undersample}
	
	\subsection{Outline of the proof strategy}
	
	As increasing $a$ and $b$ only makes $\mathcal{F}_{a,b,k,n}$ larger, it suffices to prove \prettyref{thm:lower_bnd} for $a=\xmax-\xmin$ and $b=1$. We will apply \prettyref{lmm:fano} with the following definitions to obtain the minimax lower bound in \prettyref{thm:main-genr}.
	Given $\xv \in \R^n$, let $\P_{\xv}$ denote the probability distribution of data $\ly=[\yv_1^\top, \ldots, \yv_L^\top]^\top$ generated according to the model \eqref{eq:model} with $X = \diag(\xv)$. More specifically, we choose:
	\begin{equation}
		\label{eq:prob-model}
		\begin{gathered}
			\P_{\xv} \sim \otimes_{l=1}^L N(\boldsymbol{0},\Sigma_l^{{-1}}(\xv)) = N(\boldsymbol{0},\Sig^{{-1}}(\xv)),
			\\
			\Sig(\xv):=\diag\pr{\Sigma_1(\xv),...,\Sigma_L(\xv)},
			\quad
			\Sigma_l=\Sigma_l(\xv):=(\sigma_z^2 I_m + A_l X^2 A_l^\top )^{-1}.
		\end{gathered}
	\end{equation}
	Then, the parameter $\theta(\PP(\bx))$ corresponding to the distribution $\PP(\bx)$ is chosen as
	\begin{equation*}
		\theta (\P (\xv)) = \xv,
		\quad
		d(\theta(\P_{\xv}),\theta(\P_{\xv'}))=d(\xv,\xv'):=\|\xv-\xv'\|_2,
		\xv,\xv'\in \R^n,
	\end{equation*} 
	and we select a large enough parameter set $\calC$ to provide us with the desired complexity. In order to apply \prettyref{lmm:fano} we identify a discretization $\calS_{\rm sep}\subset\Cc,|\calS_{\rm sep}|= r$ satisfying: 
	\begin{enumerate}[label=(P\arabic*)]
		\item \label{clm:alpha_separation} For any $\xv_i, \xv_j \in \calS_{\rm sep}$ with $\xv_i \neq \xv_j$, the corresponding parameters $\theta(\PP_{\bx_i}), \theta(\PP_{\bx_j})$ (which are identically $\bx_i,\bx_j$ according to our definition) are well separated. In particular, for an $\alpha_r$ to be chosen appropriately, we will establish that $$d(\theta(\PP_{\xv_i}),\theta(\PP_{\xv_j})) = \|\bx_i-\bx_j\|\ge \alpha_r,
		\quad 
		\bx_i\neq \bx_j\in \calS_{\rm sep}.$$
		
		\item  For all $1\le i\ne j\le r$, for a $\beta_r>0$ to be chosen later, the distributions $\PP_{\bx_i},\PP_{\bx_j}$ are difficult to distinguish in the Kullback-Leibler divergence, at the level $\beta_r$
		\begin{equation*}
			\KL(\P_{\xv_i} \parallel
			\P_{\xv_j}) \le \beta_r, \quad \bx_i,\bx_j\in \calS_{\rm sep}. 
		\end{equation*}
	\end{enumerate}
	Then, \prettyref{lmm:fano} directly implies that there exists a constant $C>0$ for which
	\begin{align}
		\inf_{\xvh} \sup_{1\le j\le r} \E\br{\frac{\|\xvh-\xv_j\|^2}{n}}
		\geq {  \frac{C\alpha_r^2}{n}\pr{1-\frac{\beta_r+\log 2}{\log r}}^2 }. 
	\end{align} 
	In particular, for some small $c>0$, we will end up making the following choice for $r, \alpha_r,\beta_r$
	$$
	\log r = \Theta(k\log n),
	\quad
	\alpha_r
	= \Theta_{\e,x_{\max},x_{\min}}\pr{\frac{\max(\sigma_z^4,n^2)k\log n}{m^2 L}},
	\quad
	\beta_r=c\log r.
	$$
	As a consequence, by \prettyref{lmm:fano} the following lower bounds holds for any estimator $\xvh$
	\begin{align*}
		\begin{split}
			\max_{1\le i \le r} \E\br{\frac{\|\xvh-\xv_i\|^2}{n}}
			\ge \frac{\alpha_r^2}{4n}\pr{1-\frac{\beta_r+\log 2}{\log r}}^2=\Theta_{\e,x_{\max},x_{\min}}\pr{\frac{\max(\sigma_z^4,n^2)k\log n}{m^2 n L}}.
		\end{split}
	\end{align*}
	Our proofs in the following sections will use a construction which will dictate, with a high probability, the above choices for $r,\alpha_r,\beta_r$. All the proofs are provided in the appendix.
	
	\subsection{Construction of the signal class}
	
	We will first construct the signal class $\calC$. Fix $0<x_{\min}<x_{\max}$. For any array of nonnegative integers $0=a_0<a_1<\cdots <a_k=n$, let $\cb{(a_{l-1},a_l]_{\Z}:l\in [k]}$ be an ordered $k$-partition of $[n]$.
	\begin{align*}
		[n]=\cup_{l=1}^k (a_{l-1},a_l]_{\Z}.
	\end{align*}
	Let $\mathcal F(a_0,...,a_k)$ denote the set of functions from $[n]$ to $[x_{\min},x_{\max}]$ that are constant on each integer interval $(a_{l-1},a_l]_{\Z}, l=1,2,...,l$. Define
	\begin{align*}
		\mathcal F_k:=\bigcup_{0=a_0<\cdots <a_k=n} \mathcal F(a_0,...,a_k),
		\ 
		\mathcal{X}_k:=\{(f(1),...,f(n)) \in \R^n: f\in \mathcal F_k\} \subset [x_{\min},x_{\max}]^n.
	\end{align*}
	$\mathcal{X}_k$ satisfies the polynomial complexity of order $k$ defined in \prettyref{def:polcomplex}. To see this, first note that
	$\mathcal{X}_k$ can be written as a finite union:
	\begin{align} \mathcal{X}_k=\bigcup_{0=a_0<a_1<\cdots <a_k=n} \{(f(1),...,f(n)) \in \R^n: f\in \mathcal F(a_0,...,a_k)\}.
	\end{align}
	We have
	\begin{align*}
		N_{\e}(\mathcal{X}_k)\le & \sum_{0=a_0<\cdots <a_k=n}\pr{\frac{\xmax-\xmin}{\e}}^k={n\choose k}\pr{\frac{\xmax-\xmin}{\e}}^k
		\overset{(a)}{\le} \pr{\frac{n(\xmax-\xmin)}{\e}}^k
	\end{align*}
	where for $(a)$ we used ${n\choose k}\le n^k$. Hence $\mathcal{X}_k$ satisfies \eqref{eq:degree k polynomial growth}. Intuitively, we have considered a `piecewise constant' class of signals. This is a natural and popular choice for class of images in image processing \cite{rudin1992nonlinear, jalali2016compression, donoho1999wedgelets}.  \prettyref{exam:piececonstant_1} shows that this set satisfies the polynomial complexity of order $k$.  
	
	We now pick the \emph{signal class} $\Cc$ as any subset $[x_{\min},x_{\max}]^n$ which is a superset of $\mathcal X_k$ and satisfies the polynomial complexity of order $k$. Hence $
	\mathcal X_k \subset \Cc \subset [x_{\min},x_{\max}]^n$.
	Given this choice of $\Cc$ we now would like to show that 
	\begin{align*}
		R_2(\Cc,m,n,\sigma_z)
		=&\Omega_{\e,x_{\max},x_{\min}}\pr{\frac{\max(\sigma_z^4,m^2,n^2)k\log n}{m^2 n L}}.
	\end{align*}

	\subsection{Discretization of the signal class to apply Fano's Lemma}
	\label{sec:discretization}

	To construct an $\alpha_r$-separated subset $\calS_{\rm sep}$, we first define $\mathcal X^{ \rm finite} \subset \Cc$ as follows. For $\e\in (0,1)$, denote
	\begin{equation}\label{eq:Ndiv:def}
		N_{\rm div}:=kn^{\e}. 
	\end{equation}
	For simplicity, assume that $N_{\rm div},n/N_{\rm div}$ are integers. We partition $[n]$ into $N_{\rm div}$ pieces as 
	\begin{equation*}
		[n]=\bigcup_{l=1}^{N_{\rm div}}\left(\frac{(l-1)n}{N_{\rm div}},\frac{ln}{N_{\rm div}}\right]_{\Z}
	\end{equation*}
	Note that each subscripts $j\in [n]$ of the coordinates of $\xv$ will fall into one of the intervals $\left(\frac{(l_j-1)n}{N_{\rm div}},\frac{l_j n}{N_{\rm div}}\right]_{\Z}$. Each integer interval contains $\Theta(n/N_{\rm div})$ indices of $\xv = (x_1, x_2, \ldots, x_n)$. 
	Fix $0<\de_r<\frac{x_{\max}-x_{\min}}{2}$ to be determined later and define $\Bar{x}=\frac{x_{\min}+x_{\max}}{2}$. Let $\mathcal B_{N_{\rm div}}$ denote the collection of all functions from $[n]$ to $\cb{\Bar{x},\Bar{x}+\de_r}$, that are piecewise constant on each $\left(\frac{(l-1)n}{N_{\rm div}},\frac{ln}{N_{\rm div}}\right]_{\Z}, l=1,2,...,N_{\rm div}$. We define
	\begin{align}\label{eq:m8}
		\mathcal B^{\rm finite}:=\cb{(f(1),...,f(n)): f \in \mathcal B_{N_{\rm div}}} \cap \Cc \subset \cb{\Bar{x},\Bar{x}+\de_r}^n.
	\end{align}
	We construct $\mathcal{X}^{\mathrm{finite}} \subset \mathcal{B}^{\mathrm{finite}}$ as follows: among the $N_{\mathrm{div}}$ intervals, select $k/2$ of them; for the entries of the vector corresponding to these intervals, assign the value $\bar{x} + \delta_r$, and set all remaining entries to $\bar{x}$. It is straightforward to see that 
	\begin{align}\label{lower bound on X finite}
		|\mathcal X^{\rm finite}| = {N_{\rm div} \choose k/2}.
	\end{align}
	
	Now we construct a set $\calS_{\rm sep}\subset\mathcal X^{\rm finite}$ which we choose as the hypothesis for applying Fano's Theorem, and satisfies the $\alpha_r$-separation condition \ref{clm:alpha_separation}. Set $k':=k/4$ and let $\calS_{\rm sep}$ denote the set of all vectors in $\mathcal{X}^{\rm finite}$ with the following property: If  $\xv_i$ and $\xv_{j}$ are in $\calS_{\rm sep}$ then, their $q$-th components satisfy $x_{i,q} \ne x_{j,q}$ for all $q$ in at least $k'$-many different intervals of the form $\left(\frac{(l-1)n}{N_{\rm div}},\frac{l n}{N_{\rm div}}\right]_{\Z}$. The following guarantees hold. \\
	
	\noindent \textbf{Cardinality of $\calS_{\rm sep}$:} The following lemma obtains an upper and a lower bound for $|\calS_{\rm sep}|$. 
	\begin{lemma}\label{lmm:r-bound} Let $r:=|\calS_{\rm sep}|$ denote the cardinality of the set $\calS_{\rm sep}$. Then $\pr{cN_{\rm div}\over k}^{c'k}\le r \le \pr{Cn \over k}^{C'k}$
		for some absolute constants $c, C, c', C'>0$. Consequently, $\log r= \Theta(k\log n)$. 
	\end{lemma}
	\noindent \textbf{Minimum separation of elements in $\calS_{\rm sep}$:} As the number of integer points in each interval $\left(\frac{(l-1)n}{N_{\rm div}},\frac{ln}{N_{\rm div}}\right]_{\Z}$ is bounded below by $\frac{n}{N_{\rm div}}-2$, we have $\min_{\substack{\xv_i,\xv_j\in \mathcal S\\\xv_i\ne \xv_j}}\nm{\xv_i-\xv_j}_2^2 \ge  k'\cdot \pr{\frac{n}{N_{\rm div}}-2} \cdot \delta_r^2$. In view of the above, for a small constant $c>0$, we choose $\alpha_r$ as 
	\begin{align}\label{the logic of choosing alpha}
		\alpha_r^2 = \frac{ckn\de_r^2}{N_{\rm div}}.
	\end{align} 
	
	\noindent \textbf{Uniform signal strengths for elements in $\calS_{\rm sep}$:} Let $\xv_i, \xv_j \in \calS_{\rm sep}$. Suppose that $X_i, X_j$ denote diagonal square matrices $\diag(\xv_i), \diag(\xv_j)$ respectively. Since we assumed $\frac{n}{N_{\rm div}}$ is an integer, and $\xv_i,\xv_j$ have the same number of components equal to $\bar x$ and $\bar x +\de_r$, we get
	\begin{equation}\label{eq:traceCondition}
		{\rm Tr}(X_i^2-X_j^2)=0. 
	\end{equation}
	

	\subsection{The bound for Kullback-Leibler divergence}\label{sec:The bound for Kullback-Leibler divergence}
	We will use the following result.
	
	\begin{lemma}\label{lmm: upper bound for KL divergence}
		Denote $E_{\max}:= \max_{1 \leq l \leq L} \lambda_{\max} (A_lA_l^T), E_{\min}:= \min_{1 \leq l \leq L} \lambda_{\min} (A_lA_l^T)$. On the event $\Ec_{\rm sing}$ in \eqref{eq:esing}, if $\frac{\br{\sigma_z^2+x_{\max}^2\cdot\emax}\cdot \emax}{\pr{\sigma_z^2+x_{\min}^2\cdot E_{\min}}^2}\cdot x_{\max}\delta_r<\frac 14$, we have for all $\bx_i\neq \bx_j\in \calS_{\rm sep}$
		\begin{align*}
			\KL(\P_{\xv_i}\parallel \P_{\xv_j})\le & 2{\left(\sigma_z^2+x_{\max}^2 \emax \right)^2 \over \left(\sigma_z^2+x_{\min}^2 \emin \right)^4} 
			\sum_{l=1}^L
			\nmhs{A_l(X_i^2-X_j^2)A_l^\top }^2,
		\end{align*}
		where $X_i$ and $X_j$ are diagonal matrices corresponding to the vectors $\xv_i, \xv_j \in \calS_{\rm sep}$.
	\end{lemma}
	We now apply the upper tail bound of \prettyref{lmm:decoupling-genr} to find a high-probability upper bound for 
    $$\sum_{l=1}^L \nmhs{A_l(X_i^2-X_j^2)A_l^\top }^2.$$
	
	\begin{lemma}\label{lmm:bound1}
		There is a constant $C>0$ and an event $\calG$ with $\PP[\calG]\geq 1-\frac 1{(rmL)^{8}}-2mLe^{-cn}$, such that, for any $1\leq i<j\leq r$ we have on $\calG$
		\begin{align*}
			\sum_{l=1}^L \nmhs{A_l(X_i^2-X_j^2)A_l^\top }^2
            &{\le} C \Biggl\{L m x_{\max}^2 \frac{kn}{N_{\rm div}}\de_r^2\log(mLr^2)+Lm \de_r^2\log^2(mLr^2)
			+ 4Lm(m-1)x_{\max}^2\frac{kn}{N_{\rm div}}\de_r^2\\
			&\quad +\log m\pr{x_{\max}^2 \sqrt{\frac{kn}{N_{\rm div}}}\de_r^2 \sqrt{mL(\sqrt m+\sqrt n)^4\log r^2}+\de_r^2 (\sqrt n +\sqrt m)^2 \log r^2 }\Biggr\}.
		\end{align*}
	\end{lemma}
	On the event $\calG \cap \Ec_{\rm sing}$, using\prettyref{lmm: upper bound for KL divergence} and $\emax\approx n$, for a constant $\bar C:=C_{x_{\min},x_{\max}}>0$ we have
	\begin{align}\label{eq:bound-beta}
		&~\beta_r:=\max_{1\le i< j\le r} \KL(\P_i\parallel \P_j)
		\nonumber\\
		&\le 2\frac{\pr{\sigma_z^2+x_{\max}^2\cdot\emax}^2}{\pr{\sigma_z^2+x_{\min}^2\cdot \frac{1}{4}(\sqrt{n}-\sqrt{m})^2}^4}\max_{1\le i< j\le r}\sum_{l=1}^L \nmhs{A_l(X_i^2-X_j^2)A_l^\top }^2
		\nonumber\\ 
		&\le {\bar C\over \max(\sigma_z^4,n^2)}\Bigg({L m  \frac{kn}{N_{\rm div}}\de_r^2\log(mLr^2)}+Lm \de_r^2\log^2(mLr^2) +{Lm(m-1)\frac{kn}{N_{\rm div}}\de_r^2}
		\nonumber \\
		& \qquad \quad+\log m { \sqrt{\frac{kn}{N_{\rm div}}}\de_r^2 \sqrt{mL(\sqrt m+\sqrt n)^4\log r^2}}+ {\de_r^2 (\sqrt n +\sqrt m)^2 (\log m)(\log r^2) }\Bigg)
		\nonumber\\
		&
		\leq {\bar C\delta_r^2m^2nLk\over \max(\sigma_z^4,n^2)N_{\rm div}}
		\Bigg({\log(mLr^2)\over m} + {\log^2(mLr^2)N_{\rm div}\over mnk}
		+1+{ \sqrt{\frac{n N_{\rm div}\log r^2}{Lkm^3/(\log m)^2}}}+ {(\log m)(\log r^2)N_{\rm div}\over Lm^2}\Bigg) 
		\nonumber\\
		&
		\leq
		\Theta_{x_{\min},x_{\max}}(1){m^2nLk\over \max(\sigma_z^4,n^2)N_{\rm div}}\delta_r^2,
	\end{align}
	where the last inequality followed by using the following inequalities that are consequences of \prettyref{lmm:r-bound}, alongside our assumptions $\log m =\Theta(\log n),\log L=O(\log n)$, and there exists $\e\in (0,1/2)$ such that $k \le n^{1-2\e}, \max(\sigma_z^4,m^2,n^2)k\log n \le m^2n^{1-\e}L$.
	\begin{itemize}
		\item $\log(mLr^2)\leq \log m + \log L+ 2\log r<m/3$ for all large $m,L,n$ as $\log m =\Theta(\log n),\log L=O(\log n)$ and $\log r=\Theta(\log n)$ from \prettyref{lmm:r-bound}.
		\item Similarly, we have ${\log^2(mLr^2)N_{\rm div}\over mnk}\leq 1$ for all large $m,L,n$, as ${N_{\rm div}\over nk}<n^{-(1-\epsilon)}$ from \eqref{eq:Ndiv:def}.
		\item As ${N_{\rm div}\over k}\leq {n^{\epsilon}},\frac{n N_{\rm div}\log r^2}{Lkm^3/(\log m)^2}\leq 1$, using $\max(\sigma_z^4,m^2,n^2)k\log n \le m^2n^{1-\e}L$ we get
		$$
		\frac{n N_{\rm div}\log r^2}{Lkm^3/(\log m)^2}
		\leq {n^{1+\e}\log r^2\over Lm^3/(\log m)^2}\lesssim
		{\max(\sigma_z^4,m^2,n^2)\log n\over n^{1-\e}Lm^{2.5}}
		\lesssim {1\over k\sqrt m}.
		$$
		\item Finally ${(\log m)(\log r^2)N_{\rm div}\over Lm^2}\leq 1$ for all large $m,L,k$ as $N_{\rm div}=kn^\varepsilon$ and $m^2L>kn^{1+\e}$ from the last property $\max(\sigma_z^4,m^2,n^2)k\log n \le m^2n^{1-\e}L$ as well.
	\end{itemize}
	%
	In view of the above argument, and \prettyref{lmm: upper bound for KL divergence}, our we choose $\delta_r$ satisfying the following.
	\begin{itemize}
		\item We want $\beta_r<c\log r$ for a small $c>0$ to maximize the lower bound from \prettyref{lmm:fano}.
		\item To justify the application of \prettyref{lmm: upper bound for KL divergence} we need $\frac{\br{\sigma_z^2+x_{\max}^2\cdot\emax}\cdot \emax}{\pr{\sigma_z^2+x_{\min}^2\cdot E_{\min}}^2}\cdot x_{\max}\delta_r<\frac 14$. 
		
		\item Since the maximum value the signals in $\calS_{\rm sep} \subset \mathcal{C}$ can take is $\bar{x} + \delta_r$ and the maximum value of the entries of the vectors in $\mathcal{C}$ is $x_{\max}$ we have to ensure that $\bar{x}+ \delta_r < x_{\max}$. 
	\end{itemize}
	In view of the above, for a small constant $c_\delta:=c_\delta(x_{\min},x_{\max})>0$, we choose 
	\begin{equation}\label{eq:choice of delta_r}
		\de_r^2:=c_{\de}\Delta^{-2}\cdot \frac{\max(\sigma_z^4,n^2)\log r}{nm^2L}\cdot \frac{N_{\rm div}}{k}
	\end{equation}
	where $r = |\calS_{\rm sep}|$, and  $\Delta=\Delta(m,n,x_{\max},x_{\min})
	:=\frac{\br{\sigma_z^2+x_{\max}^2\cdot\emax}\cdot \emax}{\pr{\sigma_z^2+x_{\min}^2\cdot E_{\min}}^2}\cdot 2x_{\max}.
	$
	The first two conditions can be verified using the upper bound \eqref{eq:bound-beta} and the definitions of $\delta_r^2$. To check the final condition, i.e., $\bar{x}+ \delta_r < x_{\max}$, we first note that since $E_{\max}$ and $E_{\min}$ depend on $A_l$, $\delta_r$ is a random number. Hence, to find a deterministic upper bound for $\delta_r$ we only consider the measurement matrices that belong to the event $\Ec_{\rm sing}$ defined in \eqref{eq:esing}. 
	In view of the definition of $\Ec_{\rm sing}$, in the above event we have $E_{\max}= \Theta(n)$ and $E_{\min} = \Theta (n)$. Hence,
	\begin{align}\label{eq: simplified delta r}
		\de_r^2=\Theta_{x_{\max},x_{\min}}\pr{\frac{\max(\sigma_z^4,n^2)\log r}{nm^2L}\cdot \frac{N_{\rm div}}{k}}
		=\Theta_{x_{\max},x_{\min}}\pr{\frac{\max(\sigma_z^4,n^2)\log r}{n^{1-\e}m^2L}},
	\end{align}
	where the last identity followed by plugging \eqref{eq:Ndiv:def} and the result of \prettyref{lmm:r-bound}
	By the assumptions we made in the main theorem, i.e., $\max(\sigma_z^4,m^2,n^2)k\log n \le m^2n^{1-\e}L$ and the fact that $\max(\sigma_z^4,m^2,n^2)=n^2$, we can see that $\pr{\frac{n^{1+ \epsilon}  k \log (n/k) }{m^2 L}}<1$. Hence, by choosing $c_\delta$ in the definition of $\delta_r$ in \eqref{eq:choice of delta_r} small enough we can ensure that $\bar{x}+\delta_r< x_{\max}$ holds.

	\subsection{Concluding the proof}
	
	For a small $c_{\de}>0$, \eqref{eq:bound-beta} ensures that we have $\beta_r\le \frac{1}{10}\log r$. Therefore, conditioned on $A_1,...,A_L$ being in the high probability event $\calG\cap \Ec_{\rm sing}$ as defined in \prettyref{lmm:bound1}, using \eqref{the logic of choosing alpha}, \eqref{eq:choice of delta_r} we get 
	\begin{equation}
		\label{eq:alpha-choice}
		\alpha_r^2 = \frac{ckn\de_r^2}{N_{\rm div}}
		=\Theta_{x_{\max},x_{\min}}\pr{\frac{\max(\sigma_z^4,n^2)}{m^2} \frac{k\log \pr{N_{\rm div}/k}}{L}}.
	\end{equation}
	As a consequence, by \prettyref{lmm:fano} we have for any estimator $\xvh$, the lower bound
	\begin{equation*}
		\begin{split}
			&\max_{1\le i \le r} \E\br{\frac{\|\xvh-\xv_i\|^2}{n}}
			\ge \frac{\alpha_r^2}{4n}\pr{1-\frac{\beta_r+\log 2}{\log r}}^2 \P\pr{\widetilde\Ec_{\rm maxsing}\cap \Ec_{\rm dcpl} \cap \Ec_{\rm sing}}
			=\Theta\pr{\frac{\alpha_r^2}{n}}\\
			=& \Theta_{x_{\max},x_{\min}}\pr{\frac{\max(\sigma_z^4,n^2)}{m^2 n} \frac{k\log \pr{N_{\rm div}/k}}{L}}
			=\Theta_{\e,x_{\max},x_{\min}}\pr{\frac{\max(\sigma_z^4,n^2)k\log n}{m^2 n L}}.
		\end{split}
	\end{equation*}


	\section{Proof of the lower bound in \prettyref{thm:lower_bnd} for the case $m\geq \frac n4$}\label{sec:proof of lower bound oversample}
	
	This section will primarily establish the lower bound for the sub-case $m\geq 4n,\sigma_z^2=0$ given by 
	\begin{equation}\label{eq:m12}
		R_2(\Cc_k,m,n,0)
		=\inf_{\bde}\sup_{\xv\in \Cc_k}\E\br{\frac{\nm{\bde(\ly)-\xv_o}_2^2}{n} }
		=  \Omega_{\e,\xmin\xmax}\pr{\frac{k\log n}{nL}},
		\quad 
		m\geq 4n.
	\end{equation}
	Then, the lower bound for a general $\sigma_z^2\geq 0$ and $4n\geq m\ge \frac n4$ follow from \prettyref{lmm:monotonicity in m} and \prettyref{lmm:monotonicity in sigma}, as the minimax risk is a non-increasing function of $m$ and an increasing function of  $\sigma_z^2$, the variance of the additive noise component. To see the above, we first fix $\sigma_z^2=0$. In view of \prettyref{lmm:monotonicity in m}, as the risk is a non-increasing function of $m$, the last display implies
	\begin{equation}\label{eq:m19}
		R_2(\Cc_k,m,n,0)\ge R_2(\Cc_k,4n,n,0)=\Omega_{\e,x_{\min},x_{\max}}\pr{ \frac{k\log n}{nL}},
		\quad 
		\frac n4\leq m\leq 4n.
	\end{equation}
	Hence, combining the lower bounds in \eqref{eq:m12} and \eqref{eq:m19} we get
	\begin{align}\label{eq:section on the decoupling in the upper bound}
		R_2(\Cc_k,m,n,0)
		=\Omega_{x_{\max},x_{\min}}\pr{\frac{k\log n}{nL}},
		\quad 
		m\geq \frac n4.
	\end{align}
	To achieve a lower bound for a general $\sigma_z^2 \le m$, when $m\geq \frac n4$, we first use that the miminax error is non-decreasing function in $\sigma_z$ (\prettyref{lmm:monotonicity in sigma}) to get $R_2(\Cc_k,m,n,0)\leq R_2(\Cc_k,m,n,\sigma_z)$.
	Then, for $\sigma_z^2\le m$, combining \eqref{eq:section on the decoupling in the upper bound} with the upper bound for $\sigma_z^2=m$ from \prettyref{thm:main-genr}, for $m\geq \frac n4$, we get
	\begin{equation}
		{C_1\frac{k\log n}{nL}}\leq R_2(\Cc_k,m,n,0)\leq 
		R_2(\Cc_k,m,n,\sigma_z)
		\leq 
		R_2(\Cc_k,m,n,\sqrt{m})
		\leq {C_2\frac{k\log n}{nL}},
	\end{equation}
	where $C_1,C_2$ are constants depending on $\xmin,\xmax$. Hence, the sandwich inequality implies
	\begin{equation}\label{eq:small noise case lower bound}
		R_2(\Cc_k,m,n,\sigma_z)=\Theta_{x_{\max},x_{\min}}\pr{\frac{k\log n}{nL}}, \text{ whenever } m\geq \frac n4, \sigma_z^2\le m.
	\end{equation}
	Next, suppose $m\ge \frac n4$ and $ \sigma_z^2 \ge m$. We observe that the only use of $m\leq \frac n4$ in the proofs of \prettyref{sec:proof-of-lower-bound-undersample}, was to bound $\emin$ from below. Notably, in the proof of  \prettyref{sec:proof-of-lower-bound-undersample}, only places we used the lower bound on $\emin$, are given in \prettyref{lmm: upper bound for KL divergence} and \eqref{eq:bound-beta}. Then we note that we can repeat the entire analyses of lower bound in \prettyref{sec:proof-of-lower-bound-undersample} to establish the lower bound for the case $m\ge \frac n4$ and $ \sigma_z^2 \ge m$ by replacing $\Ec_{\rm sing}$ with $\Ec_{\rm maxsing}$ and the lower bound for the singular value $\smin(\sigma_z^2 I_n+A_lX_o^2 A_l^\top )$ by $\sigma_z^2$. As a consequence, we have the lower bound
	\begin{equation}\label{eq:large noise case lower bound}
		R_2(\Cc_k,m,n,\sigma_z)=\Omega_{x_{\max},x_{\min}}\pr{\frac{\sigma_z^4}{m^2 n} \cdot \frac{k\log n}{L}}, \text{ whenever } m\ge \frac n4, \sigma_z^2 \ge m.
	\end{equation}
	Combining \eqref{eq:large noise case lower bound} and \eqref{eq:small noise case lower bound} for the case $m\geq \frac n4$ yields the desired minimax lower bound
	\begin{align*}
		R_2(\Cc_k,m,n,\sigma_z)=\Omega_{x_{\max},x_{\min}}\pr{\frac{\max(\sigma_z^4,m^2)}{m^2 n} \cdot \frac{k\log n}{L}},
		\quad m\geq \frac n4, \sigma_z\geq 0.
	\end{align*}

	We complete our proof by establishing the lower bound in \eqref{eq:m12}. We use ideas related to sufficient statistic for our proof, using the following definition throughout this section.
		\begin{definition}\cite[Definition 6.2.1]{casella2024statistical}
			A statistics $\bT(\ly)$ is sufficient for $\xv_o$ if the distribution of sample $\ly$ conditioned on $\bT(\ly)$ is independent of $\xv_o$.
		\end{definition}
		Define $\bA=\diag(A_1,\dots,A_L)\in \reals^{mL\times nL}$ and note that $\bA^\top \bA=\diag(A_1^\top A_1,\dots,A_L^\top A_L)$. Throughout the section we analyze the expected loss on the high probability event $\esingp$ as in \eqref{eq: definition of esing prime}, where all the matrices $\{A_l^\top A_l\}_{l=1}^L$ are invertible. Then we have the following result.
		\begin{proposition}\label{prop:sufficient}
			Consider the case $\sigma_z=0$ and that the event $\esingp$ holds. Then $\bT_{\bA}(\ly)=(\bA^\top  \bA)^{-1}\bA^\top \ly$ is a sufficient statistic for the parameter $\xv_o$.
		\end{proposition}
		
		We will use the Rao-Blackwell theorem to first bound the desired minimax risk from below using the squared error loss for the sufficient statistic $\bT_{\bA}(\ly)$.
		\begin{theorem}\cite[Theorem 2.5]{shao2008mathematical}\label{thm:rao}
			Let $\bde(\ly)$ be any estimator of the parameter $\xv_o$ and $\bT_{\bA}(\ly)$ is a sufficent statistic for $\xv_o$. Then $g(\bT_{\bA}(\ly)) =\E[ \bde(\ly) \mid \bT_{\bA}(\ly)]$ is also an estimator for $\xv_o$ and it provides an improved error guarantee
			\begin{equation*}
				\EE\qth{\nm{\bde(\ly)-\xv_o}_2^2|\esingp}\ge \EE\qth{\nm{g(\bT_{\bA}(\ly))-\xv_o}_2^2|\esingp}.
			\end{equation*}
		\end{theorem}
		In view of the above result, we establish a lower bound to $\EE\qth{\nm{g(\bT_{\bA}(\ly))-\xv_o}_2^2}$ to complete our analysis. This is provided in the result below.
		
		\begin{lemma}\label{lmm:lower bound for the improved estimator}
			Consider \eqref{eq:multilook speckle noise model} with $\sigma_z=0,m\ge 4n$. Then, for a constant $C\geq 0$, we have
			$$
			\inf_{g}\sup_{\xv\in \Cc_k}\E\br{\nm{g(\bT_{\bA}(\ly))-\xv_o}_2^2|\esingp} \geq C{\frac{k\log n}{L}}.
			$$
		\end{lemma}
		
		Using the last display, \eqref{eq:m12} follows from \prettyref{thm:rao} and \prettyref{lmm:lower bound for the improved estimator} as $\P\pr{\Ec_{\rm sing}'}\geq \frac 12$
		\begin{equation*}
			\inf_{\bde}\sup_{\xv\in \Cc_k}\E\br{\frac{\nm{\bde(\ly)-\xv_o}_2^2}{n} }
			\ge \inf_{\bde}\sup_{\xv\in \Cc_k}\E\br{\left.\frac{\nm{g(\bT_{\bA}(\ly))-\xv_o}_2^2}{n} \right \vert \esingp}
			=  \Omega\pr{\frac{k\log n}{nL}}.
		\end{equation*}

		\section{Proof of the upper bound in \prettyref{thm:main-genr}}
		\label{sec:upper bound}
		
		\subsection{General strategy}\label{sec:General strategy}
		
		Without a loss of generality, we assume $a=\xmax-\xmin,b=1$. We provide the proof of all the technical results in the appendix. We will show that the desired upper bound is achieved by the maximum likelihood estimator
		\begin{equation}\label{eq:mle-genr}
			\xvh_o = \arg\min_{\xv \in \Cc} f(\xv),
		\end{equation}
		where $\calC$ is the class of all possible signals, and the negative log-likelihood $f(\xv)$ is defined as
		\begin{align}\label{eq:log-likelihood}
			f (\xv) 
			= \sum_{l=1}^L \log\det \left(\sigma_z^2 I_m +  A_l X^2 A_l^\top \right) + \sum_{l=1}^L \yv_l^\top  \left( \sigma_z^2 I_m+ A_l X^2 A_l^\top \right)^{-1}  \yv_l
		\end{align}
		For the entirety of the analysis in this section, we will restrict ourselves to the following event 
		
		\begin{align}\label{eq:eigen-conditioning}
			\calS=\sth{\sth{A_l}_{i=1}^L:\sigma_{\min}(A_lA_l^\top)\geq E_{\min},\quad  \sigma_{\max}(A_lA_l^\top) \leq \emax,\quad  l=1,\dots, L}.
		\end{align}
		This is the only place where we impose restrictions on the singular values of $\sth{A_l}_{l=1}^L$. To get to the specific minimax risk guarantees in different regimes we will choose appropriate values of $\emax,\emin$. In particular, we have the following considerations.
		\begin{itemize}
			\item Case I ($n\geq 4m$): We will choose $\emax=\frac 94(\sqrt n+\sqrt m)^2$ and $\emin = \frac 14(\sqrt n-\sqrt m)^2$. In that case the event $\calS$ satisfies $\esing\subseteq \calS$, where $\esing$ is given as in \eqref{eq:esing}. This implies $\PP\qth{\calS}\geq \PP[\esing]\geq 1-Le^{-cn}$ for some constant $c>0$.
			\item Case II ($n< 4m$): We will choose $\emax=\frac 94(\sqrt n+\sqrt m)^2$ and $\emin = 0$. In that case the event $\calS$ satisfies $\emaxsing\subseteq \calS$, where $\emaxsing$ is given as in \eqref{eq:emaxsing}. This implies $\PP\qth{\calS}\geq \PP[\emaxsing]\geq 1-Le^{-cn}$ for some constant $c>0$.
		\end{itemize}
		In other words we also have
		\begin{align}\label{eq:prob-S}
			\PP\qth{\calS}\geq 1-Le^{-cn},\quad \text{for all $L,m,n$}.
		\end{align}
		Consider the following notations for simplifying the presentation. Let $\{\Sigma_l\}_{l=1}^L$ be the collection of inverses of the covariance matrix $\EE\qth{\yv_l \yv_l^\top|A_l}$ given by
		\begin{equation*}
			\Sigma_l=\Sigma_l(\xv):=(\sigma_z^2 I_m + A_l X^2 A_l^\top )^{-1},
			\quad l=1,\dots,L.
		\end{equation*} 
		Define the vector $\ly\in \reals^{mL}$ and block-diagonal matrix $\Sig(\xv)\in \reals^{mL\times mL}$ as the collection of all the observations and the inverse covariance matrices over different looks
		\begin{align}\label{eq:So}
			\ly^\top :=(\yv_1^\top ,...,\yv_L^\top ),
			\quad \Sig(\xv):=\diag\pr{\Sigma_1(\xv),...,\Sigma_L(\xv)},
			\So=\bSigma(\bx_o),
			\Sohat = \bSigma(\hat \bx_o).
		\end{align}
		In view of the above notations, we can rewrite the negative log-likelihood in \eqref{eq:log-likelihood} as
		\begin{align}
			f(\xv)= -\log \det(\Sig(\xv))+ \ly^\top  \Sig(\xv)\ly.
		\end{align}
		Our proof strategy draws inspiration from the empirical loss minimization literature, such as \cite{fan2024factor}, \cite{fan2025factor}. To achieve a parametric error rate in the sample size $L$, we compare the negative log-likelihood for the estimator $\xvh_o$ and the true parameter $\xv_o$, and this also leads to the minimax optimal rate. Since $\xvh_o$ is the minimizer of \eqref{eq:mle-genr}, we have
		\begin{equation}\label{eq:1}
			f(\xvh_o)\leq f(\xv_o).
		\end{equation}
		For a fixed $\bx$, define $\overline f(\bx)$ as the conditional expectation of $f(\bx)$ given $A_1,\dots,A_L$,
		\begin{align}
			\overline{f}(\bx):=\E[f(\bx) \mid A_1,...,A_L]
			= {-\log\det \Sig(\xv) +{\rm Tr}\pr{\Sig(\xv) \Sig(\xv_o)^{-1}}}.
		\end{align}
		Simplifying the expression for $f(\xvh_o)-f(\xv_o)$, with the notations in \eqref{eq:So} we get
		\begin{align*}
			~f(\xvh_o)-f(\xv_o)
			&= \ly^\top  \pr{\Sohat-\So} \ly
			-\Tr \qth{\So^{-1}\pr{\Sohat-\So}}+\Tr \qth{\So^{-1}\pr{\Sohat-\So}}
			-\log \det(\Sohat\So^{-1}) 
			\nonumber\\
			&= \ly^\top\pr{\Sohat-\So}\ly - \Tr \qth{\So^{-1} \pr{\Sohat-\So}}+ \overline{f}(\hat \bx_o)-\overline{f}(\bx_o).
		\end{align*}
		Therefore, in view of \eqref{eq:1} we get
		\begin{align}\label{eq:m15}
			\ly^\top  \pr{\So-\Sohat} \ly - \Tr(\So^{-1}(\So-\Sohat))\ge \overline{f}(\hat \bx_o)-\overline{f}(\bx_o).
		\end{align}
		Our following approach is to find an upper bound for the left side in terms of $\nm{\widehat \xv_o - \xv_o}_2$ and a lower bound for the right side in terms of $\nm{\widehat \xv_o - \xv_o}_2$, and simplify the inequality to get an upper bound for $\nm{\widehat \xv_o - \xv_o}_2$. Throughout the rest of the draft we will use the notation
		\begin{align}\label{eq:cnms}
			\begin{split}
				\cnms &:= C(n,m,\sigma_z,x_{\max},x_{\min}) = c{\sigma_z^2+x^2_{\max} \emax \over \sigma_z^2+x^2_{\min} \emin },
			\end{split}
		\end{align}
		where $c>0$ is a large universal constant. Note that in the regime $n\geq 4m$, $\cnms$ is of constant order as long as $x_{\max},x_{\min}$ are of constant order. We will use expressions similar to $\cnms$ in the analysis of the case $n<4m$, the related details will be presented later. \\
		
		\noindent \textbf{Establishing an upper bound on $\ly^\top  (\So-\Sohat) \ly - \Tr(\So^{-1}(\So-\Sohat))$:} The main challenge in the analysis is the dependency of $\Sohat$ on $\ly$, which prevents us from directly applying concentration inequalities to bound $\ly^\top  (\So-\Sohat) \ly$. To resolve this issue, we use a $\delta$-net argument, as will be clarified below. Consider a $\delta$-net of the set $\calC_k$, denoted by $\netdelta$, with the choice of $\delta$ to be discussed later. Define $\xvt_o$ as the closest vector in $\netdelta$ to $\bx_o$, 
		\begin{align}\label{eq:netdelta}
			\xvt_o=\argmin_{\xv\in\netdelta}\|\xvh_o-\xv\|_2.
		\end{align}
		We will use the following notations for the rest of the section
		\begin{align}
			\label{eq:m-Sotil}
			\Sotil=\Sig({\xvt}_o),\quad \Xt_o={\rm diag}({\xvt_o}),
			\quad \Sotil=\bSigma(\xvt),\quad \Xt={\rm diag}({\xvt}),\quad \xvt\in \netdelta.
		\end{align}
		Then in view of triangle inequality we get
		\begin{align}\label{eq:m1}
			&\abs{\ly^\top  (\So-\Sohat) \ly - \Tr(\So^{-1}(\So-\Sohat))}
			\nonumber\\
			&\leq \abs{\ly^\top  (\Sotil-\So) \ly - \Tr(\So^{-1}(\Sotil-\So))}
			+\abs{\ly^\top  (\Sotil-\Sohat) \ly - \Tr(\So^{-1}(\Sotil-\Sohat))}.
		\end{align}
		We use an union bound argument to control the first term above, uniformly over all possible choices of $\xvt\in \netdelta$. This is done in the following result.
		
		\begin{lemma}\label{lmm:m-error-union-bound}
			There exist constants $c_1,c_2,c_3,c_4>0$ such that
			the following holds with probability $1-Le^{-cn}-e^{-c_1k\log ((\xmax-\xmin)n/ \delta_{\rm net})}$
			$$\abs{\ly^\top  (\Stil-\So) \ly - \Tr(\So^{-1}(\Stil-\So))}
			\leq b_1 \sqrt{\mathscr Z} +b_1',\quad
			\text{for all } \xvt_o\in \netdelta,
			$$
			where, with the notation in \eqref{eq:cnms}, $b_1,b_1',\mathscr Z$ are defined as
			\begin{align*}
				\begin{gathered}
					b_1= c_3\sqrt{k\log\pr{\frac{(\xmax-\xmin)n}{\delta_{\rm net}}}},
					\quad  b_1'= c_4
					\cnms k\log\pr{\frac{(\xmax-\xmin)n}{\delta_{\rm net}}}\cdot {x_{\max}^2},\\
					\mathscr Z= \Tr(\So^{-1}(\Sotil-\So)\So^{-1}(\Sotil-\So)).
				\end{gathered} 
			\end{align*}
		\end{lemma}
		The following result controls the final term of \eqref{eq:m1}.
		\begin{lemma}\label{lmm:m-bound-deltanet-error}
			Let $\cnms$ be as in \eqref{eq:cnms} and denote $b_2=(\cnms)^2 mL \delta_{\rm net}$. There exist constants $c_1,c_2>0$ such that the following holds with probability $1- Le^{-c_1 n} - e^{-c_2mL}$
			\begin{align*}
				&\abs{\ly^\top  (\Sotil-\Sohat) \ly - \Tr(\So^{-1} (\Sotil-\Sohat))}
				\leq b_2.
			\end{align*}
		\end{lemma}
		Combining \prettyref{lmm:m-bound-deltanet-error} with \prettyref{lmm:m-error-union-bound}, in view of \eqref{eq:m1} we have
		\begin{align}\label{eq:m16}
			\abs{\ly^\top  (\So-\Sohat) \ly - \Tr(\So^{-1}(\So-\Sohat))}
			\leq 
			b_1 \sqrt{\mathscr Z} +b_1'+b_2
		\end{align}

		\noindent \textbf{Establishing a lower bound on $\overline{f}(\xvh_o)-\overline{f}(\xv_o)$:} Our analysis uses the decomposition
		\begin{align}\label{eq:m2}
			\overline{f}(\xvh_o)-\overline{f}(\xv_o)=\overline{f}(\xvh_o)-\overline{f}(\xvt_o)+\overline{f}(\xvt_o)-\overline{f}(\xv_o),
		\end{align}
		with $\xvt_o$ as in \eqref{eq:netdelta}. The first term, $\overline f(\xvh_o)-\overline f (\xvt_o)$ can be bounded by $\cnms x_{\max} n\delta_{\rm net}$ using the fact that $\tilde \bx_o$ is chosen to be at most $\delta_{\rm net}$ distance away from $\hat\bx_o$. We bound $\overline f(\xvt_o)-\overline f (\xv_o)$ using the following result.
		\begin{lemma}\label{lmm:bound-fbar-diffs}
			Let $\sigma_z^2 I_m+A_l\widetilde{X}_o^2A_l^\top$ and $\sigma_z^2 I_m+A_l X_o^2 A_l^\top, 1\le l \le L$, be invertible. Then,
			\begin{align*}
				\overline{f}(\xvt_o)-\overline{f}(\xv_o)&\geq {1\over 2(1+\widetilde\lambda_{\max})^2}{\rm Tr}\left(\So^{-1}(\Sotil-\So)\So^{-1}(\Sotil-\So)\right),
			\end{align*}
			where $\widetilde\lambda_{\max}>0$ is the maximum singular value of $\So^{-{1\over 2}}(\Sotil-\So)\So^{-{1\over 2}}$. Moreover, $\widetilde \lambda_{\max}\leq \cnms$ on the event $\calS$ in \eqref{eq:eigen-conditioning}.
		\end{lemma}
		Note that we use these results in the subcase $n\geq 4m$, and the subcase $n<4m$ with $\sigma_z^2\geq m$. Then we have $\cnms$ is at most a constant. The following result controls $|\overline f(\xvh_o)-\overline f (\xvt_o)|$ for a given $\delta_{\rm net}$.
		\begin{lemma}
			\label{lmm:f-sotil-f-sohat-bound}
			$\PP\qth{|\overline f(\xvh_o)-\overline f (\xvt_o)|
				\leq (\cnms)^2\cdot x_{\max} mL \de_{\rm net}}\geq 1-L e^{-cn}$ for some $c>0$.
		\end{lemma}
		
		%
		Combining the above results, in view of \eqref{eq:m2} we have, with probability $1-L\exp(-cn)$, for a small enough $\delta_{\rm net}$
		\begin{equation}\label{eq:m14}
			\overline{f}(\xvh_o)-\overline{f}(\xv_o) \ge \frac{\mathscr Z}{(\cnms)^2}-1.
		\end{equation}
		
			\noindent \textbf{Simplifying the quadratic inequality:} Combining \eqref{eq:m14} and \eqref{eq:m16}, in view of \eqref{eq:m15},
	we have 
	\begin{align}\label{eq:m17}
		\PP\qth{\frac{\mathscr Z}{(\cnms)^2} \le b_1 \sqrt{\mathscr Z}+b_1'+ b_2+1}
		\geq 1-e^{-c_1k\log ((\xmax-\xmin)n/ \delta_{\rm net})}-Le^{-cn}-\exp\pr{-cmL}.
	\end{align}
	%
	Rewrite the last inequality as $az^2-bz-c\leq 0$, with $z=\sqrt{\mathscr Z}, a=\frac 1{(\cnms)^2},b=b_1, c=b_1'+b_2+1
	$. Note that by our choice of $b_1,b_1'$ and $b_2$, $\frac{b^2}{a^2}\lesssim \frac{c}{a}$.
    Note $z^2$ is no more than the square of the positive root of $az^2-bz-c = 0$, which implies
	\begin{align}\label{eq:quad-ineq-1}
		\mathscr Z 
		= z^2
		\leq \pth{b+\sqrt{b^2+4ac}\over 2a}^2
		\leq \frac{2b^2+2(b^2+4ac)}{4a^2}=\frac{b^2}{a^2}+\frac{2c}{a}\le C'\cdot  \frac{c}{a}
	\end{align}
	for some constant $C'>0$. Using the notations from, \eqref{eq:cnms}, \prettyref{lmm:m-error-union-bound} and \prettyref{lmm:m-bound-deltanet-error} we get
	\begin{align*}
		\mathscr Z
		\le C'\cdot  \frac {b_1'+b_2+1}a 
		=
		C'(\cnms)^2\pth{c_3 k\log\pr{\frac{(\xmax-\xmin)n}{\delta_{\rm net}}}\cdot {x_{\max}^2}+(\cnms)^2 mL \delta_{\rm net}+1}.
	\end{align*}
	%
	Choose $\delta_{\rm net}=\frac{x_{\max}}{n^5}$ and recall $mL\le n^4 k\log n$ from \prettyref{thm:main-genr}. Then, from the last display we use \eqref{eq:m17} to get for a constant $C>0$
	\begin{align}\label{eq:m-scrZ-upperbound}
		\PP\qth{\mathscr{Z}\leq 
			C\cdot (\cnms)^2 k\log n}
		=1-O\pr{n^{-ck}+L\exp(-cn)+2\exp\pr{-cmL}}.
	\end{align}

	\noindent \textbf{Finding a lowerbound for $\mathscr Z$:}
	In view of \prettyref{lmm:bound-trace}, using the block structure of $\bSigma_o,\Stil$ given in \eqref{eq:m-Sotil}, we have on the event $\calS$,
	\begin{align}\label{eq:Z-conjugation-relation}
		\begin{split}
			\mathscr Z
			&
			= \trbr{\So^{-1}\pr{\Sotil-\So}\So^{-1}\pr{\Sotil-\So}}\\
			&  =\sum_{l=1}^L\trbr{(\Sigma_l(\xv_o)^{-1}(\Sigma_l(\xvt_o)-\Sigma_l(\xv_o))\Sigma_l(\xv_o)^{-1}(\Sigma_l(\xvt_o)-\Sigma_l(\xv_o))} 
			\\
			&
			\geq 
			{1 \over (\cnms)^2 \left(\sigma_z^2+x_{\max}^2 \emax  \right)^2} \sum_{l=1}^L\nmhs{A_l(\widehat{X}_o^2-X_o^2)A_l^\top }^2,
		\end{split}
	\end{align}
	where $\cnms$ is as in \eqref{eq:cnms}. The lower bound on $\mathscr{Z}$ is completed with the following lower bound on $\sum_{l=1}^L\nmhs{A_l(\widehat{X}_o^2-X_o^2)A_l^\top }^2$. The proof follows from \prettyref{lmm:decoupling-genr} and is given in \prettyref{app:section on the decoupling in the upper bound}.
	\begin{lemma}\label{lmm:decoupling-lower-bound}
		 
		Let $n^4k\log n>mL$. Then the following holds with a probability $1-n^{-ck}-mLe^{-cn}$
		\begin{equation*}
			\begin{split}
				&\sum_{l=1}^L\|A_l(\widetilde{X}_o^2-X_o^2)A_l^\top \|^2_{\rm HS}\\
				&\ge 4m(m-1)L x_{\min}^2 \|\xvt_o- \xv_o\|_2^2  - \breve C_{x_{\max}} \max(m,n) \{\|\xvt_o- \xv_o\|_{2} \cdot  \sqrt{mLk \log n} + k \log n\}.
			\end{split}
		\end{equation*}
		where $\breve C_{x_{\max}}$ is a constant depending on $x_{\max}$ that grows at most polynomially with $x_{\max}$.
	\end{lemma}

	\noindent { \textbf{Final upper bound on $\|\xvt_o- \xv_o\|_2^2$:}
	We combine \eqref{eq:Z-conjugation-relation}, \eqref{eq:m-scrZ-upperbound}, and \prettyref{lmm:decoupling-lower-bound} to summarize the above in terms of the following quadratic inequality with respect to $\|\xvt_o- \xv_o\|_2$, that holds with a probability $1-O\pr{n^{-ck}+L\exp(-cn)+\exp\pr{-cmL}}$
	\begin{align*}
		\begin{gathered}
			a\|\xvt_o- \xv_o\|_2^2 -b \|\xvt_o- \xv_o\|_2
			-d\leq 0\\
			a={C_1m(m-1)Lx_{\min}^2\over (\cnms)^2(\sigma_z^2+x_{\max}^2 \emax )^2},
			\quad
			b={\breve C_{x_{\max}} \max(m,n) \sqrt{mLk\log n} \over (\cnms)^2(\sigma_z^2+x_{\max}^2 \emax )^2},
			\\
			d= {\breve C_{x_{\max}} \max(m,n) k \log n \over C\cdot (\cnms)^2(\sigma_z^2+x_{\max}^2 \emax )^2}
			+C\cdot (\cnms)^2 k\log n.
		\end{gathered}
	\end{align*}
	In view of an argument similar to \eqref{eq:quad-ineq-1}, we have with a probability $1-O\pr{n^{-ck}+L\exp(-cn)+\exp\pr{-cmL}}$
	\begin{align*}
		&\frac 1n\|\xvt_o- \xv_o\|_2^2
		\leq {2\over n}\pth{{b^2\over a^2}+{ d\over a}}\\
		&\leq
		O_{x_{\min},x_{\max}}
		\Biggl(
		{\max\sth{m^2,n^2}mLk(\log n)\over nm^4L^2}
		+{\max\sth{m,n}k\log n\over nm^2L}
		+{(\cnms)^4(\sigma_z^2+x_{\max}^2 \emax )^2k\log n \over nm^2L}
		\Biggr)\\
		&\leq
		O_{x_{\min},x_{\max}}
		\Biggl(
		{\max\sth{m^2,n^2}k\log n\over nm^2L}\pth{{1\over m}+{1\over \max\sth{m,n}}}
		+{(\cnms)^4(\sigma_z^2+x_{\max}^2 \emax )^2k\log n \over nm^2L}
		\Biggr)
	\end{align*}
	In view of $\frac 1n\|\xvt_o- \xv_o\|_2^2\leq x_{\max}^2$, and $\|\xvt_o-\xvh_o\|_2\leq \delta_{\rm  net}\leq {x_{\max}\over n^5}$ from \eqref{eq:netdelta}, the above implies
	\begin{align*}
		\EE\qth{\frac 1n\|\xvh_o- \xv_o\|_2^2}
		&\leq 
		O_{x_{\min},x_{\max}}
		&&\Biggl({o(1)\cdot \max\sth{m^2,n^2}k\log n\over nm^2L}+{(\cnms)^4(\sigma_z^2+x_{\max}^2 \emax )^2 k\log n\over nm^2L}
		\nonumber\\
		& && \quad + n^{-ck}+L\exp(-cn)+\exp\pr{-cmL}\Biggr).
	\end{align*}
	%
	%
	%
	%
	%
	%
	
	}
		
		\subsection{Proof of \prettyref{thm:main-genr}}
		
		We first consider the subcase $n\geq 4m$, and the subcase $n<4m$ with $\sigma_z^2\geq m$. Then we have $\cnms$ is of a constant order and $\emax=\Theta(m+n)$. In view of the last display, the above implies
		\begin{align*}
			\EE\qth{\frac 1n\|\xvh_o- \xv_o\|_2^2}
			&\leq C_{x_{\max},x_{\min}}
			\Biggl\{\frac{\max(\sigma_z^4,m^2,n^2)k\log n}{m^2 n L}
			+n^{-ck}+Le^{-cn}+e^{-cmL}
			+{1\over n^{10}}\Biggr\},
		\end{align*}
		for a constant $C>0$ depending on $x_{\min},x_{\max}$.
			We now focus on the remaining scenario of $n<4m$ with $\sigma_z^2 < m$. Here, using the fact that the error is a non-decreasing function of $\sigma_z$ (\prettyref{lmm:monotonicity in sigma}), we obtain the upper bound
			\begin{align*}
				& R_2(\Cc_k,m,n,\sigma_z)
				\le  R_2(\Cc_k,m,n,\sqrt{m})
				=C
				\Biggl\{\frac{k\log n}{n L}
				+n^{-ck}+Le^{-cn}+e^{-cmL}+{1\over n^{10}}\Biggr\},
			\end{align*}
			for a constant $C>0$ depending on $x_{\min},x_{\max}$. The above coincides with our desired upper bound.
			
			
			
			\section{Discussion of the Appendix}
				\prettyref{app:example-proofs} provides proofs of Examples in \prettyref{sec:Cc}.
				\prettyref{app:The lower bound result for 0<b<1} extends \prettyref{thm:lower_bnd} to general $a,b>0$. \prettyref{app:proof_preliminaries}, \prettyref{app:proof_lowerbound_n_geq_4m}, \prettyref{app:lower-bound-oversample} and \prettyref{app:technical} include the
				proofs of auxiliary lemmas from \prettyref{sec:preliminaries} (preliminary results), \prettyref{sec:proof-of-lower-bound-undersample}, \prettyref{sec:proof of lower bound oversample} (lower bound results), and \prettyref{sec:upper bound} (upper bound results), respectively. We include in \prettyref{app:Comparison to the unvarying measurement model} the proof of the fixed forward operator model result \prettyref{thm:main-genr2}. We then elaborate the proof of the sparsity based result \prettyref{thm:ksparse} in
				\prettyref{app:sparse case}. \\
				
			\noindent \textbf{Acknowledgement.} The authors would like to acknowledge corrections pointed out by Songbin Liu, which helped improve the manuscript.

			
			\bibliographystyle{apalike} 
			\bibliography{references}       
			
			
			\appendix

			\section{Proofs of Examples of \prettyref{sec:Cc}}
			
			\label{app:example-proofs}
			
			\subsection*{Proof of Example
				\ref{exam:LipschitzFunc previous}}
			Note that for all $\xv,\yv\in \R^k$ we have $\|g(\xv)-g(\yv)\|_2\le M \|\xv-\yv\|_2$. If $B_1,...,B_r$ are the balls in a ball of radius $\sqrt{k}$ centered at $0$ containing $\e$-covering of $[0,1]^k$, then $g(\Cc)$ is contained in a ball of radius $M\sqrt{k}$. Hence, using \prettyref{prop:estimate of covering numbers}, we have $N_{\e}\pr{g([0,1]^k)} \leq \pr{\frac{2M\sqrt{k}}{\e}+1}^k$.
			
			\subsection*{Proof of \prettyref{exam:sparse}}
			We have
			\begin{align}
				& \mathcal{S}_k = \{\xv \in \mathbb{R}^n \ | \ \|\xv\|_0 \leq k \}= \cup_{1\le i_1 \le \cdots \le i_k \le n} \cb{\xv\in \R^n: x_j=0, j\ne i_1,\cdots,i_k}.
			\end{align}
			Hence, $\mathcal{S}_k$ is a union of ${n \choose k}$ $k$-dimensional subspaces $\cb{\xv\in \R^n: x_j=0, j\ne i_1,\cdots,i_k}$. According to  by \prettyref{prop:estimate of covering numbers} the intersection of each of these subspaces and $B_2(1)$ can be covered by at most $\pr{\frac{2}{\e}+1}^k$ balls of radius $\e$. Hence, $N_{\e}(\Cc) \leq {n \choose k} \pr{\frac{2}{\e}+1}^k$. To obtain the lower bound we notice that according to  \prettyref{prop:estimate of covering numbers}, in order to cover one of the subspaces we need $\Big( \frac{1}{\epsilon} \Big)^k$. 

			\subsection*{Proof of \prettyref{exam:singular value and Minkowski}}
			
			$f(\theta) = D\theta$ is a $\smax(D)$-Lipchitz function of $\theta$. Hence, combining \prettyref{exam:sparse} with a proof similar to the one presented for  \prettyref{exam:LipschitzFunc previous} establishes the result.

			\subsection*{Proof of  \prettyref{exam:piececonstant_1}}

			Note $\Cc \subset D^{-1}(\mathcal{S}_k \cap B_2(0,1))$. By direct calculation,
			\begin{align*}
				D^{-1} := \begin{bmatrix}
					1 & 1 & \cdots & 1 \\
					0 & 1 & \cdots & 1 \\
					0 & 0 & \ddots & 1 \\
					0 & 0 & \cdots & 1 \\
				\end{bmatrix}.
			\end{align*}
			Hence $\smax(D^{-1})\le 
			\nmhs{D^{-1}}=\sqrt{\frac{n(n+1)}{2}}<n$. So this is a special case of \prettyref{exam:singular value and Minkowski}.

			\subsection*{Proof of  \prettyref{exam:piecepolynomial_1}}
			We showed $\sigma_{\min} (D) \geq \frac{1}{n}$ in the proof of \prettyref{exam:piececonstant_1}. Hence, we have
			\[
			\sigma_{\min}(D^{M+1}) \geq (\sigma_{\min}(D))^{M+1} \geq \Big(\frac{1}{n}\Big)^{M+1}. 
			\]

			\section{The extension of \prettyref{thm:lower_bnd} to general $a,b>0$}\label{app:The lower bound result for 0<b<1}

			Suppose our signal class $\Cc \in \mathcal{F}_{a,b,k,n}$ satisfies the polynomial complexity $N_{\e}(\Cc)\le \pr{an^b \over \e}^k$ for general $a,b>0$. Note that we can make the transform
			\begin{align*}
				\pr{an^b \over \e}^k = \pr{(\xmax-\xmin)n \over \e'}^{k'},
			\end{align*}
			where $k':=bk$ and $\e':=\frac{(\xmax-\xmin)\e^{1/b}}{a^{1/b}}$. Put $C'=\frac{\e'}{\e}\Cc \subset [\xmin',\xmax']^n$ where $\xmin'=\frac{\e'}{\e}\xmin$ and $\xmax'=\frac{\e'}{\e}\xmax$, and we have
			\begin{align}
				N_{\e'}(\Cc')=N_{\e}(\Cc)\le \pr{an^b \over \e}^k= \pr{(\xmax-\xmin)n \over \e'}^{k'}.
			\end{align}
			This means $\Cc  \in \mathcal{F}_{a, b,k,n}$ if and only if $\Cc'\in \frac{\e'}{\e}\mathcal{F}_{a_0, b_0,k',n}$ where $a_0=\xmax-\xmin$ and $b_0=1$. Hence by \prettyref{thm:lower_bnd}, we have when $\e'\in (0,1/2)$, $k' \le n^{1-2\e'}$, and $\max(\sigma_z^4,m^2,n^2)k\log n \le m^2n^{1-\e'}L$
			\begin{align*}
					&\sup_{\Cc  \in \mathcal{F}_{a, b,k,n}} R_2(\Cc',m,n,\sigma_z)
					= \sup_{\Cc'  \in \frac{\e'}{\e}\mathcal{F}_{a_0, b_0,k',n}} R_2(\Cc',m,n,\sigma_z)\\
					&= \Omega_{\e',x_{\max}',x_{\min}'}\pr{\frac{\max(\sigma_z^4,m^2,n^2)k\log n}{m^2 n L}}
					= \Omega_{\e,x_{\max},x_{\min},a,b}\pr{\frac{\max(\sigma_z^4,m^2,n^2)k\log n}{m^2 n L}}.
			\end{align*}

			\section{Proof of auxiliary lemmas from \prettyref{sec:preliminaries}}
			\label{app:proof_preliminaries}

			\subsection{Proof of \prettyref{lmm:bound-trace}}

			In the following steps $\otimes$ denotes the Kronecker product. We have
			\begin{align*}
				\begin{split}
					&\trbr{\Sigma^{-1}(\widetilde\Sigma-\Sigma)\Sigma^{-1}(\widetilde\Sigma-\Sigma))}
					=  \operatorname{Vec}(\widetilde\Sigma-\Sigma)^\top   \br{\Sigma^{-1}\otimes \Sigma^{-1}}\operatorname{Vec}(\widetilde\Sigma-\Sigma)\\
					& \ge \nm{\operatorname{Vec}(\widetilde\Sigma-\Sigma)}_2^2 \lmin(\Sigma^{-1}\otimes \Sigma^{-1})\\
					& = \|\widetilde\Sigma-\Sigma\|_{\rm HS}^2 \lmin^2(\Sigma^{-1})=\|\widetilde\Sigma-\Sigma\|_{\rm HS}^2 \lmin^2(\sigma_z^2 I_m+A_lX_o^2 A_l^\top )\\
					& = \|\widetilde\Sigma-\Sigma\|_{\rm HS}^2 \left[\sigma_z^2 +\lmin(A_lX_o^2 A_l^\top ) \right]^2
					\ge \|\widetilde\Sigma-\Sigma\|_{\rm HS}^2 \left(\sigma_z^2+x_{\min}^2\lmin(A_lA_l^\top ) \right)^2.
				\end{split}
			\end{align*}
			
			On the other hand, using 
			$\widetilde\Sigma-\Sigma=\Sigma(\Sigma^{-1}-\widetilde\Sigma^{-1})\widetilde\Sigma=\Sigma A(X^2-\widetilde X^2)A^T\widetilde\Sigma$, we have
			\begin{align*}
				\|\widetilde\Sigma-\Sigma\|_{\rm HS}
				\ge & \lmin(\Sigma)\lmin(\widetilde\Sigma)\nmhs{A(\widetilde{X}_o^2-X_o^2)A^\top }\\
				\ge & \frac{\nmhs{A(\widetilde{X}^2-X^2)A^\top }}{\lmax(\sigma_z^2 I_m+A X^2 A^\top )\lmax(\sigma_z^2 I_m+A\widetilde{X}_o^2A^\top )}
				\ge \frac{\nmhs{A(\widetilde{X}^2-X^2)A^\top }}{\left(\sigma_z^2+x_{\max}^2 \lmax(AA^\top )\right)^2}.
			\end{align*}
			This proves the lower bound. The proof of the upper bound is similar. Note that 
			\begin{align*}
				\begin{gathered}
					\trbr{\Sigma^{-1}(\widetilde\Sigma-\Sigma)\Sigma^{-1}(\widetilde\Sigma-\Sigma))}
					\le \|\widetilde\Sigma-\Sigma\|_{\rm HS}^2 \left(\sigma_z^2+x_{\max}^2\lmax(A_lA_l^\top ) \right)^2,\\
					\|\widetilde\Sigma-\Sigma\|_{\rm HS}\le \frac{\nmhs{A(\widetilde{X}^2-X^2)A^\top }}{\lmin(\sigma_z^2 I_m+A X^2 A^\top )\lmin(\sigma_z^2 I_m+A\widetilde{X}_o^2A^\top )}
					\le \frac{\nmhs{A(\widetilde{X}^2-X^2)A^\top }}{\left(\sigma_z^2+x_{\min}^2 \lmin(AA^\top )\right)^2}.
				\end{gathered}
			\end{align*}

			We use the following linear algebraic result to bound the differences between the inverse of two matrices.

			\subsection{Proof of \prettyref{lmm:decoupling-genr}}\label{app:subsection on the Proof of lemma the decoupling from A}
			The proof closely follows that of  \cite[Lemma 4 and 5]{zhou2024corrections}. However, we obtain sharper results with revised techniques, which we present here. We start with a decoupling result that plays a critical role in our proof:
			\begin{lemma}{\cite[Theorem 3.4.1.]{de2012decoupling}}\label{lmm: abstract decoupling}
				Let $X_1, X_2, \ldots, X_n$ denote random variables with values in measurable space $(S, \mathcal{S})$. Let $(\widetilde{X}_1, \widetilde{X}_2, \ldots, \widetilde{X}_n)$ denote an independent copy of $X_1, X_2, \ldots, X_n$. For $ i \neq j$ let $h_{i,j} : S^2 \rightarrow \mathbb{R}$. Then, there exists a constant $C$ such that for every $t >0$ we have
				\[
				\Prob \pr{\Big|\sum_{i \neq j} h_{i,j} (X_i, X_j) \Big | > t} \leq C \Prob \pr{C \Big|\sum_{i \neq j} h_{i,j} (X_i, \widetilde{X}_j) \Big | > t} . 
				\] 
			\end{lemma}
			For $\{A_l\}_{l=1}^L\in \reals^{m\times n},D\in \reals^{L\times L}$, define 
			\begin{align}
				\bA = \diag\pth{\{A_l\}_{l=1}^L}
				\in \reals^{mL\times nL},
				\quad
				\mathbf D = \diag\pth{\{D\}_{l=1}^L}
				\in \reals^{nL\times nL}.
			\end{align}

			
			For $1\le l \le L$, let $\av_{l,i}^\top $ denote the $i^{\rm th}$ row of matrix $A_l$. We have
			\begin{align}\label{eq:concent:imp:usedecoup}
				\begin{split}
					\|\bA\mathbf D \bA\|^2_{\rm HS}=\sum_{l=1}^L\|A_l DA_l^\top \|^2_{\rm HS} 
					=\sum_{l=1}^L\sum_{i=1}^m \sum_{j=1}^m |\av_{l,i}^\top  D \av_{l,j}|^2
					= \sum_{l=1}^L \sum_{i \ne j} |\av_{l,i}^\top  D \av_{l,j}|^2+\sum_{l=1}^L\sum_{i=1}^m |\av_{l,i}^\top  D \av_{l,i}|^2.
				\end{split} 
			\end{align}
			One of the concentration results that will be used extensively in our paper is the concentration of quadratic functions, a.k.a. Hanson-Wright inequality. 
			\begin{lemma}[Hanson-Wright inequality, \cite{hansonwright71}] \label{lmm:hanson wright} 	
				Let $\bxi= (\xi_1,...,\xi_n)^\top $ be a random vector with independent components with $\E[\xi_i]=0$ and $\|\xi_i\|_{\rm subgau}\leq K$. Let A be an $n\times n$ matrix. Then, for $t>0$, 
				\begin{align}
					\Prob \pr{\br{\Big|\bxi^\top  A\bxi-\E[\bxi ^\top  A\bxi]\Big|>t}}  \leq 2\exp\left(-c\min\left({t^2\over K^4\|A\|_{\rm HS}^2 },{t \over K^2\|A\|_2 }\right)\right),
				\end{align}
				where $c$ is a constant, and $\|\xi_i\|_{\rm subgau} = \inf\{t>0 : \mathbb{E} (\exp (\xi_i^2/t^2)) \leq 2\}. 
				$
			\end{lemma}
			
			We will use the following result to often bound the $\rm HS$ norm.
			
			\begin{lemma}\label{lmm:HS-norm-up}
				Let $A$ denote an arbitrary matrix and $D$ be a diagonal matrix. Then,
				\[
				\|AD\|_{\rm HS} \leq \sigma_{\max}(A) \|D\|_{\rm HS}. 
				\]
			\end{lemma}

			\begin{proof}
				Define $B= A^TA$. Then, we have
				\begin{align}\label{eq:firstbd:la:simple}
					\|AD\|_{\rm HS}^2 = \Tr (DA^TAD)= \Tr(DBD)= \sum_{i} D_{ii}^2 B_{ii}. 
				\end{align}
				Note that if $\mathbf{e}_i$ is the unit vector with a one in the $i^{\text{th}}$ position and zeros elsewhere. Then,
				\begin{equation}\label{eq:firstbd:la:simple2}
					|B_{ii}| = e_i^T B e_i \leq \sigma_{\max}(B) = \sigma_{\max}^2(A)
				\end{equation}
				Combining \eqref{eq:firstbd:la:simple} and \eqref{eq:firstbd:la:simple2}  establishes the desired result. 
			\end{proof}
			
			We get back to the proof of \prettyref{lmm:decoupling-genr}. By using the union bound and \prettyref{lmm:hanson wright}, we have
			\begin{align}\label{eq:diagonal part}
				\begin{split}
					& \P\pr{\sum_{l=1}^L\sum_{i=1}^m |\av_{l,i}^\top  D \av_{l,i}|^2 > Lm\pr{\Tr(D)+t_1}^2}\\
					\le & Lm\P\pr{|\av_{l,i}^\top  D \av_{l,i}|^2 > \pr{\Tr(D)+t_1}^2}\\
					= & Lm\P\pr{|\av_{l,i}^\top  D \av_{l,i}| > \Tr(D)+t_1}
					\le 2mL\exp\left(-c\min\left({t_1^2\over K^4\|\bd\|_{2}^2 },{t_1 \over K^2\|\bd\|_{\infty} }\right)\right),
				\end{split}
			\end{align}
			where $K$ the subgaussian norm of each element of $\av_{l,i}$ and $K\in [1,2]$ is a fixed number. For the off-diagonal part of \eqref{eq:concent:imp:usedecoup}, first note that
			\begin{align}
				\E \br{\sum_{l=1}^L \sum_{i \ne j} |\av_{l,i}^\top  D \av_{l,j}|^2} = Lm(m-1)\sum_{i=1}^n d_i^2=Lm(m-1)\|\bf{d}\|_2^2.
			\end{align}	
			By \prettyref{lmm: abstract decoupling}, there exists a constant $C>0$ such that 
			\begin{align}\label{eq:using the abstract decoupling on off diagonal part}
				\begin{split}
					&{\Prob \Bigg(\bigg| \sum_{l=1}^L \sum_{i \ne j} |\av_{l,i}^\top  D \av_{l,j}|^2 -Lm(m-1) \sum_{i=1}^n d_i^2\bigg|> t \Bigg)} \\ 
					\leq & C \Prob \Bigg(C \bigg| \sum_{l=1}^L \sum_{i \ne j} |\av_{l,i}^\top  D \tilde{\av}_{l,j}|^2 -Lm(m-1) \sum_{i=1}^n d_i^2\bigg|> t \Bigg)  \\
					= & C  \Prob \Bigg(C \bigg| \sum_{l=1}^L \sum_{i=1}^m  \av_{l,i}^\top  D \bigg( \sum_{i \ne j=1}^m \tilde{\av}_{l,j}  \tilde{\av}_{l,j}^\top \bigg) D \av_{l,i} -Lm(m-1) \sum_{i=1}^n d_i^2\bigg|> t \Bigg), 
				\end{split}
			\end{align}
			where $\tilde \av_{l,j}$'s denote the independent copies of $\av_{l,j}$'s for $1\le l\le L$ and $1\le j \le m$. Define $\tilde{A}_l$ as the $m\times n$ Gaussian matrix whose rows are $\tilde{\av}_{l,1}^\top , \ldots, \tilde{\av}_{l,m}^\top $, i.e.,
			\begin{align*}
				\widetilde A_l:=\begin{bmatrix}
					\tilde \av_{l,1}, \hdots,
					\tilde \av_{l,m}
				\end{bmatrix}^\top.
			\end{align*}
			
			Also, let $\tilde{A}_{l,\backslash i}$ denote the matrix that is constructed by removing the $i^{\rm th}$ row of $\tilde{A}_l$. Define
			\[
			\mathbf F := \diag\pth{\{F_l\}_{l=1}^L},
			\quad
			F_l :=\diag\pth{\{D \tilde{A}_{l,\backslash j}^\top  \tilde{A}_{l,\backslash j} D\}_{j=1}^m},
			\quad 
				\lv:=[\vvv_1^\top ,...,\vvv_L^\top ],
				\quad
				\vvv_l:= [\av_{l,1}^\top , \ldots, \av_{l,m}^\top ].
			\]
			Let $\P_{\mathbf {\tilde A}}(\cdot):=\P\br{\cdot \mid \tilde A_1,...,\tilde A_L}$. Note that for any event $\Ec$, by definition $\P_{\mathbf {\tilde A}}(\Ec)=\E[\mathbf{1}_{\Ec}\mid \mathbf {\tilde A}]$. 
			It follows 
			\begin{align}\label{eq:off diagonal part 2}
				\begin{split}
					&  \P \Bigg(C \bigg|\sum_{l=1}^L \sum_{i=1}^m  \av_{l,i}^\top  D \bigg(\sum_{i \ne j=1}^m \tilde{\av}_{l,j}  \tilde{\av}_{l,j}^\top \bigg) D \av_{l,i} -Lm(m-1) \sum_{i=1}^n d_i^2\bigg|> t_2 \Bigg)\\
					=& \mathbb{E}\br{ \P_{\mathbf {\tilde A}}\pr{C\abs{\lv^\top  \mathbf{F} \lv -\E\br{\lv^\top  \mathbf{F} \lv  | \mathbf{\tilde{A}} }}\ge t_2/2}} + \P \pr{C\abs{\E\br{\lv^\top  \mathbf{F} \lv  | \mathbf{\tilde{A}} } - Lm(m-1)\sum_{i=1}^n d_i^2}\ge t_2/2}.
				\end{split}
			\end{align}
			Note that from \prettyref{lmm:hanson wright} (with fixed $\mathbf {\tilde A}$) we have
			\begin{equation}\label{eq:decouple:applied}
				\P_{\mathbf {\tilde A}}\pr{C\abs{\lv^\top  \mathbf{F} \lv -\E\br{\lv^\top  \mathbf{F} \lv  | \mathbf{\tilde{A}} }}\ge t/2} 
				\le  2\exp \left( -c \min \Bigg( \frac{t^2}{4C^2 K^4\|\mathbf{F}\|_{\rm HS}^2 }, \frac{t}{2C K^2\|\mathbf{F}\|_2 }\Bigg) \right).
			\end{equation}

			Define $\widetilde \Ec_{\rm maxsing}:=\bigcap_{l=1}^L \bigcap_{i=1}^m \cb{\smax(\tilde{A}_{l,\backslash i})\le \frac{3}{2}(\sqrt{m}+\sqrt{n})}$, and $E_{\max} = \frac{9}{4} (\sqrt{m}+ \sqrt{n})^2$. 
			By \prettyref{lmm:singvalues}, we have that $\P(\Ec_{\rm maxsing})\ge 1-mL\exp(-cn)$.
			Restricted to the event $\widetilde \Ec_{\rm maxsing}$, we have
			\begin{align}\label{eq:F2:bound}
				\|\mathbf{F}\|_2=\max_{1\le l\le L}\max_{1\le i\le m}\smax\pr{D\tilde{A}_{l,\backslash i}^\top  \tilde{A}_{l,\backslash i}D}
				\le \emax\nm{\bd}_{\infty}^2.
			\end{align}
			and
			\begin{align}\label{eq:m5}
				& \nmhs{\mathbf{F}}^2=\sum_{l=1}^L \nmhs{F_l}^2
				=\sum_{l=1}^L \sum_{i=1}^m \nmhs{D \tilde{A}_{l,\backslash i}^\top  \tilde{A}_{l,\backslash i} D}^2 
				\nonumber\\
				&
				=  \sum_{l=1}^L \sum_{i=1}^m \trbr{D \tilde{A}_{l,\backslash i}^\top  \tilde{A}_{l,\backslash i} DD \tilde{A}_{l,\backslash i}^\top  \tilde{A}_{l,\backslash i} D} 
				\nonumber\\
				&=  \sum_{l=1}^L \sum_{i=1}^m \trbr{ \tilde{A}_{l,\backslash i}^\top  \tilde{A}_{l,\backslash i} DD \tilde{A}_{l,\backslash i}^\top  \tilde{A}_{l,\backslash i} DD} 
				\nonumber\\
				&\overset{(a)}{\le} \sum_{l=1}^L \sum_{i=1}^m (\sum_p \sigma_{p}( \tilde{A}_{l,\backslash i}^\top  \tilde{A}_{l,\backslash i} DD))^2 =\sum_{l=1}^L \sum_{i=1}^m  \nmhs{\tilde{A}_{l,\backslash i}^\top  \tilde{A}_{l,\backslash i} DD}^2
				\le Lm \br{\emax}^2
				\|\bd^2\|_{2}^2, 
			\end{align}
			where the last inequality followed using  \prettyref{lmm:HS-norm-up} and (a) followed using the result below.
			\begin{lemma}\cite[Theorem 7.4.1.1]{horn2012matrix} \label{lmm:vonnueman} 
				Let $A \in \mathbb{R}^{n \times n} $ and $B \in \mathbb{R}^{n \times n}$  denote two matrices with singular values $\sigma^{A}_i$ and $\sigma^B_i$. Then, we have
				$
				{\rm Tr}(AB)  \leq \sum \sigma^{A}_i \sigma_i^B.
				$
			\end{lemma}
			Combining, \eqref{eq:decouple:applied}, \eqref{eq:F2:bound}, and \eqref{eq:m5} we have
			\begin{align}\label{eq:off_diag:2}
				\begin{split}
					\mathbb{E}\br{ \P_{\mathbf {\tilde A}}\pr{C\abs{\lv^\top  \mathbf{F} \lv -\E\br{\lv^\top  \mathbf{F} \lv  | \mathbf{\tilde{A}} }}\ge t_2/2}}
					\leq 
					2 e^{ -c \min \pr{\frac{4t_2^2}{81C^2 K^4  \|\bd^2\|_{2}^2  m L ( \sqrt{n} + \sqrt{m})^4 }, \frac{2t_2}{9C K^2 \| \bd\|_{\infty}^{2} ( \sqrt{n} + \sqrt{m})^2} }}.
				\end{split}
			\end{align}
			Now note that 
			\begin{align}\label{eq:off:diag:on:diag}
				\begin{split}
					& \P \pr{C\abs{\E\br{\lv^\top  \mathbf{F} \lv  | \mathbf{\tilde{A}} } - Lm(m-1)\sum_{i=1}^n d_i^2}\ge t/2} \\
					=& \P \pr{C\abs{\br{ \sum_{l=1}^L \sum_{i=1}^m \sum_{j \neq i, j=1}^m \tilde{\ba}_{l,j}^\top D^2 \tilde{\ba}_{l,j}} - Lm(m-1)\sum_{i=1}^n d_i^2}\ge t/2}  \\
					\overset{(a)}{\leq}& m \P \pr{\abs{\sum_{l=1}^L \sum_{j\ne i_0, j=1}^m   \tilde{\ba}_{l,j}^\top D^2 \tilde{\ba}_{l,j} - L(m-1)\sum_{i=1}^n d_i^2}\ge t/2m}  \\
					\overset{(b)}{\leq}& 2 m  \exp\left(-c\min\left({t^2\over 4K^4m^2\nmhs{\mathbf{D}_{(m-1)L}^2}^2 },{t  \over 2K^2m\nm{\mathbf{D}_{(m-1)L}^2}_2 }\right)\right)  \\
					{\leq}& 2 m  \exp\left(-c\min\left({t^2\over 4K^4L m^3 \|\bd^2\|_2^2 },{t  \over 2K^2m \|\bd^2\|_\infty }\right)\right),
				\end{split}
			\end{align}
			where to obtain inequality (a) we have used the union bound (the distribution for $\sum_{j\ne i_0, j=1}^m   \tilde{\ba}_{l,j}^\top D^2 \tilde{\ba}_{l,j}$ is the same for all $1\le i_0 \le m$), and to obtain (b) we have used \prettyref{lmm:hanson wright} for the $L(m-1)n \times L(m-1)n$ matrix
			\begin{align*}
				\mathbf D_{(m-1)L}^2 := \begin{bmatrix}
					D^2 & 0 & \cdots & 0 \\
					0 & D^2 & \cdots & 0 \\
					0 & 0 & \cdots & 0 \\
					0 & 0 & \cdots & D^2 \\
				\end{bmatrix}.
			\end{align*}

			Finally, combining \eqref{eq:diagonal part},\eqref{eq:using the abstract decoupling on off diagonal part}, \eqref{eq:off diagonal part 2}, \eqref{eq:off_diag:2}, and \eqref{eq:off:diag:on:diag} we have that
			\begin{align} 
				\begin{split}
					&\Prob \pr{ \br{\sum_{l=1}^L\|A_l DA_l^\top \|^2_{\rm HS} > Lm\pr{\Tr(D)+t_1}^2+Lm(m-1)\|\bd\|_{2}^2 + t_2} \cap \widetilde \Ec_{\rm maxsing}}\\
					\le & \P\pr{\sum_{l=1}^L\sum_{i=1}^m |\av_{l,i}^\top  D \av_{l,i}|^2 > Lm\pr{\Tr(D)+t_1}^2}\\
					&+\Prob \pr{\br{~\abs{ \sum_{l=1}^L \sum_{i \ne j} |\av_{l,i}^\top  D \av_{l,j}|^2 -Lm(m-1) \sum_{i=1}^n d_i^2}> t_2 } \cap \widetilde \Ec_{\rm maxsing} }\\
					\leq & 2mL\exp\left(-c\min\left({t_1^2\over K^4\|\bd\|_{2}^2 },{t_1 \over K^2\|\bd\|_{\infty} }\right)\right)\\
					&+2C  \exp \left( -c \min \pr{\frac{4t_2^2}{81C^2 K^4   \|\bd^2\|_{2}^2  m L ( \sqrt{n} + \sqrt{m})^4 }, \frac{2t_2}{9CK^2 \| \bd\|_{\infty}^{2} ( \sqrt{n} + \sqrt{m})^2} }\right)  \nonumber \\
					&+ 2 m  \exp\left(-c\min\left({t_2^2\over 4K^4L m^3 \|\bd^2\|_2^2 },{t_2  \over 2K^2Lm^2 \|\bd^2\|_\infty }\right)\right),
				\end{split}
			\end{align}
			
			for some constants $c$ and $C$. On the other hand,
			\begin{align*}
				&\Prob \pr{ \br{\sum_{l=1}^L\|A_l DA_l^\top \|^2_{\rm HS} < Lm(m-1) \|\bd\|_2^2 - t} \cap \widetilde \Ec_{\rm maxsing}}\\
				\le & \Prob \Bigg(\bigg[ \sum_{l=1}^L \sum_{i \ne j} |\av_{l,i}^\top  D \av_{l,j}|^2 -Lm(m-1) \sum_{i=1}^n d_i^2< -t \bigg] \cap \widetilde \Ec_{\rm maxsing} \Bigg)\\
				\le &\Prob \Bigg(\bigg[\bigg| \sum_{l=1}^L \sum_{i \ne j} |\av_{l,i}^\top  D \av_{l,j}|^2 -Lm(m-1) \sum_{i=1}^n d_i^2\bigg|> t \bigg] \cap \widetilde \Ec_{\rm maxsing} \Bigg)\\
				\le & 2C  \exp \left( -c \min \pr{\frac{4t^2}{81C^2K^4   \|\bd^2\|_{2}^2  m L ( \sqrt{n} + \sqrt{m})^4 }, \frac{2t}{9CK^2 \| \bd\|_{\infty}^{2} ( \sqrt{n} + \sqrt{m})^2} }\right)  \nonumber \\
				&+ 2 m  \exp\left(-c\min\left({t^2\over 4K^4L m^3 \|\bd^2\|_2^2 },{t  \over 2K^2Lm^2 \|\bd^2\|_\infty }\right)\right),
			\end{align*}
			for some constants $c$ and $C$.

			\subsection{Proof of  \prettyref{lmm:monotonicity in m} and \prettyref{lmm:monotonicity in sigma}}\label{app:subsection on the Proof of monotonicity lemmas}
			\begin{proof}[Proof of \prettyref{lmm:monotonicity in m}]
				Let $m<m'$ be two positive integers and consider two scenarios of our model:
				\begin{align}
					\begin{split}
						\yv_l =& A_l X_o \wv_l + \zv_l, \text{ 
							for } l=1,\ldots,L;
							\quad
						\yv_l' = A_l' X_o \wv_l' + \zv_l', \text{ 
							for } l=1,\ldots,L.
					\end{split}
				\end{align}
				where we have independent $\wv_l \sim \mathcal N(0,I_n)$, $\wv_l' \sim \mathcal N(0,I_n)$, $\zv_l \sim \mathcal N(0,\sigma_z I_m)$, and $\zv_l' \sim \mathcal N(0,\sigma_z I_{m'})$. We would like to show
				\begin{align}\label{eq:R_2 m m'}
					R_2(\Cc_k,m',n,\sigma_z)\le R_2(\Cc_k,m,n,\sigma_z).
				\end{align}
				Indeed, for each $1\le l \le L$, we look at each component of vectors $\yv_l$ and $\yv_l'$: For $1\le i,i' \le m$,
				\begin{align}
					y_{l,i}=& \sum_j A_{l,ij}x_{o,j}w_{l,j}+z_{l,i},
					\quad
					y_{l,i'}'= \sum_j A_{l,i'j}'x_{o,j}w_{l,j}'+z_{l,i'}'.
				\end{align}
				For $1\le l \le L$, let $\yv_l'|_m$ denote the truncation of the $m'$-dimensional vector on its first $m$ components.
				
				Let $\mathscr E_m$ denote the collection of all estimators for the first scenario, let $\mathscr E_{m'}'$ denote the collection of all estimators for the second scenario, and let $\mathscr E_{m}'\subset \mathscr E_{m'}'$ denote the subcollection of estimators for the second scenario that only use the information of $\yv_l'|_m$ for $1\le l \le L$.
				
				In these two scenarios, we construct estimators $\xvh(\yv_1,...,\yv_L)$ and $\xvh'(\yv_1',...,\yv_L')$, and we can always view $\xvh(\yv_1,...,\yv_L)$ as a special case of $\xvh'(\yv_1',...,\yv_L')$ where we only use the information of truncations $\yv_1'|_m,...,\yv_L'|_m$. Therefore, \eqref{eq:R_2 m m'} follows from the following relation
				\begin{align*}
					\begin{split}
						&\inf_{\xvh\in \mathscr E_m} \sup_{\xv_o\in \Cc_k}\E \br{  \| \xvh(\yv_1,  \ldots, \yv_L) - \xv_o\|_2^2}\\
						= & \inf_{\xvh'\in \mathscr E_m'} \sup_{\xv_o\in \Cc_k}\E \br{  \| \xvh'(\yv_1'|_m,  \ldots, \yv_L'|_m) - \xv_o\|_2^2}
						\ge  \inf_{\xvh'\in \mathscr E_{m'}'} \sup_{\xv_o\in \Cc_k}\E \br{  \| \xvh'(\yv_1',  \ldots, \yv_L') - \xv_o\|_2^2}.
					\end{split}
				\end{align*}
			\end{proof}

			\begin{proof}[Proof of \prettyref{lmm:monotonicity in sigma}]
				The argument here follows closely the proof of \cite[Lemma 3.1]{malekian2025speckle}. Consider two scenarios of our spec6
				e noise model:
				\begin{align*}
					\begin{split}
						\yv_l =& A_l X_o \wv_l + \zv_l, \text{ 
							for } l=1,\ldots,L,
							\quad
						\yv_l' = A_l X_o \wv_l' + \zv_l', \text{ 
							for } l=1,\ldots,L.
					\end{split}
				\end{align*}
				where we have independent $\wv_l \sim \mathcal N(0,I_n)$, $\wv_l' \sim \mathcal N(0,I_n)$, $\zv_l \sim \mathcal N(0,\sigma_z I_m)$, and $\zv_l' \sim \mathcal N(0,\sigma_z' I_m)$ with $\sigma_z'>\sigma_z>0$. We would like to show
				\begin{align*}
					R_2(\Cc_k,m,n,\sigma_z)\le R_2(\Cc_k,m,n,\sigma_z').
				\end{align*}
				Let $\xvh(\yv_1',...,\yv_L')$ be any estimator of $\xvh_o$ with observations $\yv_1',...,\yv_L'$. Conditioning on $A_1,...,A_L$, we have
				\begin{align*}
					(\yv_1',...,\yv_L')\overset{d}{=}(\yv_1+\uv_1,...,\yv_L+\uv_L),
				\end{align*}
				where $\uv_1,...,\uv_L$ are i.i.d. $\mathcal N(0,\sqrt{\sigma_z'^2-\sigma_z^2}\cdot I_m)$. It follows that
				\begin{align*}
					\E\br{\frac{\nm{\xvh(\yv_1',...,\yv_L')-\xv_o}^2}{n}}=\E\br{\frac{\nm{\xvh(\yv_1+\uv_1,...,\yv_L+\uv_L)-\xv_o}^2}{n}}.
				\end{align*}
				Furthermore, if we let $\E_Y$ denote the conditional expectation $\E_Y[\cdot]:=\E[\cdot \mid A_1,...,A_L, \yv_1,...,\yv_L]$, then by the tower rule and Jensen's inequality we have
				\begin{align*}
					\begin{split}
						& \E \br{\|\xvh (\yv_1+\uv_1,  \ldots, \yv_L+\uv_L)- \xv_o\|_2^2}  \\
						= & \E \br{\E_Y\|\xvh (\yv_1+\uv_1,  \ldots, \yv_L+\uv_L)- \xv_o\|_2^2}  
						\geq  \E \br{  \| \E_Y\br{\xvh (\yv_1+\uv_1,  \ldots, \yv_L+\uv_L)} - \xv_o\|_2^2}.
					\end{split}
				\end{align*}
				It follows that
				\begin{align*}
					\begin{split}
						&\sup_{\xv_o\in \Cc_k}\E \br{\|\xvh (\yv_1+\uv_1,  \ldots, \yv_L+\uv_L)- \xv_o\|_2^2}
						\geq \sup_{\xv_o\in \Cc_k}\E \br{  \| \E_Y\br{\xvh (\yv_1+\uv_1,  \ldots, \yv_L+\uv_L)} - \xv_o\|_2^2}.
					\end{split}
				\end{align*}

				Now, treating $\widehat{\zv}_o (\yv_1,  \ldots, \yv_L):=\E_Y\br{\xvh (\yv_1+\uv_1,  \ldots, \yv_L+\uv_L)}$ as an estimator of $\xv_o$ using only the information of $\yv_1,...,\yv_L$, we have for every estimator $\widehat{\zv}(\yv_1',...,\yv_L')$ of the second scenario
				\begin{align*}
					\begin{split}
						\sup_{\xv_o\in \Cc_k}\E \br{\|\widehat{\zv} (\yv_1',  \ldots, \yv_L')- \xv_o\|_2^2}
						\geq \sup_{\xv_o\in \Cc_k}\E \br{  \| \widehat{\zv}_o (\yv_1,  \ldots, \yv_L) - \xv_o\|_2^2}
						\geq \inf_{\widehat{\zv}} \sup_{\xv_o\in \Cc_k}\E \br{  \| \widehat{\zv} (\yv_1,  \ldots, \yv_L) - \xv_o\|_2^2}.
					\end{split}
				\end{align*}
				Therefore
				\begin{align*}
					\begin{split}
						& R_2(\Cc_k,m,n,\sigma_z')\\
						= &\frac{1}{n}\inf_{\widehat{\zv}}\sup_{\xv_o\in \Cc_k}\E \br{\|\widehat{\zv} (\yv_1',  \ldots, \yv_L')- \xv_o\|_2^2}
						\geq  \frac{1}{n}\inf_{\widehat{\zv}} \sup_{\xv_o\in \Cc_k}\E \br{  \| \widehat{\zv} (\yv_1,  \ldots, \yv_L) - \xv_o\|_2^2}
						\geq R_2(\Cc_k,m,n,\sigma_z).
					\end{split}
				\end{align*}
				
			\end{proof}

			\section{Proofs of results of \prettyref{sec:proof-of-lower-bound-undersample}}
			
			\label{app:proof_lowerbound_n_geq_4m}
			
			\subsection{Proof of \prettyref{lmm:r-bound}}
			We first establish the lower bound. Let $\xv \in \calS_{\rm sep}$. Consider all vectors obtained from $\xv$ by selecting $k/4$ intervals and flip the values of the entries whose indices belong to those intervals; if the value is $\bar{x}+ \delta_r$ switch that to $\bar{x}$, and if the value is $\bar{x}$ switch that to $\bar{x}+ \delta_r$.  Denote the collection of all such vectors by $B(\xv)$. Note that if $\cup_{\xv \in \calS_{\rm sep}} B(\xv)$ does not cover $\mathcal{X}^{\rm finite}$, it means that there exists another vector $\xv \in \mathcal{X}^{\rm finite}$ that is different from all elements of $\calS_{\rm sep}$ in at least $k/4$ intervals, which is in contradiction with the fact that $\calS_{\rm sep}$ includes \textbf{all} such vectors. Hence, we can conclude that
			\begin{equation}\label{eq:XsubsetB}
				\mathcal{X}^{\rm finite} \subset \cup_{\xv \in \calS_{\rm sep}} B(\xv). 
			\end{equation}
			Note that 
			\begin{equation}\label{eq:size:Bx}
				B(\xv) = {N_{\rm div} \choose k/4}.
			\end{equation}
			Combining \eqref{lower bound on X finite}, \eqref{eq:XsubsetB}, and \eqref{eq:size:Bx} and assuming that $r$ is the size of $\calS_{\rm sep}$ we have
			\begin{align}\label{eq:m9}
				r\ge {{N_{\rm div} \choose k}}/{{N_{\rm div} \choose k/4}}.
			\end{align}
			Using the following classical bounds for ${n \choose k}$
			\begin{align}\label{eq:classical estimate for binormial coefficients}
				\left(\frac{n}{k}\right)^k \leq {n \choose k} \leq \left( \frac{en}{k}\right)^k,
			\end{align}
			we have 
			\begin{align}\label{eq:10}
				r\ge {{N_{\rm div} \choose k/2}}/{{N_{\rm div} \choose k/4}} \geq { \pth{\frac{N_{\rm div}}{k/2}}^{\frac{k}{2}}}/{\pth{\frac{e N_{\rm div}}{k/4}}^{\frac{k}{4}} }  = (\frac{N_{\rm div}}{ e k})^{\frac{k}{4}}.  
			\end{align}
			We get the desired upper bound by combining \eqref{lower bound on X finite}, \eqref{eq:classical estimate for binormial coefficients} with $r\leq |\mathcal{X}^{\rm finite}|$, as $\calS_{\rm sep} \subset \mathcal{X}^{\rm finite}$.

			\subsection{Proof of \prettyref{lmm: upper bound for KL divergence}}
			As we discussed in \prettyref{sec:discretization} we have
			\begin{equation*}
				\P_{\xv} \sim \otimes_{l=1}^L N(\boldsymbol{0},\Sigma_l^{{-1}}(\xv)) = N(\boldsymbol{0},\Sig^{{-1}}(\xv)),
				\quad
				\Sigma_l=\Sigma_l(\xv):=(\sigma_z^2 I_m + A_l X^2 A_l^\top )^{-1},
				\quad 
				\Sig(\xv):=\diag\pr{\{\Sigma_l(\xv)\}_{i=1}^L}.
			\end{equation*}

			We use the following result to quantify the KL divergence of two Gaussian distributions:
			\begin{proposition} \cite{duchi2007derivations}\label{prop:Classical formula for KL divergence of normal distributions}
				If $\mathbb Q_j\sim N(\boldsymbol{\mu}_j,\boldsymbol{\Lambda}_j), j=1,2$, are two $d$-dimensional multivariate normal distributions, then
				\begin{align*}
					\KL(\Q_1\parallel \Q_2)=\frac{1}{2}\left[\log\frac{\det \boldsymbol{\Lambda}_2}{\det \boldsymbol{\Lambda}_1} - d + \Tr \pr{ \boldsymbol{\Lambda}_2^{-1}\boldsymbol{\Lambda}_1 } + (\boldsymbol{\mu}_2 - \boldsymbol{\mu}_1)^\top  \boldsymbol{\Lambda}_2^{-1}(\boldsymbol{\mu}_2 - \boldsymbol{\mu}_1)\right].
				\end{align*}
			\end{proposition} 
			
			Using \prettyref{prop:Classical formula for KL divergence of normal distributions}, we condition on $A_1,...,A_L$, with $\Lambda_1=\Sig(\xv_i)^{-1},\Lambda_2=\Sig(\xv_j)^{-1}$, to get
			\begin{align}\label{eq:upperKL1}
				\begin{split}
					\KL(\P_{\xv_i}\parallel \P_{\xv_j})=&\frac{1}{2}\left[\log\frac{\det \Sig(\xv_j)^{-1}}{\det \Sig(\xv_i)^{-1}} - mL + \Tr \pr{ \Sig(\xv_j)\Sig(\xv_i)^{-1}}\right] \\
					=& \frac{1}{2}\left[\log \det \Sig(\xv_j)^{-1}\Sig(\xv_i)  + \Tr \pr{ \br{\Sig(\xv_j)-\Sig(\xv_i)}\Sig(\xv_i)^{-1}}\right]\\
					=& \frac{1}{2}\left[-\log \det \Sig(\xv_j)\Sig(\xv_i)^{-1}  + \Tr \pr{ \br{\Sig(\xv_j)-\Sig(\xv_i)}\Sig(\xv_i)^{-1}}\right].
				\end{split}
			\end{align}
			In order to find an upper bound for the KL divergence, we use the mean value theorem. 
			By applying the mean value theorem to $\log \det \Sig(\xv_j)\Sig(\xv_i)^{-1}$ and defining $\lambda_q$ as the $q^{\rm th}$ eigenvalue of the matrix $\Sig(\xv_i)^{-\frac{1}{2}} \pr {\Sig(\xv_j)-\Sig(\xv_i)}\Sig(\xv_i)^{-\frac{1}{2}}$, we obtain
			
			\begin{align}\label{eq:mvt:logdet1}
				\begin{split}
					&-\log \det \Sig(\xv_j)\Sig(\xv_i)^{-1}= -\log\det \br{\Sig(\xv_i)^{-1/2}\Sig(\xv_j)\Sig(\xv_i)^{-1/2}}\\
					=&- \log \det\br{\Sig(\xv_i)^{-\frac{1}{2}}\pr{\Sig(\xv_j)-\Sig(\xv_i)}\Sig(\xv_i)^{-\frac{1}{2}}+I_{mL}}\\
					=& -\sum_{i=1}^{mL}\log(1+\la_q)
					\overset{(a)}{=}- \sum_{i=1}^{mL}\pr{\la_q-\frac{\la_q^2}{2(1+\la_q')^2}}\\
					=& -\Tr\pr{\Sig(\xv_i)^{-\frac{1}{2}}\pr{\Sig(\xv_j)-\Sig(\xv_i)}\Sig(\xv_i)^{-\frac{1}{2}}}+\sum_{i=1}^{mL}\frac{\la_q^2}{2(1+\la_q')^2},
				\end{split}
			\end{align}
			where to obtain (a) we have used the Taylor expansion for $\log(1+ \lambda_q)$, and defined $\lambda_q'$ as a point between zero and $\lambda_q$. Note that this eigenvalue can be negative. Since we have $\Tr \pr{ \br{\Sig(\xv_j)-\Sig(\xv_i)}\Sig(\xv_i)^{-1}}=\Tr \pr{\Sig(\xv_i)^{-\frac{1}{2}}\pr{\Sig(\xv_j)-\Sig(\xv_i)}\Sig(\xv_i)^{-\frac{1}{2}}}$, by combining \eqref{eq:upperKL1} and \eqref{eq:mvt:logdet1} we obtain
			\begin{equation}\label{eq:upper:KL2}
				\KL(\P_{\xv_i} \parallel \P_{\xv_j} ) = \sum_{i=1}^{mL} \frac{\lambda_q^2}{4(1+ \lambda_q')^2}. 
			\end{equation}
			Since $\lambda_q'$ can be a negative number, in order to obtain a useful upper bound for $\KL(\P_{\xv_i} \parallel \P_{\xv_j} )$ as is required by Fano's inequality we need to find an upper bound for $|\lambda'_q|$. But since this quantity is between zero and $\lambda_q$ we can bound $|\lambda_q|$ instead. We have
			
			\begin{align}\label{eq:lambdaqupper}
				\begin{split}
					&\left|\la_q\pr{\Sig(\xv_i)^{-\frac{1}{2}}\pr{\Sig(\xv_j)-\Sig(\xv_i)}\Sig(\xv_i)^{-\frac{1}{2}}}\right| \\
					{\leq} &  \frac{\sigma_{\max} \pr{\Sig(\xv_j)-\Sig(\xv_i)}}{\sigma_{\min}(\Sig(\xv_i))}  \\
					\overset{(a)}{\le} & \frac{\sigma_{\max} \pr{\Sig(\xv_j)^{-1} -\Sig(\xv_i)^{-1}}}{\sigma_{\min}(\Sig(\xv_i)) \sigma_{\min}\pr{\Sig(\xv_i)^{-1}} \sigma_{\min}\pr{\Sig(\xv_j)^{-1}} }  \\
					= & \frac{\max_{1\le l\le L}\br{\sigma_{\max}(\sigma_z^2 I_m + A_lX_i^2 A_l^\top )}  \max_{1\le l\le L}\br{\sigma_{\max} (A_lX_i^2 A_l^\top - A_lX_j^2 A_l^\top ) }}{\sigma_{\min}\pr{\Sig(\xv_i)^{-1}} \sigma_{\min}\pr{\Sig(\xv_j)^{-1}}}.  \\
					\leq & \frac{\br{\sigma_z^2+x_{\max}^2 \max_{1\le l\le L}\lambda_{\max}(A_lA_l^\top )} \max_{1\le l\le L}\lmax(A_lA_l^\top )  \|{\xv_i}^2 - \xv_j^2\|_\infty}{\left(\sigma_z^2+x^2_{\min} \min_{1\le l \le L}\lambda_{\min} (A_lA_l^\top )\right)^2},
				\end{split}
			\end{align}
			where $(a)$ follows from the following result.
			\begin{lemma}{\cite[Lemma 6.1]{chen2024bagged}}\label{lmm:boundeigenvalues}
				Let $B,C\in \reals^{n\times n}$ be symmetric, invertible matrices. Then $\nm{B^{-1}- C^{-1}}_2\leq \sigma_{\max}\pr{B^{-1}- C^{-1}}\le \frac{\sigma_{\max} (B-C)}{\sigma_{\min}(B) \sigma_{\min}(C) }.$
			\end{lemma} 
			
			As we discussed before we would like to use our upper bounds for $\KL(\P_{\xv_i} \parallel \P_{\xv_j})$ for the Fano's inequality. Hence, in the rest of the proof, we assume that $\xv_i, \xv_j \in \calS_{\rm sep}$. This implies that $\|\xv_i-\xv_j \|_{\infty} \leq \delta_r$. Hence, we can simplify \eqref{eq:lambdaqupper} in the following way
			\begin{align}
				\left|\la_q\pr{\Sig(\xv_i)^{-\frac{1}{2}}\pr{\Sig(\xv_j)-\Sig(\xv_i)}\Sig(\xv_i)^{-\frac{1}{2}}}\right|
				\le &\frac{\br{\sigma_z^2+x_{\max}^2\cdot\emax}\cdot \emax}{\pr{\sigma_z^2+x_{\min}^2\cdot E_{\min}}^2}\cdot 2x_{\max}\|\xv_i-\xv_j\|_{\infty}\\
				\le &\frac{\br{\sigma_z^2+x_{\max}^2\cdot\emax}\cdot \emax}{\pr{\sigma_z^2+x_{\min}^2\cdot E_{\min}}^2}\cdot 2x_{\max}\de_r. \label{eq:upp:lambdaq}
			\end{align}
			Hence, choosing $\delta_r$ such that $\frac{\br{\sigma_z^2+x_{\max}^2\cdot\emax}\cdot \emax}{\pr{\sigma_z^2+x_{\min}^2\cdot E_{\min}}^2}\cdot 2x_{\max}\de_r < \frac{1}{2}$, we use \eqref{eq:upper:KL2} to get
			\begin{align}\label{bound for the KL divergence before restricting A}
				\begin{split}
					\KL(\P_i\parallel \P_j)
					\le & \sum_{i=1}^{mL}\frac{\la_q^2}{4(1+\la_q')^2}\le \sum_{i=1}^{mL}\la_q^2 \\
					=& \Tr \pr{ \Sig(\xv_i)^{-1}\br{\Sig(\xv_j)-\Sig(\xv_i)} \Sig(\xv_i)^{-1}\br{\Sig(\xv_j)-\Sig(\xv_i)}} \\
					{\le} & {\left(\sigma_z^2+x_{\max}^2 \max_{ 1\le l \le L} \lmax(A_lA_l^\top ) \right)^2 \over \left(\sigma_z^2+x_{\min}^2 \min_{1\le l \le L} \lmin(A_lA_l^\top )\right)^4} 
					\sum_{l=1}^L
					\nmhs{A_l(X_i^2-X_j^2)A_l^\top }^2,
				\end{split}
			\end{align}
			where to obtain the last inequality we have used \prettyref{lmm:bound-trace}.
			
			\subsection{Proof of \prettyref{lmm:bound1}}
			
			We set $\bd_{i,j}:= \bx_i^2-\bx_j^2$, and define $D_{i,j}=\diag(\bd_{i,j})$) to get:
			\begin{align}
				\begin{split}
					& \|\bd_{i,j}\|_{\infty} = \max_p |x_{i,p}^2 - x_{j,p}^2| \leq 2x_{\max} \|\bx_i -\bx_j\|_\infty, 
					\quad
					\|\bd_{i,j}^2\|_{\infty} \leq 4x_{\max}^2 \|\bx_i -\bx_j\|_\infty^2,\\
					& \|\bd_{i,j}\|_2^2=\sum_{p=1}^n \pr{x_{i,p}^2-x_{j,p}^2}^2\le 4x_{\max}^2\sum_{p=1}^n \pr{x_{i,p}-x_{j,p}}^2=4x_{\max}^2\|\xv_i-\xv_j\|_2^2,\\
					& \|\bd_{i,j}^2\|_2^2=\sum_{p=1}^n \pr{x_{i,p}^2-x_{j,p}^2}^4\le 16x_{\max}^4\sum_{p=1}^n \pr{x_{i,p}-x_{j,p}}^4=16x_{\max}^4\|\xv_i-\xv_j\|_4^4.
				\end{split}
			\end{align}
			We choose the following values of $t_{1,i,j}$ and $t_{2,i,j}$ to apply the upper tail bound in \prettyref{lmm:decoupling-genr}
			\begin{align*}
				t_{1,i,j}:=& C_{t_1}\pr{x_{\max}\nm{\xv_i-\xv_j}_{2}\sqrt{\log(mLr^2)}+ x_{\max}\nm{\xv_i-\xv_j}_{\infty}\log(mLr^2)};\\
				t_{2,i,j}:=& C_{t_2}\log m \pr{x_{\max}^2 \nm{\xv_i-\xv_j}_{4}^2 \sqrt{mL(\sqrt m+\sqrt n)^4\log r^2}+x_{\max}^2\nm{\xv_i-\xv_j}_{\infty}^2 (\sqrt n +\sqrt m)^2 \log r^2 },
			\end{align*}
			where $C_{t_1}$ and $C_{t_2}$ are two constants. To apply \prettyref{lmm:decoupling-genr}, we note that the non-constant terms appearing in the exponent of \prettyref{lmm:decoupling-genr} obey the following lower bounds
			\begin{enumerate}[label=(\alph*)]
				\item $\frac{t_{1,i,j}^2}{\|\bd_{i,j}\|_2^2}$ is bounded from below by
				$\frac{t_{1,i,j}^2}{4 x_{\max}^2\|\bx_i -\bx_j\|_2^2} = \Omega( \log (mL r^2))$,
				\item $\frac{t_{1,i,j}}{\|\bd_{i,j}\|_{\infty}}$ is bounded from below by 
				$\frac{t_{1,i,j}}{2 x_{\max} \|\bx_i-\bx_j\|_{\infty}} = \Omega( \log (mL r^2))$,
				\item $\frac{t_{2,i,j}^2}{mL(\sqrt{n} + \sqrt{m})^4 \|\bd_{i,j}^2\|_2^2}$ is bounded from below by 
				$\frac{t_{2,i,j}^2}{16 x_{\max}^4 mL(\sqrt{n} + \sqrt{m})^4 \|\bx_i-\bx_j\|_4^4} = \Omega (\log (mLr^2))$,
				\item $\frac{t_{2,i,j}}{(\sqrt{n} + \sqrt{m})^2 \|\bd_{i,j}\|_{\infty}^2}$ is bounded from below by 
				$\Omega (\log (mLr^2))$,
				\item $\frac{t_{2,i,j}^2}{Lm^3 \|\bd_{i,j}\|_2^2 }$ is bounded from below by $\Omega ((\log m)^2 \log r)$,
				\item $\frac{t_{2,i,j}}{m \|\bd_{i,j}\|_{\infty}}$ is bounded from below by $\Omega ((\log m)(\log r)))$.
			\end{enumerate}
			In view of the above definition, consider the event 
			\begin{equation}\label{eq:dcpl}
				\Ec_{\rm dcpl}:=\bigcap_{1\le i< j\le r}\br{\sum_{l=1}^L\|A_l (X_i^2-X_j^2)A_l^\top \|^2_{\rm HS} 
					< Lmt_{1,i,j}^2+Lm(m-1)\|\bd_{i,j}\|_{2}^2 + t_{2,i,j}}.
			\end{equation}
			Then, using \prettyref{lmm:decoupling-genr}, a union bound for all $1\le i<j\le r$, and $\Tr(D_{ij})=0, D_{ij}=\diag(\bd_{ij}), \bd_{ij}=\bx_i^2-\bx_j^2$ (see \eqref{eq:traceCondition}), the above display implies for large constants $C_{t_1},C_{t_2}$
			\begin{align}
				&\P(\Ec_{\rm dcpl}^c \cap \widetilde\Ec_{\rm maxsing})
				\nonumber\\
				&\le \sum_{1\leq i<j\leq r}
				\Prob \pr{ \br{\sum_{l=1}^L\|A_l D_{ij}A_l^\top \|^2_{\rm HS} > Lm\pr{\Tr(D_{ij})+t_{1,i,j}}^2+Lm(m-1)\|\bd_{i,j}\|_{2}^2 + t_{2,i,j}}\cap \widetilde \Ec_{\rm maxsing}}
				\nonumber\\
				&
				\leq r^2\exp\sth{-\tilde C\pth{\log(mLr^2)+(\log m)(\log r)}},
			\end{align}
			for a large constant $\tilde C$. Hence, by making $C_{t_1},C_{t_2}$ large enough such that $\tilde C>10$, we have
			\begin{align*}
				&\P(\Ec_{\rm dcpl}\cap \widetilde\Ec_{\rm maxsing})
				\ge   1-\P(\Ec_{\rm dcpl}^c\cap \widetilde\Ec_{\rm maxsing})-\P(\widetilde\Ec_{\rm maxsing}^c) 
				\nonumber\\
				&\ge 1-r^2\exp\sth{-10\pth{\log(mLr^2)+(\log m)(\log r)}}-2mLe^{-cn}
				\overset{(a)}{\ge} 1-\frac 1{(rmL)^{8}}-2mLe^{-cn}.
			\end{align*}
			On the above high-probability event $\Ec_{\rm dcpl}\cap \widetilde\Ec_{\rm maxsing}$, we have for each $1\leq i<j\leq r$
			\begin{align*}
				& \sum_{l=1}^L \nmhs{A_l(X_i^2-X_j^2)A_l^\top }^2 
				\leq Lmt_{1,i,j}^2+Lm(m-1)\|\bd_{i,j}\|_{2}^2 + t_{2,i,j}\\
				&\overset{(a)}{\le} C L m x_{\max}^2 \nm{\xv_i-\xv_j}_{2}^2\log(mLr^2)+CLm \nm{\xv_i-\xv_j}_{\infty}^2\log^2(mLr^2)
				+ 4Lm(m-1)x_{\max}^2\|\xv_i-\xv_j\|_2^2\\
				&\quad +C\log m\pr{x_{\max}^2 \nm{\xv_i-\xv_j}_{4}^2 \sqrt{mL(\sqrt m+\sqrt n)^4\log r^2}+\nm{\xv_i-\xv_j}_{\infty}^2 (\sqrt n +\sqrt m)^2 \log r^2 } \\
				&\overset{(b)}{\le} C L m x_{\max}^2 \frac{kn}{N_{\rm div}}\de_r^2\log(mLr^2)+CLm \de_r^2\log^2(mLr^2)
				+ 4Lm(m-1)x_{\max}^2\frac{kn}{N_{\rm div}}\de_r^2\\
				&\quad +C\log m\pr{x_{\max}^2 \sqrt{\frac{kn}{N_{\rm div}}}\de_r^2 \sqrt{mL(\sqrt m+\sqrt n)^4\log r^2}+\de_r^2 (\sqrt n +\sqrt m)^2 \log r^2 },
			\end{align*}
			where $(a)$ followed by using the inequality $(a+b)^2\le 2a^2+2b^2$, and for $(b)$ we have used the fact that by our construction of $\mathcal X_k^{\rm finite}$, $\xv_i-\xv_j$ is $\frac{2kn}{N_{\rm div}}$-sparse for any $\xv_i, \xv_j\in \calS_{\rm sep}$, and thus $\nm{\xv_i-\xv_j}_{\infty}\le \de_r$, $\nm{\xv_i-\xv_j}_2^2\le \frac{kn}{N_{\rm div}}\de_r^2$, and $\nm{\xv_i-\xv_j}_4^2\le \sqrt{\frac{kn}{N_{\rm div}}}\de_r^2$.
			
			\section{Proof of results in \prettyref{sec:proof of lower bound oversample}}
			\label{app:lower-bound-oversample}
			
			\subsection*{Proof of \prettyref{prop:sufficient}}
			$\bT_{A}(\ly)$ is an one to one transformation of $\ly$ whenever $ \bA^\top \bA$ is invertible, as $\bA \bT_{ \bA}(\ly)=\ly$. As any one to one transformation of a sufficient statistic is also a sufficient statistic, and $\ly$ is itself a sufficient statistic, we get that $\bT_{\bA}(\ly)$ is a sufficient statistic.

			\subsection{Proof of \prettyref{lmm:lower bound for the improved estimator}}
			By our construction, we have that $\bT_{\bA}(\ly)
			=[\uv_1^\top ,...,\uv_L^\top ]^\top$, where $\uv_l=(A_l^\top  A_l)^{-1}A_l\yv_l\in \R^n$. In the above optimization problem, as the estimators we consider are all of the form $g(\bT_{\bA}(\ly))$, we may treat $\uv_1,...,\uv_L$ as our observations and the problem transformed into recovering $\xv_o\in \Cc_k$ from the simplified model
			\begin{equation*}
				\uv_l=X_o \wv_l, \quad \bw_l\sim N(0,I_n),l=1,2,...,L.
			\end{equation*}
			Hence, we will apply \prettyref{lmm:fano}, to derive the lower bound, with the new data distribution
			\begin{equation*}
				\begin{gathered}
					\P_{\xv} \sim \otimes_{l=1}^L N(\boldsymbol{0},\Sigma_l^{{-1}}(\xv)) = N(\boldsymbol{0},\Sig^{{-1}}(\xv)),
					\\
					\Sig(\xv):=\diag\pr{\Sigma_1(\xv),...,\Sigma_L(\xv)},
					\quad
					\Sigma_l=\Sigma_l(\xv):=(\sigma_z^2 I_n + X^2)^{-1}.
				\end{gathered}
			\end{equation*}
			The proof of lower bound uses \prettyref{lmm:fano} and a similar discretization and construction of $\alpha_r$-separated set $\{\xv_1,...,\xv_r\}$ as in \prettyref{sec:discretization}, with $N_{\rm div}=kn^{\e}$ and the choice for $\delta_r^2,\alpha_r^2$
			\begin{equation}\label{eq:new choice of delta}
				\de_r^2:=c_{\de}\frac{\log r}{nL}\cdot \frac{N_{\rm div}}{k}, \text{ and } \frac{\alpha_r^2}{n}:=\frac{k}{N_{\rm div}}\de_r^2,
			\end{equation}
			for a small constant $c_\delta>0$. Denote $\PP_{\bx_i}:=\PP_i$ for $i\in \{1,\dots,r\}$. As our earlier construction of $\{\bx_1,\dots,\bx_r\}$ implies $\frac{kn}{N_{\rm div}}$ and $\xv_i-\xv_j$ is $\frac{2kn}{N_{\rm div}}$-sparse for any $\xv_i\neq \xv_j$, which implies $\nm{\xv_i-\xv_j}_2^2\le \frac{2kn}{N_{\rm div}}\de_r^2$. Hence, using \prettyref{lmm: upper bound for KL divergence}, with $A_l=I_n,\emax=\emin=1$, we get
			\begin{equation}\label{eq:m6}
				\begin{split}
					&~\beta_r:=\max_{1\le i<j\le r}\KL(\P_i\parallel \P_j)
					\le  L\frac{\xmax^4}{\xmin^8}
					\max_{1\le i<j\le r}\nmhs{X_i^2-X_j^2}^2 \\
					& \le  \frac{4\xmax^6}{\xmin^8} L\max_{1\le i<j\le r} \nm{\xv_i-\xv_j}_2^2 
					\le \frac{4\xmax^6}{\xmin^8} \frac{Lkn\de_r^2}{N_{\rm div}}\le \frac{1}{10} \log r,
				\end{split}
			\end{equation}
			for sufficiently small $c_{\de}$ in \eqref{eq:new choice of delta} to guarantee a small $\delta_r^2$. Now it follows from \prettyref{lmm:fano} that
			\begin{equation}\label{eq:m7}
				\begin{split}
					\inf_{\bde}\sup_{\xv\in \Cc_k}\E\br{\left.\frac{\nm{g(\bT_{\bA}(\ly))-\xv_o}_2^2}{n} \right \vert \esingp} 
					\ge \frac{\alpha_r^2}{4n}\pr{1-\frac{\beta_r+\log 2}{\log r}}^2 \PP[\esingp]
					=\Theta_{\e,x_{\min},x_{\max}}\pr{ \frac{k\log n}{nL}}.
				\end{split}
			\end{equation}
			where the last equality uses \eqref{eq:new choice of delta}, with ${N_{\rm div}\over k}=n^\e,
			{\alpha_r^2\over n}
			={k\de_r^2 \over N_{\rm div}}
			=\Theta_{x_{\min},x_{\max}}\pr{ \frac{k\log \pr{N_{\rm div} \over k}}{nL}}$.

			\section{Proof of technical results in \prettyref{sec:upper bound}}\label{app:technical}
			
			\subsection{Proof of \prettyref{lmm:m-error-union-bound}}
			
			As stated from the beginning of \prettyref{sec:General strategy}, without loss of generality, we may assume $\Cc\in \mathcal{F}_{a,b,k,n}$ with $a=\xmax-\xmin$ and $b=1$. Namely \begin{align}\label{eq:covering-bound}
				|\netdelta|\le \pr{(\xmax-\xmin) n \over \delta_{\rm net}}^k.
			\end{align}
			
			Fix a general $\xvt_o\in \netdelta$. For $1\le l \le L$, define the matrices
			\begin{align}
				\begin{gathered}
					A_{\sigma_z, l}
					=\qth{
						\sigma_z I_m ~ A_l X_o}
					\in \reals^{m\times (m+n)},
					\quad
					B_l=A_{\sigma_z, l}^\top \pr{\Sigma_l(\xvt_o)-\Sigma_l(\xv_o)}A_{\sigma_z, l}
					\in \reals^{(m+n)\times (m+n)},
					\\
					\bA_{\sigma_z}
					=\diag(A_{\sigma_z, 1},\dots, A_{\sigma_z, L})\in \reals^{mL\times (m+n)L},
					\quad
					\bB = \diag(B_1,\dots,B_L)\in \reals^{L(m+n)\times L(m+n)}.
				\end{gathered}
			\end{align}
			In view of the notations \eqref{eq:So} and \eqref{eq:m-Sotil}, the above display provides us with the following identities 
			\begin{align}
				\bA_{\sigma_z}\bA_{\sigma_z}^\top =\So^{-1},\quad
				\bB = \bA_{\sigma_z}^\top \pr{\Sotil-\So}\bA_{\sigma_z}.
			\end{align}
			Define the $L(m+n)$ dimensional vector, $\lw^\top  = [\zv_1^\top /\sigma_z,\wv_1^\top ,\ldots, \zv_L^\top /\sigma_z, \wv_L^\top ]$. It follows that
			\begin{align}
				\ly^\top  \pr{\Sotil-\So} \ly=\lw^\top  \bB \lw.
			\end{align}
			Then, conditioning on $A_1,...,A_L$, by the Hanson-Wright inequality (\prettyref{lmm:hanson wright}), we have 
			\begin{align}\label{eq:m11}
				\Prob \qth{|\lw^\top \bB\lw- \Tr(\So^{-1}(\Sotil-\So))| >t} \leq 2\exp\left(-c\min \left({t^2\over 4\|\bB\|_{\rm HS}^2 },{ t \over 2\|\bB\|_2 } \right)\right).
			\end{align}
			We simplify the above using upper bounds on $\|\bB\|_{\rm HS},\|\bB\|_2$. In view of \eqref{eq:cnms} we bound $\|\bB\|_2$ as
			\begin{align*}
				\begin{split}
					& \|\bB\|_2 = \max\pr{\|B_l\|_2}_{l=1}^L\\
					\le & \pr{\sigma_z^2+x_{\max}^2\max_{1\le l \le L} \smax(A_l)^2}\max_{1\le l \le L} \smax\pr{\Sigma_l(\xvt_o)-\Sigma_l(\xv_o)} \\
					\le
					& \pr{\sigma_z^2+x_{\max}^2\emax}\max_{1\le l \le L} (\smax(\Sigma_l(\xvt_o))+\smax(\Sigma_l(\xv_o))) \\
					= &
					\pr{\sigma_z^2+x_{\max}^2\emax}\max_{1\le l \le L} (\{\smin(\Sigma_l(\xvt_o)^{-1})\}^{-1}+\{\smin(\Sigma_l(\xv_o)^{-1})\}^{-1})
					\le 
					\cnms,
				\end{split}
			\end{align*}
			where the last inequality followed by noting that on the event $\calS$ in \eqref{eq:eigen-conditioning} we have for each $\xv=\xv_o,\xvt_o$
			$$
			\smin(\Sigma_l(\xv)^{-1})
			=\sigma_z^2+\sigma_{\min}(A_lX^2A_l^\top)
			\geq 
			\sigma_z^2+x_{\min}^2\sigma_{\min}(A_lA_l^\top)
			\geq 
			\sigma_z^2+x_{\min}^2\emin,
			\quad 1\leq l\leq L.
			$$
			Next, using the identity $\|\bB\|_{\rm HS}^2=\Tr(\bB^2)$ we get 
			\begin{align}
				\|\bB\|_{\rm HS}^2 =\trbr{\So^{-1}\pr{\Sotil-\So}\So^{-1}\pr{\Sotil-\So}}.
			\end{align}
			which we have defined as $\mathscr Z$.
			Hence, conditioned on the event $\calS$ in \eqref{eq:eigen-conditioning}, we simplify \eqref{eq:m11} to get
			\begin{align}
				\Prob \qth{\left. |\lw^\top \bB\lw- \Tr(\So^{-1}(\Sotil-\So))| >t \right| \calS} 
				\leq 2\exp\left(-c\min \left({t^2\over 4 \mathscr Z },{t \over 2\cnms} \right)\right).
			\end{align}
			To extend the above probability statement for all possible choice of $\xvt_o\in \netdelta$ we use a union bound argument.
			
			The above implies that by choosing $t = b_1\sqrt {\mathscr Z}+b_1'$ as defined in the result statement, we get
			\begin{align}
				&\Prob \qth{\left.|\lw^\top \bB\lw- \Tr(\So^{-1}(\Sotil-\So))| >t \text{ for some $\xvt_o\in \netdelta$}\right|\calS}
				\nonumber\\
				&\leq 
				e^{k\log((x_{\max}-x_{\min}) n/\delta_{\rm net})} 
				\exp\pth{-c\min\sth{{b_1^2\over 4},{b_1'\over 2\cnms}}}
				\leq e^{-\tilde c k\log((x_{\max}-x_{\min}) n/\delta_{\rm net})},
			\end{align}
			for some constant $\tilde c>0$. In view of \eqref{eq:eigen-conditioning} and $\PP[A]\leq \PP[A|\calS]+\PP[\calS^c]$, the above display implies

			\begin{align*}
				&\Prob \qth{|\lw^\top \bB\lw- \Tr(\So^{-1}(\Sotil-\So))| >t \text{ for some $\xvt_o\in \netdelta$}}
				\leq 
				e^{-\tilde c k\log\pth{(x_{\max}-x_{\min}) n/\delta_{\rm net}}}
				+Le^{-cn}.
			\end{align*}

			%
			%
			\subsection{Proof of \prettyref{lmm:m-bound-deltanet-error}}
			
			Our entire argument is conditioning on the high-probability event $\calS$ in \eqref{eq:eigen-conditioning}. We will explain at the end of the section how the conditioning is removed to get the final result. In view of \prettyref{lmm:boundeigenvalues} we first note that
			\begin{align}\label{eq:m13}
				\sigma_{\max}(\Sotil-\Sohat)
				&\leq 
				\frac{\sigma_{\max} \pr{\Sohat^{-1} -\Sotil^{-1}}}{\sigma_{\min}\pr{\Sohat^{-1}} \sigma_{\min}\pr{\Sotil^{-1}} }
				\leq 
				\frac{\max_{1\le l\le L}\br{\sigma_{\max} (A_l(\widehat{X}_o^2-\widetilde{X}_o^2) A_l^\top) }}{\sigma_{\min}\pr{\Sohat^{-1}} \sigma_{\min}\pr{\Sotil^{-1}}}
				\nonumber \\
				&\leq 
				\frac{ \max_{1\le l\le L}\lmax(A_lA_l^\top )  \|{\xvh_o}^2 - {\xvt}_o^2\|_\infty}{\left(\sigma_z^2+x^2_{\min} \min_{1\le l \le L}\lambda_{\min} (A_lA_l^\top )\right)^2}
				\le \frac{2\emax x_{\max}\delta_{\rm net}}{(\sigma_z^2+x_{\min}^2\cdot \emin)^2}.
			\end{align} 
			In view of $\sigma_{\max}(\So^{-1})\leq \sigma_z^2 + x_{\max}^2\emax$, for any $1\le q \le mL$ we use the last display to get
			\begin{align}
				\left|\lambda_q (\So^{-1}(\Sohat-\Sotil))\right|
				\leq 
				\frac{\sigma_{\max} \pr{\Sotil-\Sohat}}{\sigma_{\min}(\So)} 
				\le \frac{[\sigma_z^2+x_{\max}^2\emax]\cdot \emax}{\pr{\sigma_z^2+x_{\min}^2\cdot \emin}^2}\cdot 2x_{\max}\delta_{\rm net}, 
			\end{align}
			with probability $1-Le^{-cn}$.
			Summing up over $1 \leq q\leq mL$ and using triangular inequality, we have
			\begin{align}\label{eq: sigma inverse times difference of sigmas}
				\abs{\Tr[\So^{-1}(\Sohat-\Sotil)]}
				\le 
				mL\cdot \frac{[\sigma_z^2+x_{\max}^2\cdot \emax]\cdot \emax}{\pr{\sigma_z^2+x_{\min}^2\cdot \emin}^2}\cdot 2x_{\max}\delta_{\rm net}.
			\end{align}
			To bound the term $|\ly^\top(\Sohat-\Sotil)\ly|$, we again use \eqref{eq:m13} to get
			\begin{align}\label{eq:y sigma y}
				\abs{\ly^\top \pr{\Sohat-\Sotil}\ly}
				\le  \sigma_{\max}\pr{\Sohat-\Sotil} \ly^\top  \ly
				\le  
				\frac{2\emax x_{\max}\delta_{\rm net}}{\pr{\sigma_z^2+x_{\min}^2\cdot \emin}^2} \ly^\top  \ly.
			\end{align}
			We use the following lemma to bound $\ly^\top  \ly$. A proof is provided at the end of this section. 
			
			\begin{lemma}\label{lmm:bound-yyt}
				For $t>0$, we have 
				\begin{align}
					\begin{split}
						&~\Prob\pr{\ly^\top  \ly\ge mL \br{\emax x_{\max}^2+\sigma_z^2}+t }\\
						&\le
						2\exp\left(-c\min \left({t^2\over 4mL\left(\sigma_z^2 + x_{\max}^2\cdot \emax \right)^2 },{t \over 2\left(\sigma_z^2 + x_{\max}^2 \cdot \emax \right) } \right)\right)+Le^{-cn}.
					\end{split}
				\end{align}
			\end{lemma} 
			We choose
			$
			t=CmL\br{\sigma_z^2 + x_{\max}^2 \cdot \emax}$ 
			for a large constant $C>0$. Then, in view of \prettyref{lmm:bound-yyt}, we we continue \eqref{eq:y sigma y} to get that with probability $1-Le^{-cn}-\exp\pr{-cmL}$ the following holds
			\begin{align}
				\abs{\ly^\top \pr{\Sohat-\Sotil}\ly} 
				\le \frac{2x_{\max}\delta_{\rm net} \emax}{\pr{\sigma_z^2+x_{\min}^2\cdot \emin}^2} \cdot \breve C mL\br{\sigma_z^2 + x_{\max}^2 \cdot \emax},
			\end{align}
			for a constant $\breve C>0$. Together with \eqref{eq: sigma inverse times difference of sigmas}, we conclude that with probability $1-Le^{-c_1n}-e^{-c_2mL}$ that
			
			\begin{align}\label{definition of b_2}
				\begin{split}
					\abs{\ly^\top  (\Sotil-\Sohat) \ly - \Tr(\So^{-1}(\Sotil-\Sohat))}
					\le
					\breve C'mLx_{\max}\delta_{\rm net}\cdot \frac{[\sigma_z^2+x_{\max}^2\cdot \emax]\cdot \emax}{\pr{\sigma_z^2+x_{\min}^2\cdot \emin}^2},
				\end{split}
			\end{align}
			for a constant $\breve C'>0$. As $mLx_{\max}\delta_{\rm net}\cdot \frac{[\sigma_z^2+x_{\max}^2\cdot \emax]\cdot \emax}{\pr{\sigma_z^2+x_{\min}^2\cdot \emin}^2}
			\leq mL\delta_{\rm  net}(\cnms)^2$ we get the desired result. 
			
			\begin{proof}[Proof of \prettyref{lmm:bound-yyt}]
				
				Define
				\begin{align}
					\begin{gathered}
						M_l=A_{\sigma_z,l}^\top  A_{\sigma_z,l},	\quad
						l\in [L], \quad
						\bM =\diag (M_1,\dots, M_L)\in \mathbb{R}^{L(m+n)\times L(m+n)}\\
						\lw^\top  = [\zv_1^\top /\sigma_z,\wv_1^\top , \ldots, \zv_L^\top /\sigma_z,\wv_L^\top ]
					\end{gathered} 
				\end{align}
				Then we can write $\ly^\top  \ly= \lw^\top  \bM  \lw^\top $. To obtain a tail bound, we first recall $\ly = [\yv_1,\dots, \yv_L]$, and observe that by conditioning on the event $\calS$ as in \eqref{eq:eigen-conditioning}, we have for each $l=1,\dots, L$,
				\begin{align}
					\begin{split}
						\E[\yv_l^\top  \yv_l \mid \calS]
						= \trbr{X_oA_l^\top A_lX_o} + m \sigma_z^2
						\le m \pr{\lmax(A_l^\top A_l) x_{\max}^2+\sigma_z^2}
						\le m \pr{\emax x_{\max}^2+\sigma_z^2}.
					\end{split}
				\end{align}
				Then, by the Hanson-Wright inequality (\prettyref{lmm:hanson wright}), we have conditioned on the event $\calS$ that,
				\begin{align}
					\Prob \pr{\left.\ly^\top  \ly> mL \pr{\emax x_{\max}^2+\sigma_z^2}+t
						\right|\calS}
					\leq & 2\exp\left(-c\min \left({t^2\over \|\bM \|_{\rm HS}^2 },{ t \over \|\bM \|_2 } \right)\right),
				\end{align}
				for some constant $c>0$. To bound the terms $\|\bM\|_{\rm HS},\|\bM\|_2$, we note that on the event $\calS$, 
				\begin{align}
					\begin{gathered}
						\|\bM \|_2  = \max_{1\le l\le L} \|M_l\|_2
						\le \sigma_z^2 + x_{\max}^2 \max_{1\le l \le m}\lmax(A_l^\top  A_l)
						\leq \sigma_z^2 + x_{\max}^2 \emax,
						\\
						\|\bM \|_{\rm HS}^2
						=\sum_{l=1}^L {\rm Tr}(M_l^2)\le \sum_{l=1}^L\sum_{i=1}^m\lambda_i^2(M_l) \leq \sum_{l=1}^L m\lambda^2_{\max}(M_l)
						\le Lm \left(\sigma_z^2 + x_{\max}^2 \emax \right)^2,
					\end{gathered}
				\end{align}
				where the inequality $\sum_{l=1}^L {\rm Tr}(M_l^2)\le \sum_{l=1}^L\sum_{i=1}^m\lambda_i^2(M_l)$ follows from the fact that ${\rm rank}(M_l)= m$ and there are at most $m$ nonzero eigenvalues. Finally we remove the condition on $\calS$ by using $\PP[A]\leq \PP[A|\calS]+\PP[\calS^c]$ with $\PP[\calS^c]\leq Le^{-cn}$ from \eqref{eq:prob-S} to get the desired result.
				
			\end{proof}

			\subsection{Proof of \prettyref{lmm:bound-fbar-diffs}}
			
			\label{app:proof-of-bound-fbar-diffs}	
			
			The proof follows that of \cite[Lemma VIII.4]{zhou2022compressed} with minor differences, which we clarify below.
			Recall, 
			\begin{align*}
				\overline{f}(\xvt_o)-\overline{f}(\xv_o)
				=& \br{-\log\det \Sotil +{\rm Tr}\pr{\Sotil \So^{-1}}}-\br{-\log\det \So +{\rm Tr}\pr{\So \So^{-1}}}\\
				=&-\br{\log \det \Sotil - \log \det \So}+\Tr([\Sotil-\So]\So^{-1}).
			\end{align*}
			For $1\le q \le mL$, let $\la_q$ denote the $q$-th eigenvalue of $\So^{-\frac{1}{2}}\pr{\Sotil-\So}\So^{-\frac{1}{2}}$, and $\la_{\max}=\max_{1\le q\le mL} |\lambda_q|$. Following the proof strategy of \cite[Lemma VIII.4]{zhou2022compressed} we can show that
			\begin{align}
				\log\frac{\det \Sotil}{\det \So}
				\leq  \Tr\pr{\So^{-\frac{1}{2}}\pr{\Sotil-\So}\So^{-\frac{1}{2}}}
				-\frac{{\rm Tr}\left(\So^{-1}(\Sotil-\So)\So^{-1}(\Sotil-\So)\right)}{2(1+\tilde \la_{\max})^2}.
			\end{align}
			We will use the following inequalities to bound $\widetilde\la_{\max}$, on the event $\calS$.
			\begin{itemize}
				
				\item $ \|\So^{-1}\|_2=\max_{1\le l \le L}\|\sigma_z^2 I_m+ A_l X_o^2 A_l^\top \|_2
				\leq \sigma_z^2+x^2_{\max}\max_{1\le l \le L}\smax(A_lA_l^\top )
				\lesssim \sigma_z^2 + x_{\max}^2 E_{\max}$. 
				
				\item Using $\|E\|^{-1}=\smin(E^{-1})$ for any invertible matrix $E$, for a constant $C$, we get
				\begin{align*}
					\|\So\|^{-1}& =  \sigma_z^2+\min_{1\le l \le L}\smin(A_l X_o^2 A_l^\top)
					\geq \sigma_z^2+x^2_{\min}\min_{1\le l \leq L}\smin(A_lA_l^\top )
					\geq \sigma_z^2 + x_{\min}^2 \emin;\\
					\|\Sotil\|^{-1}&= \sigma_z^2+\min_{1\le l \le L} \smin(A_l \tilde X_o^2 A_l^\top)
					\geq \sigma_z^2+x^2_{\min}\min_{1\le l \le L}\smin(A_lA_l^\top )
					\geq \sigma_z^2 + x_{\min}^2 \emin.
				\end{align*}
			\end{itemize}
			Using $\widetilde \lambda_{\max}=\|\So^{-{1\over 2}}(\Sotil-\So)\So^{-{1\over 2}}\|_2 \leq \|\Sotil-\So\|_2 \|\So^{-{1\over 2}}\|_2^2 = \|\Sotil-\So\|_2 \|\So^{-1}\|_2$,
			we get
			\begin{align}\label{eq:m3}
				\widetilde \lambda_{\max}
				&\leq (\|\Sotil\|_2 + \|\So\|_2)\|\So^{-1}\|_2
				\leq 2{\sigma_z^2 + x_{\max}^2 \emax\over \sigma_z^2 + x_{\min}^2 \emin} \lesssim \cnms.
			\end{align}

			\subsection{Proof of \prettyref{lmm:f-sotil-f-sohat-bound}}
			\label{app:proof-f-sotil-f-sohat-bound}
			Note that we will use these results in the subcase $n\geq 4m$, and the subcase $n<4m$ with $\sigma_z^2\geq m$. Then we have $\cnms$ is of a constant order. Define $|\lambda_{\max}|$ as the largest absolute value of eigenvalues of $\Sohat^{-{1\over 2}}(\Sotil-\Sohat)\Sohat^{-{1\over 2}}$. Then, using $|\lambda_{\max}|\le \|\Sotil-\Sohat\|_2 \|\Sohat^{-1}\|_2$, with probability $1-Le^{-cn}$ we get
			\begin{align*}
				|\lambda_{\max}|
				&\overset{(a)}{\le} \frac{\smax(\Sohat^{-1}-\Sotil^{-1})\|\Sohat^{-1}\|_2}{\smin(\Sohat^{-1})\smin(\Sotil^{-1})}\\
				&\le {\sigma_z^2+x^2_{\max}\cdot \emax\over {\sigma_z^2+x^2_{\min}\cdot \emin}^2}\max_{1\le l\le L}\smax\pr{A_l(\Xto^2-\Xho^2) A_l^\top}
				\le \cnms\cdot x_{\max} n \de_{\rm net}
				\ll 1
			\end{align*}
			where we use $\de_{\rm net}={x_{\max}\over n^5}$, $n^4k\log n\geq mL$ from the theorem statement, and $(a)$ follows from \prettyref{lmm:boundeigenvalues}. 
			
			We now get back to the proof of the main result. For the first term $\overline{f}(\xvh_o)-\overline{f}(\xvt_o)$, using an argument similar to the proof of \prettyref{lmm:bound-fbar-diffs} we can show
			\begin{align}
				&{\overline{f}(\xvt_o)-\overline{f}(\xvh_o)}
				=
				-\br{\log \det \Sotil - \log \det \Sohat}+\Tr([\Sotil-\Sohat]\So^{-1})
			\end{align}
			The above implies 
			\begin{align}
				\abs{{\overline{f}(\xvt_o)-\overline{f}(\xvh_o)}}
				&\le {1\over 2(1-|\lambda_{\max}|)^2}{\rm Tr}\left(\Sohat^{-1}(\Sotil-\Sohat)\Sohat^{-1}(\Sotil-\Sohat)\right)
				\notag\\
				&\quad +\abs{\Tr([\Sotil-\Sohat]\So^{-1})}
				+\abs{\Tr([\Sotil-\Sohat]\Sohat^{-1})}.
			\end{align}
			In view of \eqref{eq: sigma inverse times difference of sigmas}, noting that $\sigma_{\max}(\Sohat^{-1})\leq \sigma_z^2 + x_{\max}^2\emax$ we get for large enough $n$ such that $n^4 k\log n\geq mL$, we get that $\delta_{\rm net}$ is small enough to satisfy
			\begin{align*}
				\max\sth{\abs{\Tr([\Sotil-\Sohat]\So^{-1})}
				,\abs{\Tr([\Sotil-\Sohat]\Sohat^{-1})}}
				\leq (\cnms)^2mLx_{\max}\delta_{\rm net}\ll 1.
			\end{align*}
			Similarly, the bound on ${\rm Tr}\left(\Sohat^{-1}(\Sotil-\Sohat)\Sohat^{-1}(\Sotil-\Sohat)\right)$ also follows. This completes the proof.
			
			\subsection{Proof of \prettyref{lmm:decoupling-lower-bound}}
			\label{app:section on the decoupling in the upper bound}
			{
				 
				Pick $\xvt_o\in \netdelta$. Define $\bd=\xvt_o^2-\xv_o^2$ and $D=\diag(\bd)$. In view of the above, we have
			\begin{align}
				\| \bd\|_{\infty}\le x_{\max}^2, \quad \| \bd^2\|_{2}\le 4x_{\max}^3\|\xvt_o-\xv_o\|_2.
			\end{align}
			Using \prettyref{lmm:decoupling-genr} part (2), conditional on the event $\calS$ in \eqref{eq:eigen-conditioning} and $\breve c>0$ depending only on $x_{\max},K$
			\begin{align*} 
					&~\Prob \pr{\br{\sum_{l=1}^L\|A_l (\widetilde{X}_o^2-X_o^2)A_l^\top \|^2_{\rm HS} < Lm(m-1) \sum_{i=1}^n d_i^2 - t}\Biggr|\calS}
					\\
					&\leq \exp \left( -\breve c \cdot \min \pr{\frac{t^2}{ \|\xvt_o- \xv_o\|_{2}^{2}  mL(\sqrt n+\sqrt m)^4}, \frac{t}{ (\sqrt n+\sqrt m)^2} }\right)\\
					&\quad+ 2 m  \exp\left(-\breve c\min\left({t^2\over L m^3 \|\xvt_o- \xv_o\|_{2}^{2}},{t \over m  }\right)\right).
			\end{align*}
			Define 
			$$
			t_o= \breve C_{x_{\max}} \max(m,n) \pth{  \|\xvt_o- \xv_o\|_{2} \sqrt{mL} \sqrt{k \log \frac{(\xmax-\xmin)n}{\delta_{\rm net}}}
				+ k \log \frac{(\xmax-\xmin)n}{\delta_{\rm net}}},$$ 
			for a large constant $\breve C_{x_{\max}}$ as required. Then, the above display implies for a large constant $C_1$
			\begin{align*}
				\PP\qth{\br{\sum_{l=1}^L\|A_l (\widetilde{X}_o^2-X_o^2)A_l^\top \|^2_{\rm HS} < Lm(m-1) \sum_{i=1}^n d_i^2 - t_0}\Biggr|\calS}
				\leq e^{- C_1 k \log {(\xmax-\xmin)n\over \delta_{\rm net}}}
			\end{align*}
			Then using an union bound over the total possible choices of $\xvt_o\in \netdelta$, with $|\netdelta|\leq \pth{(x_{\max}-x_{\min})\sqrt n\over \delta_{\rm net}}^k$ as in \eqref{eq:covering-bound}, we get that as $C_1$ is large enough, for some large constant $C_2>0$
			\begin{align*}
				\PP\qth{\br{\sum_{l=1}^L\|A_l (\widetilde{X}_o^2-X_o^2)A_l^\top \|^2_{\rm HS} < Lm(m-1) \sum_{i=1}^n d_i^2 - t_o}
					\text{ for any $\xvt_o\in \netdelta$}\Biggr|\calS}
				\leq e^{-C_2k\log n}.
			\end{align*}
			Using $\PP[B]\leq \PP[B|\calS]+\PP[\calS^c]$, with $B$ as the event $\sum_{l=1}^L\|A_l (\widetilde{X}_o^2-X_o^2)A_l^\top \|^2_{\rm HS} < Lm(m-1) \sum_{i=1}^n d_i^2 - t_o$, and $\PP[\calS^c]\leq Le^{-cn}$ as in \eqref{eq:prob-S} we conclude $\PP[B]\leq e^{-ck\log n}+Le^{-cn}$. This implies our the desired result.
			
			}
			
			%
			%

			\section{Comparison to the fixed forward operator model}\label{app:Comparison to the unvarying measurement model}

			\subsection{Proof of \prettyref{thm:main-genr2}}
			
			
			The proofs in this section uses a similar approach to the proofs in \prettyref{sec:upper bound}, with the key modifications $A_1=\cdots=A_L=A$. To proceed with the details, we first define the notations we use throughout the section, and then point out the differences with the proofs in \prettyref{sec:upper bound}. The proof of the related technical results, particularly \prettyref{lmm:m-error-union-bound fixed}, \prettyref{lmm:m-bound-deltanet-error fixed}, \prettyref{lmm:bound-fbar-diffs fixed}, \prettyref{lmm:f-sotil-f-sohat-bound fixed}, \prettyref{lmm:decoupling-lower-bound fixed}, follow from the proofs of the results in \prettyref{sec:upper bound} by noting that the related proof in the multilook setting uses bounds on the singular values of the sensor matrices $A_1,\dots, A_L$, for which we used a common bound that also holds true for the fixed measurement matrix $A$. In addition, of results in the fixed sensor case provide guarantees with a higher probability as we do not need to have a uniform control of the singular values of the sensor matrix, as required in the independent multilook setup. This will also hold true for the subsequent results. We omit the technical details. 
			
			Similar to before, we will show the desired upper bound is achieved by the maximum likelihood estimator, over is the class of all possible signals $\calC$
			\begin{align}\label{eq:log-likelihood fixed}
				\xvh_o = \arg\min_{\xv \in \Cc} f(\xv),
				\quad 
				f (\xv) 
				= L \log\det \left(\sigma_z^2 I_m +  A X^2 A^\top \right) + \sum_{l=1}^L \yv_l^\top  \left( \sigma_z^2 I_m+ A X^2 A^\top \right)^{-1}  \yv_l
			\end{align}
			For the entirety of the analysis in this section, we will restrict ourselves to the following event 
			\begin{align}\label{eq:eigen-conditioning fixed}
				\calS_{\rm fix}=\sth{A:\sigma_{\min}(AA^\top)\geq E_{\min},\quad  \sigma_{\max}(AA^\top) \leq \emax},
			\end{align}
			where $\emax$ and $\emin$ according to the following rules
			\begin{itemize}
				\item Case I ($n\geq 4m$): We will choose $\emax=\frac 94(\sqrt n+\sqrt m)^2$ and $\emin = \frac 14(\sqrt n-\sqrt m)^2$. In that case the event $\calS_{\rm fix}$ satisfies $\esing\subseteq \calS_{\rm fix}$, where $\esing$ is given as in \eqref{eq:esing}. This implies $\PP\qth{\calS_{\rm fix}}\geq \PP[\esing]\geq 1-e^{-cn}$ for some constant $c>0$.
				\item Case II ($n< 4m$): We will choose $\emax=\frac 94(\sqrt n+\sqrt m)^2$ and $\emin = 0$. In that case the event $\calS_{\rm fix}$ satisfies $\emaxsing\subseteq \calS_{\rm fix}$, where $\emaxsing$ is given as in \eqref{eq:emaxsing}. This implies $\PP\qth{\calS_{\rm fix}}\geq \PP[\emaxsing]\geq 1-e^{-cn}$ for some constant $c>0$.
			\end{itemize}
			Consider the following notations for simplifying the presentation. Let $\Sigma$ be the inverse of the covariance matrix $\EE\qth{\yv_l \yv_l^\top|A}$ given by
			\begin{equation*}
				\Sigma=\Sigma(\xv):=(\sigma_z^2 I_m + A X^2 A^\top )^{-1}.
			\end{equation*} 
			Define the vector $\ly\in \reals^{mL}$ and block-diagonal matrix $\Sig(\xv)\in \reals^{mL\times mL}$ as the collection of all the observations and the inverse covariance matrices over different looks
			\begin{align}\label{eq:fixed A So}
				\ly^\top :=(\yv^\top_1,...,\yv^\top_L )^{\top},
				\quad \Sig(\xv):=\diag\pr{\Sigma(\xv),...,\Sigma(\xv)},
				\quad \So=\bSigma(\bx_o),
				\quad
				\Sohat = \bSigma(\hat \bx_o).
			\end{align}
			In view of the above notations, we can rewrite the negative log-likelihood in \eqref{eq:log-likelihood fixed} as
			\begin{align}
				f(\xv)= -\log \det(\Sig(\xv))+ \ly^\top  \Sig(\xv)\ly
			\end{align}
			Now we proceed with the proof. Since $\xvh_o$ is the minimizer from \eqref{eq:log-likelihood fixed}, we have
			\begin{align}\label{eq:1 fixed}
				f(\xvh_o)\leq f(\xv_o).
			\end{align}
			For a fixed $\bx$, define $\overline f(\bx)$ as the function of conditional expectation of $f(\bx)$ given $A$
			\begin{align}
				\overline{f}(\bx):=\E[f(\bx) \mid A]
				= {-\log\det \Sig(\xv) +{\rm Tr}\pr{\Sig(\xv) \Sig(\xv_o)^{-1}}}.
			\end{align}
			Simplifying the expression for $f(\xvh_o)-f(\xv_o)$, with the above notations, we get
			\begin{align}
				& ~f(\xvh_o)-f(\xv_o) 
				\nonumber\\
				&= \ly^\top  \pr{\Sohat-\So} \ly
				-\Tr \qth{\So^{-1}\pr{\Sohat-\So}}+\Tr \qth{\So^{-1}\pr{\Sohat-\So}}-\log \det(\Sohat)+\log \det(\So) 
				\nonumber\\
				&= \ly^\top\pr{\Sohat-\So}\ly - \Tr \qth{\So^{-1} \pr{\Sohat-\So}}+ \overline{f}(\hat \bx_o)-\overline{f}(\bx_o).
			\end{align}
			Therefore, in view of \eqref{eq:1 fixed} we get
			\begin{align}\label{eq:m15 fixed}
				\ly^\top  \pr{\So-\Sohat} \ly - \Tr(\So^{-1}(\So-\Sohat))\ge \overline{f}(\hat \bx_o)-\overline{f}(\bx_o).
			\end{align}
			
			Our following approach is to find an upper bound for the left side in terms of $\nm{\widehat \xv_o - \xv_o}_2$ and a lower bound for the right side in terms of $\nm{\widehat \xv_o - \xv_o}_2$, and simplify the inequality to get an upper bound for $\nm{\widehat \xv_o - \xv_o}_2$. We will use the following notation, similar to \eqref{eq:cnms}, for a constant $c>0$
			\begin{equation}\label{eq:cnms fixed}
				\cnms = C(n,m,\sigma_z,x_{\max},x_{\min}) = c{\sigma_z^2+x^2_{\max} \emax \over \sigma_z^2+x^2_{\min} \emin }.
			\end{equation}
			
			\noindent \textbf{Establishing an upper bound on $\ly^\top  (\So-\Sohat) \ly - \Tr(\So^{-1}(\So-\Sohat))$:} We use the same $\delta$-net argument as in \prettyref{sec:upper bound}. Consider a $\delta$-net of the set $\calC_k$, denoted by $\netdelta$, with the choice of $\delta_{\rm net}$ to be discussed later. Define $\xvt_o$ as the closest vector in $\netdelta$ to $\bx_o$, i.e.,
			\begin{align}\label{eq:netdelta fixed}
				\xvt_o=\argmin_{\xv\in\netdelta}\|\xvh_o-\xv\|_2.
			\end{align}
			We will use the following notations for the rest of the section
			\begin{align}
				\label{eq:m-Sotil fixed}
				\Sotil=\Sig({\xvt}_o),\quad \Xt_o={\rm diag}({\xvt_o}),
				\quad \Sotil=\bSigma(\xvt),\quad \Xt={\rm diag}({\xvt}),\quad \xvt\in \netdelta.
			\end{align}
			Then in view of triangle inequality we get
			\begin{align}\label{eq:m1 fixed}
				&\abs{\ly^\top  (\So-\Sohat) \ly - \Tr(\So^{-1}(\So-\Sohat))}
				\nonumber\\
				&\leq \abs{\ly^\top  (\Sotil-\So) \ly - \Tr(\So^{-1}(\Sotil-\So))}
				+\abs{\ly^\top  (\Sotil-\Sohat) \ly - \Tr(\So^{-1}(\Sotil-\Sohat))}.
			\end{align}
			
			We use an union bound argument to control the first term above, uniformly over all choices of $\xvt\in \netdelta$. This is done in the following result, which is the fixed $A$ version of \prettyref{lmm:m-error-union-bound}.

			\begin{lemma}\label{lmm:m-error-union-bound fixed}
				Consider the definitions in \eqref{eq:fixed A So} and \eqref{eq:m-Sotil fixed}. There exist constants $c_1,c_2,c_3,c_4>0$ such that
				the following holds with probability $1-e^{-cn}-e^{-c_1k\log ((\xmax-\xmin)n/ \delta_{\rm net})}$
				$$\abs{\ly^\top  (\Stil-\So) \ly - \Tr(\So^{-1}(\Stil-\So))}
				\leq b_1 \sqrt{\mathscr Z} +b_1',\quad
				\text{for all } \xvt_o\in \netdelta,
				$$
				where, with the notation in \eqref{eq:cnms fixed}, $b_1,b_1',\mathscr Z$ are defined as
				\begin{align}
					\begin{gathered}
						b_1= c_3\sqrt{k\log\pr{\frac{(\xmax-\xmin)n}{\delta_{\rm net}}}},
						\quad  b_1'= c_4
						\cnms k\log\pr{\frac{(\xmax-\xmin)n}{\delta_{\rm net}}}\cdot {x_{\max}^2},\\
						\mathscr Z= \Tr(\So^{-1}(\Sotil-\So)\So^{-1}(\Sotil-\So)).
					\end{gathered} 
				\end{align}
			\end{lemma}
			The following result, a counterpart to \prettyref{lmm:m-error-union-bound} for fixed $A$,  controls the final term of \eqref{eq:m1 fixed}. 
			\begin{lemma}\label{lmm:m-bound-deltanet-error fixed}
				Let $\cnms$ be as in \eqref{eq:cnms fixed} and denote $b_2=(\cnms)^2 mL \delta_{\rm net}$. There exist constants $c_1,c_2>0$ such that the following holds with probability $1- e^{-c_1 n} - e^{-c_2mL}$
				\begin{align*}
					&\abs{\ly^\top  (\Sotil-\Sohat) \ly - \Tr(\So^{-1} (\Sotil-\Sohat))}
					\leq b_2.
				\end{align*}
			\end{lemma}
			Combining \prettyref{lmm:m-bound-deltanet-error fixed} with \prettyref{lmm:m-error-union-bound fixed}, in view of \eqref{eq:m1 fixed} we have
			\begin{align}\label{eq:m16 fixed}
				\abs{\ly^\top  (\So-\Sohat) \ly - \Tr(\So^{-1}(\So-\Sohat))}
				\leq 
				b_1 \sqrt{\mathscr Z} +b_1'+b_2
			\end{align}
			
			\noindent \textbf{Establishing a lower bound on $\overline{f}(\xvh_o)-\overline{f}(\xv_o)$:} To find the lower bound, we use the decomposition
			\begin{align}\label{eq:m2 fixed}
				\overline{f}(\xvh_o)-\overline{f}(\xv_o)=\overline{f}(\xvh_o)-\overline{f}(\xvt_o)+\overline{f}(\xvt_o)-\overline{f}(\xv_o),
			\end{align}
			with $\xvt_o$ as in \eqref{eq:netdelta fixed}. The first term, $\overline f(\xvh_o)-\overline f (\xvt_o)$ can be bounded by $\cnms x_{\max} n\delta_{\rm net}$ using the fact that $\tilde \bx_o$ is chosen to be at most $\delta_{\rm net}$ distance away from $\hat\bx_o$. We bound $\overline f(\xvt_o)-\overline f (\xv_o)$ using the following result, which is the fixed $A$ version of \prettyref{lmm:bound-fbar-diffs}.
			\begin{lemma}\label{lmm:bound-fbar-diffs fixed}
				Assume that $\sigma_z^2 I_m+A\widetilde{X}_o^2A^\top$ and $\sigma_z^2 I_m+A X_o^2 A^\top$, are invertible. Then,
				\begin{align}\label{eq:avg-trace-difference fixed}
					\overline{f}(\xvt_o)-\overline{f}(\xv_o)&\geq {1\over 2(1+\widetilde\lambda_{\max})^2}{\rm Tr}\left(\So^{-1}(\Sotil-\So)\So^{-1}(\Sotil-\So)\right),
				\end{align}
				where $\widetilde\lambda_{\max}>0$ is the maximum singular value of $\So^{-{1\over 2}}(\Sotil-\So)\So^{-{1\over 2}}$. Moreover, $\widetilde \lambda_{\max}\leq \cnms$ on the event $\calS_{\rm fixed}$ in \eqref{eq:eigen-conditioning fixed}.
			\end{lemma}
			The following result, a fixed $A$ version of \prettyref{lmm:f-sotil-f-sohat-bound}, controls $|\overline f(\xvh_o)-\overline f (\xvt_o)|$ for a given $\delta_{\rm net}$.
			\begin{lemma}
				\label{lmm:f-sotil-f-sohat-bound fixed}
				$|\overline f(\xvh_o)-\overline f (\xvt_o)|
				\leq \cnms\cdot x_{\max} n \de_{\rm net}
				\ll 1$ with probability $1- e^{-cn}$ for some $c>0$.
			\end{lemma}
			
			Combining the above results, in view of \eqref{eq:m2 fixed} we have, with probability $1-\exp(-cn)$,
			\begin{equation}\label{eq:m14 fixed}
				\overline{f}(\xvh_o)-\overline{f}(\xv_o) \ge \frac{\mathscr Z}{(\cnms)^2}-1.
			\end{equation}
			
					\noindent \textbf{Simplifying the quadratic inequality:} Combining \eqref{eq:m14 fixed} and \eqref{eq:m16 fixed}, in view of \eqref{eq:m15 fixed},
		we have 
		\begin{align}\label{eq:m17 fixed}
			\PP\qth{\frac{\mathscr Z}{(\cnms)^2} \le b_1 \sqrt{\mathscr Z}+b_1'+ b_2+1}
			\geq 1-e^{-c_1k\log ((\xmax-\xmin)n/ \delta_{\rm net})}-e^{-cn}-\exp\pr{-cmL}.
		\end{align}
		
		Rewrite the last inequality as $az^2-bz-c\leq 0$, with $z=\sqrt{\mathscr Z}, a=\frac 1{(\cnms)^2},b=b_1, c=b_1'+b_2+1
		$. Note that by our choice of $b_1,b_1'$ and $b_2$, $\frac{b^2}{a^2}\lesssim \frac{c}{a}$.
    Note $z^2$ is no more than the square of the positive root of $az^2-bz-c = 0$, which implies
		\begin{align}\label{eq:quad-ineq-1 fixed}
			\mathscr Z 
		= z^2
		\leq \pth{b+\sqrt{b^2+4ac}\over 2a}^2
		\leq \frac{2b^2+2(b^2+4ac)}{4a^2}=\frac{b^2}{a^2}+\frac{2c}{a}\le C'\cdot  \frac{c}{a}
		\end{align}
		for some constant $C'>0$. Using the notations from, \eqref{eq:cnms fixed}, \prettyref{lmm:m-error-union-bound fixed} and \prettyref{lmm:m-bound-deltanet-error fixed} we get
		\begin{align}\label{upper bound for Z fixed}
			\mathscr Z
			\le C' \cdot \frac {b_1'+b_2+1}a 
			=
			C'(\cnms)^2\pth{c_3 k\log\pr{\frac{(\xmax-\xmin)n}{\delta_{\rm net}}}\cdot {x_{\max}^2}+(\cnms)^2 mL \delta_{\rm net}+1}.
		\end{align}
		%
		Choose $\delta_{\rm net}=\frac{x_{\max}}{n^5}$ and recall $mL\le n^4 k\log n$ from \prettyref{thm:main-genr2}. Then, from the last display we use \eqref{eq:m17 fixed} to get for a constant $C>0$
		\begin{align}\label{eq:m-scrZ-upperbound fixed}
			\PP\qth{\mathscr{Z}\leq 
				C\cdot (\cnms)^2 k\log n}
			=1-O\pr{n^{-ck}+\exp(-cn)+2\exp\pr{-cmL}}.
		\end{align}
		
		\noindent \textbf{Finding a lowerbound for $\mathscr Z$:}
		In view of \prettyref{lmm:bound-trace}, using the block structure of $\bSigma_o,\Stil$ given in \eqref{eq:fixed A So}, we have on the event $\Ec_{\rm sing}$,
		\begin{align}\label{eq:Z-conjugation-relation fixed}
			\begin{split}
				\mathscr Z
				&
				= \trbr{\So^{-1}\pr{\Sotil-\So}\So^{-1}\pr{\Sotil-\So}}\\
				&  =L\trbr{(\Sigma(\xv_o)^{-1}(\Sigma(\xvt_o)-\Sigma(\xv_o))\Sigma(\xv_o)^{-1}(\Sigma(\xvt_o)-\Sigma(\xv_o))} 
				\\
				&
				\geq 
				{L \over (\cnms)^2 \left(\sigma_z^2+x_{\max}^2 \emax  \right)^2} \nmhs{A(\widetilde{X}_o^2-X_o^2)A^\top }^2,
			\end{split}
		\end{align}
		where $\cnms$ is as in \eqref{eq:cnms fixed}. The lower bound on $\mathscr{Z}$ is completed with the following lower bound on $\nmhs{A(\widetilde{X}_o^2-X_o^2)A^\top }^2$. The result is the fixed $A$ and $L=1$ version of  \prettyref{lmm:decoupling-lower-bound}.
		%
		%
		\begin{lemma}\label{lmm:decoupling-lower-bound fixed}
			 
			Let $n^4k\log n>mL$. Then the following holds with a probability $1-n^{-ck}-mLe^{-cn}$
			\begin{align*}
				\begin{split}
					&\|A(\widetilde{X}_o^2-X_o^2)A^\top \|^2_{\rm HS}\\
					&\ge 4m(m-1) x_{\min}^2 \|\xvt_o- \xv_o\|_2^2  - 
					\breve C_{x_{\max}} \max(m,n) \{\|\xvt_o- \xv_o\|_{2}  \sqrt{m}\sqrt{k \log n}
					+ k \log n\}.
				\end{split}
			\end{align*}
			Here, $\breve C_{x_{\max}}$ is a constant depending only on $x_{\max}$ that grows at most polynomially in $x_{\max}$.
		\end{lemma}

		{

			\noindent \textbf{Final upper bound on $\|\xvt_o- \xv_o\|_2^2$:}
		We combine \eqref{eq:Z-conjugation-relation fixed}, \eqref{eq:m-scrZ-upperbound fixed}, and \prettyref{lmm:decoupling-lower-bound fixed} to summarize the above in terms of the following quadratic inequality with respect to $\|\xvt_o- \xv_o\|_2$, that holds with a probability $1-O\pr{n^{-ck}+\exp(-cn)+\exp\pr{-cmL}}$
		\begin{align}
			\begin{gathered}
				a\|\xvt_o- \xv_o\|_2^2 -b \|\xvt_o- \xv_o\|_2
				-d\leq 0\\
				a={C_1m(m-1)x_{\min}^2\over (\cnms)^2(\sigma_z^2+x_{\max}^2 \emax )^2},
				\quad
				b={\breve C_{x_{\max}} \max(m,n) \sqrt{mk\log n} \over (\cnms)^2(\sigma_z^2+x_{\max}^2 \emax )^2},
				\\
				d= {\breve C_{x_{\max}} \max(m,n) k \log n \over C\cdot (\cnms)^2(\sigma_z^2+x_{\max}^2 \emax )^2}
				+C\cdot (\cnms)^2 {k\log n \over L}.
			\end{gathered}
		\end{align}
		In view of an argument similar to \eqref{eq:quad-ineq-1} we have with a probability $1-O\pr{n^{-ck}+\exp(-cn)+\exp\pr{-cmL}}$
		\begin{align*}
			&\frac 1n\|\xvt_o- \xv_o\|_2^2
			\leq {2\over n}\pth{{d\over a}+{b^2 \over a^2}}\\
			&
			\leq
			O_{x_{\min},x_{\max}}
			\pth{{k\log n\over nm^2}\sth{{\max(m,n)}
			+{(\cnms)^4(\sigma_z^2+x_{\max}^2 \emax )^2
			\over L}
            +{\frac{\max(m^2,n^2)}{m}}}}.\\
        	& \leq
        	O_{x_{\min},x_{\max}}
        	\pth{{k\log n\over nm^2}\sth{
        	{(\cnms)^4(\sigma_z^2+x_{\max}^2 \emax )^2 \over L}
        	+{\frac{\max(m^2,n^2)}{m}}}}.
		\end{align*}
		This implies, in view of $\frac 1n\|\xvt_o- \xv_o\|_2^2\leq x_{\max}^2$, and using $\|\xvt_o-\xvh_o\|_2\leq \delta_{\rm  net}\leq {x_{\max}\over n^5}$ from \eqref{eq:netdelta}, we get
		\begin{align*}
			\EE\qth{\frac 1n\|\xvh_o- \xv_o\|_2^2}
			&\leq 
			O_{x_{\min},x_{\max}}
			\Biggl({k\log n\over nm^2}\sth{
					{(\cnms)^4(\sigma_z^2+x_{\max}^2 \emax )^2 \over L}
					+{\frac{\max(m^2,n^2)}{m}}}
			\nonumber\\
			&\quad +x_{\max}^2(n^{-ck}+L\exp(-cn)+\exp\pr{-cmL})\Biggr).
		\end{align*}
		
	}
		


		{	
			To summarize, we have in the regime $n\ge 4m$, or in the regime $n<4m$ but $\sigma_z^2\ge m$ that 
			\begin{align}
				R_2^{\dagger}(\Cc, m,n, \sigma_z^2)=O_{x_{\max},x_{\min}}\pr{\frac{\max(\sigma_z^4,n^2)}{nm^2} \frac{k\log n}{L}
				+{\frac{\max(m^2,n^2)k \log n}{nm^3}}}.
			\end{align}
			For the case $n<4m,\sigma_z^2\le m$, we have by monotonicity of the risk in $\sigma_Z$ (\prettyref{lmm:monotonicity in sigma}) that
			\begin{align*}
				R_2^{\dagger}(\Cc, m,n, \sigma_z^2)\le R_2^{\dagger}(\Cc, m,n, m)=O_{x_{\max},x_{\min}}\pr{\frac{\max(m^2,n^2)}{m^2 n} \frac{k\log n}{L}+ {\frac{\max(m^2,n^2)k \log n}{nm^3}}}.
			\end{align*}
			Combining these, we have 
			\begin{align*}
				R_2^{\dagger}(\Cc, m,n, \sigma_z^2)=O_{x_{\max},x_{\min}}\pr{\frac{\max(\sigma_z^4,m^2,n^2)}{m^2 n} \frac{k\log n}{L}+ 
				{\frac{\max(\sigma_z^4,m^2,n^2)k \log n}{nm^3}}}.
			\end{align*}
			
		}
			

			\subsection{Proof of \prettyref{thm:lower_bnd2} in the case $n\ge 4m$}
			
			The proof follows the same approach as in \prettyref{sec:proof-of-lower-bound-undersample} and the only difference here is the argument in bounding the KL divergence and $\beta_r$. We shall point out the difference in below. 
			
			Let $A=A_1=\cdots =A_L$ be an $m\times n$ Gaussian matrix, and let $\mathcal{X}_k, \mathcal{X}^{\rm finite}, \calS_{\rm sep}, r,  N_{\rm div}, \alpha_r$ and $\delta_r$ be chosen as in \prettyref{sec:proof-of-lower-bound-undersample}. In place of \prettyref{lmm: upper bound for KL divergence}, we have (with the same proof)
			
			\begin{lemma}\label{lmm: upper bound for KL divergence fixed}
				Denote $E_{\max}:= \lambda_{\max} (AA^T), E_{\min}:= \lambda_{\min} (AA^T)$. On the event $\Ec_{\rm sing}$, defined in \eqref{eq:esing}, if $\frac{\br{\sigma_z^2+x_{\max}^2\cdot\emax}\cdot \emax}{\pr{\sigma_z^2+x_{\min}^2\cdot E_{\min}}^2}\cdot x_{\max}\delta_r<\frac 14$, we have for all $\bx_i\neq \bx_j\in \calS_{\rm sep}$
				\begin{align*}
					\KL(\P_{\xv_i}\parallel \P_{\xv_j})\le & 2{\left(\sigma_z^2+x_{\max}^2 \emax \right)^2 \over \left(\sigma_z^2+x_{\min}^2 \emin \right)^4} 
					\cdot L
					\nmhs{A(X_i^2-X_j^2)A^\top }^2,
				\end{align*}
				where $X_i$ and $X_j$ are diagonal matrices corresponding to the vectors $\xv_i, \xv_j \in \calS_{\rm sep}$.
			\end{lemma}
			We now apply the upper tail bound of \prettyref{lmm:decoupling-genr} for the case when $L=1$ to find a deterministic upper bound for $L \nmhs{A(X_i^2-X_j^2)A^\top }^2$. We set $\bd_{i,j}:= \bx_i^2-\bx_j^2$, and define $D_{i,j}=\diag(\bd_{i,j})$. Choose
			We choose the following values of $t_{1,i,j}$ and $t_{2,i,j}$ to apply the upper tail bound in \prettyref{lmm:decoupling-genr}
			\begin{align*}
				t_{1,i,j}:=& C_{t_1}\pr{x_{\max}\nm{\xv_i-\xv_j}_{2}\sqrt{\log(mr^2)}+ x_{\max}\nm{\xv_i-\xv_j}_{\infty}\log(mr^2)};\\
				t_{2,i,j}:=& C_{t_2}\log m \pr{x_{\max}^2 \nm{\xv_i-\xv_j}_{4}^2 \sqrt{m(\sqrt m+\sqrt n)^4\log r^2}+x_{\max}^2\nm{\xv_i-\xv_j}_{\infty}^2 (\sqrt n +\sqrt m)^2 \log r^2 },
			\end{align*}
			where $C_{t_1}$ and $C_{t_2}$ are two constants. In view of the above definition, consider the event 
			\begin{equation}\label{eq:dcpl fixed}
				\Ec_{\rm dcpl}:=\bigcap_{1\le i< j\le r}\br{\|A (X_i^2-X_j^2)A^\top \|^2_{\rm HS} 
					< mt_{1,i,j}^2+m(m-1)\|\bd_{i,j}\|_{2}^2 + t_{2,i,j}}.
			\end{equation}
			
			Using the same argument as in \prettyref{sec:The bound for Kullback-Leibler divergence}, we have
			
			\begin{align}
				\begin{split}
					\P(\Ec_{\rm dcpl}^c \cap \widetilde\Ec_{\rm maxsing})	
					\leq  r^2\exp\sth{-\tilde C\pth{\log(mr^2)+(\log m)(\log r)}},
				\end{split}
			\end{align}
			and for sufficiently large $C_{t_1},C_{t_2}$, we have
			\begin{align}\label{eq:m10 fixed}
				\P(\Ec_{\rm dcpl}\cap \widetilde\Ec_{\rm maxsing})
				\ge 1-\frac 1{(rm)^{8}}-2m\exp\pr{-cn}.
			\end{align}
			In view of the above, on the high-probability event $\Ec_{\rm dcpl}\cap \widetilde\Ec_{\rm maxsing}$, we have for each $1\leq i<j\leq r$,
			\begin{align*}
				\nmhs{A_l(X_i^2-X_j^2)A_l^\top }^2 
				\leq & mt_{1,i,j}^2+m(m-1)\|\bd_{i,j}\|_{2}^2 + t_{2,i,j}\\
				& + C\log m\pr{x_{\max}^2 \sqrt{\frac{kn}{N_{\rm div}}}\de_r^2 \sqrt{m(\sqrt m+\sqrt n)^4\log r^2}+\de_r^2 (\sqrt n +\sqrt m)^2 \log r^2 },
			\end{align*}
			Hence, restricting to the event $\Ec_{\rm dcpl}\cap \widetilde\Ec_{\rm maxsing}\cap \Ec_{\rm sing}$, using \prettyref{lmm: upper bound for KL divergence}, we have for $\bar C:=C_{x_{\min},x_{\max}}>0$
			\begin{align}\label{eq:bound-beta fixed}
				&~\beta_r:=\max_{1\le i< j\le r} \KL(\P_i\parallel \P_j)
				\nonumber\\
				&\le 2\frac{\pr{\sigma_z^2+x_{\max}^2\cdot\emax}^2 L}{\pr{\sigma_z^2+x_{\min}^2\cdot \frac{1}{4}(\sqrt{n}-\sqrt{m})^2}^4}\max_{1\le i< j\le r} \nmhs{A(X_i^2-X_j^2)A^\top }^2
				\nonumber\\ 
				&\le {\bar C L \over \max(\sigma_z^4,n^2)}\Bigg({ m  \frac{kn}{N_{\rm div}}\de_r^2\log(mr^2)}+m \de_r^2\log^2(mr^2) +{m(m-1)\frac{kn}{N_{\rm div}}\de_r^2}
				\nonumber \\
				&\qquad \qquad \quad+\log m { \sqrt{\frac{kn}{N_{\rm div}}}\de_r^2 \sqrt{m(\sqrt m+\sqrt n)^4\log r^2}}+ {\de_r^2 (\sqrt n +\sqrt m)^2 (\log m)(\log r^2) }\Bigg)
				\nonumber\\
				&
				\leq {\bar C\delta_r^2m^2nLk\over \max(\sigma_z^4,n^2)N_{\rm div}}
				\Bigg({\log(mr^2)\over m} + {\log^2(mr^2)N_{\rm div}\over mnk}
				+1 +{ \sqrt{\frac{n N_{\rm div}\log r^2}{km^3/(\log m)^2}}}+ {(\log m)(\log r^2)N_{\rm div}\over m^2}\Bigg)
				\nonumber \\
				&\leq
				\Theta_{x_{\min},x_{\max}}(1){m^2nLk\over \max(\sigma_z^4,n^2)N_{\rm div}}\delta_r^2,
			\end{align}
			where the last inequality followed by factoring out $\delta_r^2$ and using the following inequalities that are consequences of \prettyref{lmm:r-bound}, alongside our assumptions $\log m =\Theta(\log n),\log L=O(\log n)$, and there exists $\e\in (0,1/2)$ such that $k \le n^{1-2\e}, \max(\sigma_z^4,m^2,n^2)k\log n \le m^2n^{1-\e}L$. Note that the bound \eqref{eq:bound-beta fixed} is the same as the bound \eqref{eq:bound-beta} for $\beta_r$. 
			
			As a consequence, by \prettyref{lmm:fano} we have for any estimator $\xvh$, the same lower bound
			\begin{equation*}
				\begin{split}
					\max_{1\le i \le r} \E\br{\frac{\|\xvh-\xv_i\|^2}{n}}
					\ge & \frac{\alpha_r^2}{4n}\pr{1-\frac{\beta_r+\log 2}{\log r}}^2 \P\pr{\widetilde\Ec_{\rm maxsing}\cap \Ec_{\rm dcpl} \cap \Ec_{\rm sing}}
					=\Theta\pr{\frac{\alpha_r^2}{n}}\\
					=& \Theta_{x_{\max},x_{\min}}\pr{\frac{\max(\sigma_z^4,n^2)}{m^2 n} \frac{k\log \pr{N_{\rm div}/k}}{L}}
					=\Theta_{\e,x_{\max},x_{\min}}\pr{\frac{\max(\sigma_z^4,n^2)k\log n}{m^2 n L}}.
				\end{split}
			\end{equation*}
			
			\subsection{Proof of \prettyref{thm:lower_bnd2} in the case $n\le 4m$}
			
			The proof follows the same approach as in \prettyref{sec:proof of lower bound oversample} and here we point out the main differences. This section will primarily establish the lower bound for the sub-case $m\geq 4n,\sigma_z^2=0$ given by 
			\begin{equation}\label{eq:m12 fixed}
				R_2^{\dagger}(\Cc_k,m,n,0)
				=\inf_{\bde}\sup_{\xv\in \Cc_k}\E\br{\frac{\nm{\bde(\ly)-\xv_o}_2^2}{n} }
				=  \Omega_{\e,\xmin\xmax}\pr{\frac{k\log n}{nL}},
				\quad 
				m\geq 4n.
			\end{equation}
			Then, the lower bound for a general $\sigma_z^2\geq 0$ and $4n\geq m\ge \frac n4$ follow from \prettyref{lmm:monotonicity in m} and \prettyref{lmm:monotonicity in sigma} by the same monotonicity argument as in \prettyref{sec:proof of lower bound oversample} before \eqref{eq:m19}. Note that, in view of \prettyref{lmm:monotonicity in m}, for any $\frac n4\leq m\leq 4n$, the last display implies
			\begin{equation}\label{eq:m19 fixed}
				R_2^{\dagger}(\Cc_k,m,n,0)\ge
				R_2^{\dagger}(\Cc_k,n/4,n,0)\ge R_2^{\dagger}(\Cc_k,4n,n,0)=\Omega_{\e,x_{\min},x_{\max}}\pr{ \frac{k\log n}{nL}}.
			\end{equation}
			By the same argument as in the paragraph before \eqref{eq:large noise case lower bound},  in the case $4m\ge n$ and $ \sigma_z^2 \ge m$ we have the lower bound
			\begin{equation}\label{eq:large noise case lower bound fixed}
				R_2^{\dagger}(\Cc_k,m,n,\sigma_z)=\Omega_{x_{\max},x_{\min}}\pr{\frac{\sigma_z^4}{m^2 n} \cdot \frac{k\log n}{L}}, \text{ whenever } \sigma_z^2\ge m.
			\end{equation}
			To achieve a lower bound for the sub-case $m\geq \frac n4, \sigma_z^2 \le m$, we first use that the miminax error is non-decreasing function in $\sigma_z$ (\prettyref{lmm:monotonicity in sigma}) to get $R_2^{\dagger}(\Cc_k,m,n,0)\leq R_2^{\dagger}(\Cc_k,m,n,\sigma_z)$. Then, to achieve a lower bound to $R_2^{\dagger}(\Cc_k,m,n,0)$ we combine the lower bounds in \eqref{eq:m12 fixed} and \eqref{eq:m19 fixed} to get
			\begin{align}\label{eq:lower bound for sigma equal to zero fixed}
				R_2(\Cc_k,m,n,0)=\Omega_{x_{\max},x_{\min}}\pr{\frac{m^2}{m^2 n} \cdot \frac{k\log n}{L}}=\Omega_{x_{\max},x_{\min}}\pr{\frac{k\log n}{nL}}.
			\end{align}
			Then, for $\sigma_z^2\le m$, by \eqref{eq:section on the decoupling in the upper bound} we get for a constant $C_1$ depending on $\xmin,\xmax$
			\begin{equation*}
				{C_1\frac{k\log n}{nL}}\leq R_2^{\dagger}(\Cc_k,m,n,0)\leq 
				R_2^{\dagger}(\Cc_k,m,n,\sigma_z).
			\end{equation*}
			Combining \eqref{eq:large noise case lower bound fixed} and \eqref{eq:lower bound for sigma equal to zero fixed} for the case $n\le 4m$ yields the desired minimax lower bound
			\begin{align*}
				R_2^{\dagger}(\Cc_k,m,n,\sigma_z)=\Omega_{x_{\max},x_{\min}}\pr{\frac{\max(\sigma_z^4,m^2)}{m^2 n} \cdot \frac{k\log n}{L}},
				\quad m\geq \frac n4, \sigma_z\geq 0.
			\end{align*}
			
			Now it remains to establish \eqref{eq:m12 fixed}. We will show this with a sufficient statistic as in \prettyref{sec:proof of lower bound oversample}. Define $\bA=\diag(A,\dots,A)\in \reals^{mL\times nL}$ and note that $\bA^\top \bA=\diag(A^\top A,\dots,A^\top A)$. Throughout the section we analyze the expected loss on the high probability event $\esingp$ defined in \eqref{eq: definition of esing prime}, where $A=A_1=\cdots=A_L$ and $A^{\top}A$ is invertible. Then with the same proof as \prettyref{prop:sufficient}, we have
			\begin{proposition}\label{prop:sufficient fixed}
				Consider the case $\sigma_z=0$ and that the event $\esingp$ holds. Then $\bT_{\bA}(\ly)=(\bA^\top  \bA)^{-1}\bA^\top \ly$ is a sufficient statistic for the parameter $\xv_o$.
			\end{proposition}
			By Rao-Blackwell theorem (\prettyref{thm:rao}), we have,
			\begin{equation*}
				\EE\qth{\nm{\bde(\ly)-\xv_o}_2^2|\esingp}\ge \EE\qth{\nm{g(\bT_{\bA}(\ly))-\xv_o}_2^2|\esingp}.
			\end{equation*}
			
			Therefore it suffices to prove the following lemma
			\begin{lemma}\label{lmm:lower bound for the improved estimator fixed}
				Consider the model \eqref{eq:multilook speckle noise model} with $\sigma_z=0,m\ge 4n$. Then, there exists a constant $C\geq 0$, we have
				$$
				\inf_{g}\sup_{\xv\in \Cc_k}\E\br{\nm{g(\bT_{\bA}(\ly))-\xv_o}_2^2|\esingp} \geq C{\frac{k\log n}{L}}.
				$$
			\end{lemma}
			The proof of \prettyref{lmm:lower bound for the improved estimator fixed} is the same as that of \prettyref{lmm:lower bound for the improved estimator}. Finally, it follows from it follows from \prettyref{thm:rao} and \prettyref{lmm:lower bound for the improved estimator fixed} that
			\begin{equation*}
				\inf_{\bde}\sup_{\xv\in \Cc_k}\E\br{\frac{\nm{\bde(\ly)-\xv_o}_2^2}{n} }
				\ge \inf_{\bde}\sup_{\xv\in \Cc_k}\E\br{\left.\frac{\nm{g(\bT_{\bA}(\ly))-\xv_o}_2^2}{n} \right \vert \esingp}\P\pr{\Ec_{\rm sing}'}
				=  \Omega\pr{\frac{k\log n}{nL}}.
			\end{equation*}

			\section{Proof of \prettyref{thm:ksparse}}\label{app:sparse case}
			
			\subsection{Proof of the lower bound}
			For $\mathcal{S}_k^{\rm bdd}:=\cb{\xv\in \R^n: \nm{\xv}_o\le k, x_i=0 \text{ or } 0<x_{\min}\le x_i \le x_{\max}},$
			and any $S\subset [n]$ and $|S|\le k$,
			\begin{align*}
				\Cc_S:=\cb{\xv_S\in\R^{|S|}: \xmin \le x_i \le \xmax, i\in S}=[\xmin,\xmax]^{|S|}.
			\end{align*}
			Also, for a fixed $S \subset [n]$ with $|S| \leq k$, define 
			\[
			\mathcal{S}_{k, S}^{\rm bdd} = \{\xv \in \mathbb{R}^n : x_i = 0, i \notin S, \ {\rm and} \ x_i \in [x_{\min} , x_{\max} ] \ {\rm for} \ i \in S  \}. 
			\]
			In this section, we aim to prove the lower bound, i.e we aim to prove that     \begin{align}
				\begin{split}
					& R_2(\mathcal{S}_k^{\rm bdd}, m,n,k,L, \sigma_z)
					= \inf_{\xvh \in \R^n} \sup_{\xv_o\in \mathcal{S}_k^{\rm bdd}} \E\br{{\nm{\xvh-\xv_{o}}_2^2 \over n}}=\Omega_{\xmax,\xmin}\pr{k \over nL}.
				\end{split}
			\end{align}
			
			First note that by the monotonicity result proved in \prettyref{lmm:monotonicity in sigma}, we have $R_2(\mathcal{S}_k^{\rm bdd}, m,n,k,L, \sigma_z)\ge R_2(\mathcal{S}_k^{\rm bdd}, m,n,k,L, 0)$. Therefore it suffices to establish
			\begin{align}
				R_2(\mathcal{S}_k^{\rm bdd}, m,n,k,L, 0)=\Omega_{\xmax,\xmin}\pr{k \over nL}.
			\end{align}
			
			For any $S\subset [n]$ and $|S|\le k$, since $\mathcal{S}_{k, S}^{\rm bdd} \subset \mathcal{S}_{k}^{\rm bdd} $
			\begin{align*}
				\inf_{\xvh \in \R^n} \sup_{\xv_o\in \mathcal{S}_k^{\rm bdd}} \E\br{{\nm{\xvh-\xv_{o}}_2^2 \over n}} \geq \inf_{\xvh \in \R^n} \sup_{\xv_o\in \mathcal{S}_{k, S}^{\rm bdd}} \E\br{{\nm{\xvh-\xv_{o}}_2^2 \over n}}
				\inf_{\xvh_S \in \R^{|S|}} \sup_{\xv_{o,S}\in \Cc_S} \E\br{{\nm{\xvh_S -\xv_{o,S}}_2^2 \over n}}.
			\end{align*}
			To obtain the last equality we have noted that since we know the exact location of non-zero elements of $\xv_o$ for any $\xv_o \in \mathcal{S}_{k,S}^{\rm bdd}$, we have set the value of $\hat{\xv}$ to zero at those locations, and have reduced the problem to estimating the nonzero elements of $\xv_o$.  
			We now claim that, in particular, if $|S|=k$, we have
			\begin{align}\label{reduced sparse model lower bound}
				R_2(\Cc_S, m,n,k,L, 0)=\Omega_{x_{\max},x_{\min}}\pr{\frac{k}{n L}}.
			\end{align}
			Indeed, in this case our model reduces to 
			$
				\yv_l=A_{l,S} X_{o,S} \wv_{l,S}, l=1,2,...,L,
			$
			where $A_{l,S}$ is the $m\times k$ matrix of whose columns are those of $A_l$ with indices in $S$, $X_{o,S}$ is the $k\times k$ diagonal matrix in which all the diagonal locations with indices in $S$- are nonzero, and $\wv_{l,S}$ is the $k$-dimensional Gaussian vector. Since $m\gg k$, we know that $\widetilde \yv_l := (A_{l,S}^{\top} A_{l,S})^{-1} A_{l,S}^{\top}\yv_l$ is a sufficient statistic (see \prettyref{prop:sufficient}) and by using Rao-Blackwell theorem (\prettyref{thm:rao}), it suffices to obtain a lower bound for the minimax risk of estimation $\xv_{S}$ from $\widetilde \yv_1, \widetilde \yv_2, \ldots, \widetilde \yv_L$:
			\begin{align}
				\widetilde \yv_l:=X_{o,S} \wv_{l,S},\quad l=1,2,...,L.
			\end{align}
			We shall prove this by \prettyref{lmm:fano}. Note that the joint distribution of $\widetilde \yv_1, \widetilde \yv_2, \ldots, \widetilde \yv_L$ is given by:
			\begin{equation*}
				\P_{\xv} \sim \otimes_{l=1}^L N(\boldsymbol{0},X_{o,S}^2).
			\end{equation*}
			For $\xv,\xv'\in \R^n$, we define the pseudo-metric,
			\begin{equation*}
				d(\theta(\P_{\xv}),\theta(\P_{\xv'}))=d(\xv,\xv'):=\|\xv-\xv'\|_2.
			\end{equation*}
			
			Next, we shall construct an $\alpha_r$-separated set $\{\xv_1,...,\xv_r\}$ in $[\xmin,\xmax]^{|S|}=[\xmin,\xmax]^{k}$. To construct this subset, we use the following steps. Consider the following definition.
			\begin{definition}\label{def:deltapacking}
				
				A $\alpha$-\emph{separating} subset of $S$ is a finite or countable collection $\cb{\xv_i}$ of points of $\text{X}$ satisfying ${\rm dist}(\xv_i,\xv_j)\ge \alpha$ for any $i\ne j$. We call the largest possible cardinality among all $\alpha$-separated subsets of $S$ the $\alpha$-\emph{packing number} of $S$, denoted $P_{\alpha}(S)$. In other words,
				\begin{equation*}
					P_{\delta}(S):=\sup\cb{n:\text{there exists a $\alpha$-separated subset of $S$ of cardinality $n$}}.
				\end{equation*}
			\end{definition}
			Using the  $\delta$-packing defined in \prettyref{def:deltapacking} we have:

			\begin{lemma}\label{lmm:Hamming cube packing number}
				For any $p< \frac{k}{2}$, the $k$-dimensional Hamming cube $\cb{0,1}^k$ has $p$-packing number at least $\frac{2^k}{\sum_{i=0}^{p}{k \choose i}}$ with respect to the Hamming distance.
			\end{lemma}
			
			\begin{proof}
				Let $\mathcal{S}$ be a $p$-packing subset of $\cb{0,1}^k$. By definition we must have
				\begin{align}\label{eq: containment hamming}
					\cb{0,1}^k \subset \cup_{\xv\in \mathcal{S}} B(\xv, p),
				\end{align}
				where $B(\xv, p)$ is the ball centered at $\xv$ with radius $p$ with respect to the Hamming distance.  In other words, $B(\xv, p)$ is the collection of points in $\cb{0,1}^k$ whose at most $p$ coordinates are different from those of $\xv$.
				A direct counting gives $|B(\xv, p)|\le \sum_{i=0}^{p}{k \choose i}$.
				Taking cardinality of both sides of \eqref{eq: containment hamming} yields
				\begin{align*}
					2^k \le \sum_{i=0}^{p}{k \choose i} |\mathcal S|.
				\end{align*}
				This completes the proof.
			\end{proof}
			To obtain a simpler lower bound for the $p$-packing that can be used in our arguments, we assume that $p = k/7$ and establish the following lemma:
			\begin{lemma}\label{lmm:lower bound for packing number}
				We have $\frac{2^k}{\sum_{i=0}^{\fl{k/7}}{k \choose i}}\ge a^k$ for some $a>1$.
			\end{lemma}
			\begin{proof}
				Indeed, for $p:=\fl{k/7} \le \frac{k}{2}$, since ${k \choose i-1}/{k \choose i}=\frac{i}{k-i+1}\le \frac{1}{2}$, for $i\le p \le k/2$, we have 
				\begin{equation*}
					\sum_{i=0}^{p} {k \choose i}\le \pr{1+\frac{1}{2}+\cdots \frac{1}{2^{p-1}}}{k \choose p} \le  2{k \choose p} \le 2 \pr{ek \over p}^p.
				\end{equation*}
				So $\log \frac{2^k}{\sum_{i=0}^{p}{k \choose i}}\ge k\log 2 - p \log \pr{ek\over p} - \log 2 \gg k \cdot b$ for some $b>0$. Then the result follows for $a=e^b$.
			\end{proof}
			
			Let $\mathcal{H}\subset \cb{0,1}^k$ be an $\fl{k/7}$-packing set with respect to the Hamming distance from \prettyref{lmm:Hamming cube packing number}. By \prettyref{lmm:lower bound for packing number}, $|\mathcal{H}|\ge a^k$.
			
			
			Now we aim to use $\mathcal{H}$ to create a subset of $\mathcal{C}_S$.
			Define $\bar{x} := \frac{x_{\min}+ x_{\max}}{2}. $ We create the subset by scaling and translation of $\mathcal{H}$ in the following way:
			\begin{align}\label{eq:alpha separation set in the sparse case}
				c\frac{\xmax-\xmin}{\sqrt{L}} \cdot \mathcal{H}+\overline{\xv} =\cb{\xv_1,...,\xv_r}\subset B_{\ell_2}\pr{\overline{\xv},c\sqrt{\frac{k}{L}}\cdot \pr{\xmax-\xmin}} \cap [\xmin,\xmax]^k.
			\end{align}
			Note that this set
			forms a $c\frac{\xmax-\xmin}{\sqrt{L}}\sqrt{\fl{k/7}}$ -separated set for $[\xmin,\xmax]^{k}$ with respect to the $\ell_2$-distance,
			with cardinality $r\ge a^k$ for some $a>1$. To use \prettyref{lmm:fano} set $\alpha_r:=c\frac{\xmax-\xmin}{\sqrt{L}}\sqrt{\fl{k/7}}$.
			It follows from \prettyref{lmm: upper bound for KL divergence} (with $A_l$ replaced by $I_n$) that
			\begin{align*}
				\begin{split}
					& \beta_r:=\max_{1\le i<j\le r}\KL(\P_i\parallel \P_j)
					\le  L\frac{\xmax^4}{\xmin^8}
					\max_{1\le i<j\le r}\nmhs{X_i^2-X_j^2}^2 \\
					\le &  \frac{4\xmax^6}{\xmin^8} L\max_{1\le i<j\le r} \nm{\xv_i-\xv_j}_2^2 
					\le \frac{4\xmax^6}{7\xmin^8} (\xmax-\xmin)^2 c^2 k \le \frac{1}{10} \log r,
				\end{split}
			\end{align*}
			for sufficiently small $c$ from \eqref{eq:alpha separation set in the sparse case}. 
			Hence, by \prettyref{lmm:fano},
			\begin{align*}
				\inf_{\xvh_S \in \R^{|S|}} \sup_{\xv_{o,S}\in \Cc_S} \E\br{{\nm{\xvh_S -\xv_{o,S}}_2^2 \over n}}
				\ge  \frac{\alpha_r^2}{4n}\pr{1-\frac{\beta_r+\log 2}{\log r}}^2 \P\pr{\Ec'_{\rm sing}}
				=\Theta\pr{\frac{\alpha_r^2}{n}}=\Theta_{x_{\min},x_{\max}}\pr{ \frac{k}{nL}}.
			\end{align*}

			\subsection{Proof of the upper bound}
			In this section we aim to prove the upper bound. In other words, we aim to prove that

			\begin{align*}
				R_2(\mathcal{S}_k^{\rm bdd}, m,n,k,L, \sigma_z)
				= \inf_{\xvh \in \R^n} \sup_{\xv_o\in \mathcal{S}_k^{\rm bdd}} \E\br{{\nm{\xvh-\xv_{o}}_2^2 \over n}}=O_{\xmax,\xmin}\pr{{k \over nL} + \frac{\sigma_z^2 k\log(n/k)}{mn}}
			\end{align*}
			Consider the model
			\begin{equation}\label{eq:appendixD}
				\yv_l = A_l X_o \wv_l + \zv_l, \text{ 
					for } l=1,\ldots,L.
			\end{equation}
			Define, $\uv_l:=X_o \wv_l \sim N(0,X_o^2)$. Observe that $\uv_l$ is unbounded, but still $k$-sparse. If we think of $u_l$ as our new unkown data, the model \eqref{eq:appendixD} reduces to $L$ copies of classical sparse linear regression models
			\[
			\yv_l = A_l \uv_l + \zv_l, \quad l = 1, \ldots, L,
			\]
			where $\uv_l\in \mathcal{S}_k:=\{\xv\in \R^k:\nm{\xv}_0\le k \}$.
			Our strategy for obtaining the upper bound is to first estimate each $\bu_\ell$ separately from $\yv_l$. This is a standard problem in sparse linear regression. We use one of the classic results for obtaining an upper bound for the minimax risk of sparse linear regression model:

			\begin{theorem}\label{thm:Verzelen12}\cite[Equation (3.9) and Proposition 6.4]{Verzelen12}
				Let $A$ be an $m\times n$ Gaussian matrix and $\zv\sim N(0,\sigma_z^2 I_n)$. Consider the sparse linear regression model $\yv=A\uv+\zv$ where the unknown signal $\uv \in \mathcal{S}_k$. Then for $k\log(en/k)\le m$, we have minimax risk estimate 
				\begin{align*}
					\inf_{\uvh} \sup_{\uv\in \mathcal{S}_k} \E\br{\nm{\uvh-\uv}_2^2}=\Theta\pr{\frac{\sigma_z^2 k\log(en/k)}{m}}.
				\end{align*}
				
			\end{theorem}

			For each $1\le l \le L$, let $\widehat{\uv}_l$ be the minimax estimator from \prettyref{thm:Verzelen12}, and let $\xvh^2:=\frac{1}{L}\sum_{l=1}^L \widehat{\uv}_l^2$ be our estimator for the unknown signal $\xv_o\in \mathcal{S}_k^{\rm bdd}$. It follows that
			\begin{align*}
				\begin{split}
					\E\br{\nm{\xvh^2-\xv_o^2}_2^2}
					=&\E\br{\nm{\xvh^2-\frac{1}{L}\sum_{l=1}^L \uv_l^2+\frac{1}{L}\sum_{l=1}^L \uv_l^2-\xv_o^2}_2^2}\\
					\overset{(a)}{\le} &2\E\br{\nm{\xvh^2-\frac{1}{L}\sum_{l=1}^L \uv_l^2}_2^2+\nm{\frac{1}{L}\sum_{l=1}^L \uv_l^2-\xv_o^2}_2^2}\\
					=& {2\over L^2}\E\br{\nm{\sum_{l=1}^L\pr{\uvh_l^2-\uv_l^2}}_2^2}+{2\over L^2}\E\br{\nm{\sum_{l=1}^L\pr{\uv_l^2-\xv_o^2}}_2^2}.
				\end{split}
			\end{align*}
			where for (a) we used $(a+b)^2\le 2(a^2+b^2)$ in the vector form.
			
			Since $(\uv_l^2-\xv_o^2)$'s are $k$-sparse, independent and mean zero, for the second term we have 
			\begin{align*}
				{2\over L^2}\E\br{\nm{\sum_{l=1}^L\pr{\uv_l^2-\xv_o^2}}_2^2}= {2\over L^2}\sum_{l=1}^L\E\br{\nm{\uv_l^2-\xv_o^2}_2^2}\le \frac{8\xmax^2 \E\br{(\zeta^2-1)^2} k}{L}.
			\end{align*}
			where $\zeta\sim N(0,1)$. 
			To treat the first term, we notice that by \prettyref{thm:Verzelen12}
			\begin{align*}
				{2\over L^2}\E\br{\nm{\sum_{l=1}^L\pr{\uvh_l^2-\uv_l^2}}_2^2}\overset{(a)}{\le} {2\over L^2}\cdot L\sum_{l=1}^L \E\br{\nm{\uvh_l^2-\uv_l^2}_2^2}\le \frac{2\sigma_z^2 k\log(en/k)}{m}.
			\end{align*}
			
			where for (a) we have used the elementary inequality $(a_1+\cdots +a_L)^2\le L(a_1^2+\cdots+a_L^2)$.

\end{document}